\documentclass[reqno,12pt]{article}

\def\today{${\scriptscriptstyle\number\day-\number\month-\number\year}$}
\usepackage{graphicx}
\usepackage{fancyhdr}
\usepackage{amsmath}
\usepackage{amsthm}
\usepackage{amsfonts}
\theoremstyle{plain}
\newtheorem{theorem}{Theorem}[section]

\newtheorem{lemma}[theorem]{Lemma}
\newtheorem{corollary}[theorem]{Corollary}
\newtheorem{definition}[theorem]{Definition}

\newtheorem{remark}[theorem]{Remark}

\def\address#1{{\center{#1}}}
\date{}
\def\m@th{\mathsurround=0pt}
\def\eqal#1{\null\,\vcenter{\openup\jot\m@th
 \ialign{\strut\hfil$\displaystyle{##}$&&$\displaystyle{{}##}$\hfil
 \crcr#1\crcr}}\,}
\def\matrix#1{\null\,\vcenter{\normalbaselines\m@th
 \ialign{\hfil$##$\hfil&&\quad\hfil$##$\hfil\crcr
 \mathstrut\crcr\noalign{\kern-\baselineskip}
 #1\crcr\mathstrut\crcr\noalign{\kern-\baselineskip}}}\,}
\def\lc#1{\hbox to .89\hsize{\qquad$\displaystyle{#1}\hfill$}}

\def\ro{\varrho}
\def\D{{\mathbb D}}

\def\N{{\mathbb N}}
\def\R{{\mathbb R}}
\def\T{{\mathbb T}}
\def\divv{{\rm div}\,}

\def\inf{{\rm inf}}
\def\supp{{\rm supp\,}}
\def\diam{{\rm diam\,}}
\def\const{{\rm const}}
\def\mize{\begin{itemize}\parskip-4pt}

\def\bye{\end{document}}
\numberwithin{equation}{section}

\title{Local solutions for nonhomogeneous Navier-Stokes equations with large flux}
\author{Joanna Renc\l awowicz$^{1)}$\quad Wojciech M. Zaj\c{a}czkowski$^{1), 2)}$}

\begin{document}
\input amssym.def
\input amssym.tex
\maketitle
\thispagestyle{fancy}

\address{$^{1)}$\ Institute of Mathematics, Polish Academy of Sciences,\\
\'Sniadeckich 8, 00-656 Warsaw, Poland\\ e-mail: jr@impan.pl, {\it corresponding author,}\\ $^{2)}$\ Institute of Mathematics and Cryptology, Cybernetics Faculty, \\
Military University of Technology,\\
S. Kaliskiego 2, 00-908 Warsaw, Poland\\
e-mail: wz@impan.pl\\}

\let\thefootnote\relax

\begin{abstract}
The local existence of solutions to nonhomogeneous Navier-Stokes equations in cylindrical domains with arbitrary large flux is demonstrated. The existence is proved by the method of successive approximations. To show the existence with the lowest possible regularity the special Besov spaces called the Sobolev-Slobodetskii spaces are used. The inflow and outflow are prescribed on the parts of the boundary which are perpendicular to the $x_3$-axis. Since the inflow and outflow are positive the crucial point of this paper is to verify that $x_3$-coordinate of velocity is also positive.

Finally, we conclude the local existence such that the velocity belongs to $W_{\sigma} ^{2+{s},1+{s}/2}(\Omega^t)$, the gradient of pressure to $W_{\sigma} ^{{s},{s}/2}(\Omega^t)$ and the density to $W_{r,\infty}^{1,1}(\Omega^t)$, where ${s}\in(0,1)$, ${\sigma} >3/{s}$, $r>5/{s}$, $r>{\sigma} $.
\end{abstract}

{\bf Abbreviated title:} Local solutions for nonhomogeneous NS.

{\bf Statement relating to the ethics and integrity policies:}

\noindent On behalf of all authors, the corresponding author states that the manuscript is not published elsewhere and there is no conflict of interest.

\noindent The datasets generated during and/or analysed during the current study are available from the corresponding author on reasonable request.

{\bf Funding statement:} There are no funders to report for this submission.

\section{Introduction}\label{s1}

In the paper, we analyze the local existence of solutions to nonhomogeneous Navier-Stokes equations in cylindrical domains with large inflow and outflow. With "nonhomogeneous" we mean a density dependent system.

Let $\Omega\subset\R^3$ be a cylindrical domain parallel to the $x_3$-axis of Cartesian coordinates which is located inside $\Omega$. The boundary of $\Omega$ denoted by $S$ is composed of two parts, $S_1$ and $S_2$, where $S_1$ is parallel to the $x_3$-axis and $S_2$ is perpendicular to it. 

Let the real number $a>0$ be given. Then $S_2(-a)$ meats the $x_3$-axis at $x_3=-a$ and $S_2(a)$ at $x_3=a$. In $\Omega^T=\Omega\times(0,T)$, where $T>0$ is given, we consider the following initial-boundary value problem
\begin{equation}\eqal{
&\varrho v_t+\varrho v\cdot\nabla v-\divv\T(v,p)=\varrho f\quad &{\rm in}\ \ \Omega^T,\cr
&\divv v=0\quad &{\rm in}\ \ \Omega^T,\cr
&\varrho_t+v\cdot\nabla\varrho=0\quad &{\rm in}\ \ \Omega^T,\cr
&v\cdot\bar n=0\quad &{\rm on}\ \ S_1^T=S_1\times(0,T),\cr
&\nu\bar n\cdot\D(v)\cdot\bar\tau_\alpha+\gamma v\cdot\bar\tau_\alpha=0,\ \ \alpha=1,2\quad &{\rm on}\ \ S_1^T,\cr
&v\cdot\bar n=d\quad &{\rm on}\ \ S_2^T=S_2\times(0,T),\cr
&\varrho=\varrho_1\quad &{\rm on}\ \ S_2^T(-a),\cr
&\bar n\cdot\D(v)\cdot\bar\tau_\alpha=0,\ \ \alpha=1,2\quad &{\rm on}\ \ S_2^T,\cr
& \int_{\Omega} p\, dx = 0, \cr
&v|_{t=0}=v(0),\ \ \varrho|_{t=0}=\varrho(0)\quad &{\rm in}\ \ \Omega,\cr}
\label{1.1}
\end{equation}
where $v$ is the velocity of the fluid with
$v(x,t)=(v_1(x,t),v_2(x,t),v_3(x,t))\in\R^3,$
$p=p(x,t)\in\R^1$ denotes the pressure, $\ro=\ro(x,t)\in\R^1$ -the density,
$f=f(x,t)=(f_1(x,t),f_2(x,t),f_3(x,t))\in\R^3$ -- the external
force field, $x=(x_1, x_2, x_3)$ are the Cartesian coordinates.

By $\nu>0$ we denote the constant viscosity coefficient, $\gamma>0$ is the slip coefficient, $\bar n$ is the unit outward vector normal to $S$, $\bar\tau_\alpha$, $\alpha=1,2$, are vectors tangent to $S$, $I$ is the unit matrix, the stress tensor has the form
$$
\T(v,p)=\nu\D(v)-pI,
$$and $\D(v)$ is the dilatation tensor
$$
\D(v)=\{v_{i,x_j}+v_{j,x_i}\}_{i,j=1,2,3}.
$$
We define the inflow and outflow with $d=(d_1,d_2), d_i\ge 0,$ $i=1,2,$ then the boundary condition $(\ref{1.1})_6$ reads
\begin{equation}\eqal{
&d_1=-v\cdot\bar n|_{S_2(-a)},\cr
&d_2=v\cdot\bar n|_{S_2(a)}\cr}
\label{1.2}
\end{equation}
and $\bar n$ is the unit outward vector normal to $S_2$. Since incompressible motions are considered the following compatibility condition holds:
\begin{equation}
\intop_{S_2(-a)}d_1dS_2=\intop_{S_2(a)}d_2dS_2.
\label{1.3}
\end{equation}
We also set $a_1=-a$, $a_2=a$.

Using Cartesian coordinates and assuming that $\varphi_0(x_1,x_2)=c_0$ is a sufficiently smooth closed curve in the plane $x_3=\const\in(-a,a)$, passing around the $x_3$-axis we have
\begin{equation}\eqal{
&\Omega=\{x\in\R^3\colon\varphi_0(x_1,x_2)<c_0,-a<x_3<a\},\cr
&S_1=\{x\in\R^3\colon\varphi_0(x_1,x_2)=c_0,-a<x_3<a\},\cr
&S_2(-a)=\{x\in\R^3\colon\varphi_0(x_1,x_2)<c_0,x_3=-a\},\cr
&S_2(a)=\{x\in\R^3\colon\varphi_0(x_1,x_2)<c_0,x_3=a\}.\cr}
\label{1.4}
\end{equation}

\begin{theorem}\label{t1.1}
Assume
\mize
\item[1.] Parameters ${s},{\sigma},r$ satisfy ${s}\in(0,1)$, $5/r<{s}$, ${3/{\sigma}}<{s}$, $r>{\sigma}$.
\item[2.] Data functions are such that for the density $\varrho_1\in W_r^{1,1}(S_2^t(-a))$, the initial density $\varrho(0)\in W_r^1(\Omega)$, the initial velocity $v(0)\in W_{\sigma}^{2+{s}-2/{\sigma}}(\Omega),$ the inflow $d_1$ and the outflow $d_2$ belong to $W_{\sigma}^{2+{s}-1/{\sigma} ,1+{s}/2}(S_2^t)$, the external force $f\in W_{\sigma}^{{s},{s}/2}(\Omega^t)$.
\item[3.] There exist positive constants $\bar d_0$, $d_0$, $\bar d_0>d_0$, $d_\infty$ and $b_0$, $b_1$ $\bar b_0$, $\bar b_1$ such that $\bar d_0\ge v_3(0)\ge d_0$, $d_i\ge d_\infty$, $i=1,2$, $\bar b_0\ge\varrho(0)\ge b_0$, $\bar b_1\ge\varrho_1\ge b_1$.
\item[4.] Quantities
$$\eqal{
&\bar d_1=|d_1|_{\infty,S_2^t(-a)}(|\varrho_{1,x}|_{r,S_2^t(-a)}+ |\varrho_{1,t}|_{r,S_2^t(-a)})\cr
&\quad+\|\varrho_x(0)\|_{L_r(\Omega)}(1+\|v(0)\|_{W_r^1(\Omega)}),\cr
&\bar d_2=\|f\|_{W_{\sigma}^{{s},{s}/2}(\Omega^t)}+\sum_{i=1}^2 \|d_i\|_{W_{\sigma}^{2+{s}-1/{\sigma},1+{s}/2}(S_2^t(a_i))}+ \|\varrho_1\|_{W_r^{1,1}(S_2^t(-a))}\cr
&\quad+\|v(0)\|_{W_{\sigma}^{2+{s}-2/{\sigma}}(\Omega)}\cr}
$$
are finite.
\end{itemize}
Then there exists a local solution $(v,p,\ro)$ to the nonhomogeneous Navier-Stokes problem (\ref{1.1}) such that
$$
v\in W_{\sigma}^{2+{s},1+{s}/2}(\Omega^t),\quad \nabla p\in W_{\sigma}^{{s},{s}/2}(\Omega^t),\quad \varrho\in W_{r,\infty}^{1,1}(\Omega^t)
$$
with $$\|\ro\|_{W_{r,\infty}^{1,1}(\Omega^t)}=\|\ro\|_{L_\infty(0,t;L_r(\Omega))}+ \|\ro_x\|_{L_\infty(0,t;L_r(\Omega))} \\ +\|\ro_t\|_{L_\infty(0,t;L_r(\Omega))}$$
and 
$$
\|\ro_1\|_{W_r^{1,1}(S_2^t(-a))}=\|\ro_1\|_{L_r(S_2^t(-a))}+\|\ro_{1,x}\|_{L_r(S_2^t(-a))}+\|\ro_{1,t}\|_{L_r(S_2^t(-a))}.
$$
Moreover, the density remains bounded
$$
\varrho_*\equiv\frac{b_1b_0}{b_1+b_0}\le\varrho(x,t)\le\bar b_0+\bar b_1\equiv\varrho^*
$$
and the velocity and the pressure satisfy
$$
\|v\|_{W_{\sigma}^{2+{s},1+{s}/2}(\Omega^t)}+\|\nabla p\|_{W_{\sigma}^{{s},{s}/2}(\Omega^t)}\le\phi({\rm data}),
$$
where data are described by assumptions 1, 2, 3.\\
Finally, the $x_3$-coordinate of velocity is positive
$$
v_3\ge d_*=d_*(d_0,\bar d_0,\bar d_1,\bar d_2,d_\infty, \|f_3\|_{L_1(0,t;L_\infty(\Omega))})>0.
$$
\end{theorem}
\begin{remark}
Note that $s\in (0,1),$ $\sigma >3$, $r>5.$
\end{remark}
The paper is organized in the following way. In Section \ref{s2}, a simplified notation, a partition of unity and the Besov and Sobolev-Slobodetskii spaces are introduced. Moreover, we prove a lower and an upper bounds for the density. A method of successive approximations, that we apply to prove the existence of solutions in Section~\ref{s9}, is defined in Section \ref{s3}. The crucial point of this technique is the assumption that the third component of velocity on the previous $n$-step is positive. Hence it is shown in the proof of Lemma \ref{l9.1} that the same lower bound for the third component of velocity holds for the $n+1$-step of the method of successive approximations. This means that at each step of the method of successive approximations the third component of velocity is bounded from below by the same positive constant. In Section \ref{s4} we find an estimate for $\varrho_n$ in terms of some norms of velocity. In Sections \ref{s5} and \ref{s6} we find an estimate for $v_{n+1}$ and $p_{n+1}$ in terms of norms of $v_n$, $p_n$ and $\varrho_n$. The relation is precisely described by Lemma \ref{l7.2} in Section \ref{s6}. Then we can show both boundedness and convergence of the sequence and this implies the existence. Finally, the existence of local solutions to the main problem (\ref{1.1}) is proved in Section~\ref{s9}.

The theory of Besov and Sobolev-Slobodetskii spaces can be found in \cite{Ad, BIN, G, N, Tr1}.

Nonhomogeneous Navier-Stokes equations problem (i.e. the density dependent Navier-Stokes system) without flux has been considered by some authors. We mention \cite{AKM}, where the existence of weak solutions was proved and strong solutions were established with $v\in W_2^{2,1}$, $\nabla p\in L_2$ and bounded positive density $\varrho$ for small times in $\R^3$ and large times in $\R^2$. In \cite{LS}, Ladyzhenskaya and Solonnikov have obtained existence results for $v\in W_q^{2,1}$, $\nabla p\in L_q$, $q>n$ and $\varrho\in C^1$, for small times with arbitrary $v_0$ and $f$ and for any given time interval with sufficiently small $v_0$ and $f$. The problem was analyzed in a bounded domain in $\R^n$ with boundary $S\in C^2$ and $v|_{S^T}=0$. In exterior domains, Padula in \cite{P} proved the existence of unique local strong solutions for the initial density $\varrho_0\in L_1(\Omega)\cap W_\infty^1(\Omega)$ satisfying the additional property:
$$
\intop_{\Omega'}\varrho_0dx>0
$$
for any $\Omega'\subset\Omega$ with positive measure.

In \cite{D}, the local existence in bounded domains in $\R^n$, $n\ge 2$ with $C^{2+\varepsilon}$ boundary was obtained for $\varrho_0\in W_q^1$, $q>n$ and small $v_0$ for $n\ge 3$. The assumption on lower bound for the density was relaxed in \cite{Ge}, where $v_0\in H^s(\R^3)$, $s\in(2,5/2)$ and the reciprocal of the density belongs to some Sobolev space embeded in BMO. In \cite{CK}, authors could remove the requirement of non-vanishing of the density by adding the compatibility condition
$$
v\Delta v_0-\nabla p_0=\varrho_0^{1/2}g
$$
for some $(p_0,g)\in H^1\times L_2$ and assuming $v_0\in H^2$, $\varrho_0\in H^1\cap L_\infty$.

In \cite{Z3}, the author considered the equations in a bounded cylinder under boundary slip conditions. Assuming that the derivatives of initial density, initial velocity, external force with respect to the third co-ordinate (see notation (\ref{1.4})) are sufficiently small in some norms, the existence of large time regular solutions in Sobolev-Slobodetskii spaces has been proved, namely $v\in H^{s+2,s/2+1}(\Omega^t)$, $s\in(1/2,1)$. Danchin an Mucha in \cite{DM} considered solutions $v$, $p$, $\varrho$ in some Besov spaces on $\R^n$ assuming additionally that the velocity $v$ tends to 0 at infinity and the density $\varrho$ tends to some positive constant $\varrho^*$ at infinity. In particular, they admit piece-wise-constant initial densities provided the jump at the interface is small enough. In \cite{DZ}, the global existence and uniqueness of solutions in the half-space $\R_+^n$, $n\ge 2$, has been established, with the initial density bounded and close enough to a positive constant, the initial velocity belonging to some critical Besov space and some smallness of data.

We also mention, that the Navier-Stokes system in a cylindrical domain with or without flux has been treated in some papers (see \cite{RZ1}, \cite{RZ2}, \cite{Z1}, \cite{Z2}), however, density dependent equations present a quite different problem and we can not transfer in an easy way such considerations and results. On the other hand, involving density means more physically realistic model and gives possibilities of applications in modelling, for example, a blood flow in veins or arteries.

\section{Notation and auxiliary results}\label{s2}

First we introduce the simplified notation.

\begin{definition}\label{d2.1}
For Lebesque and Sobolev spaces we set the following notation
$$\eqal{
&\|u\|_{L_p(Q)}=|u|_{p,Q},\quad \|u\|_{L_p(Q^t)}=|u|_{p,Q^t},\cr
&\|u\|_{L_q(0,t;L_p(Q))}=|u|_{p,q,Q^t},\cr}
$$
where $Q=\Omega$ or $Q=S\equiv\partial\Omega$ or $Q=\R^n$, $Q^t=Q\times(0,t)$, $p,q\in[1,\infty]$.\\
Let $H^s(Q)=W_2^s(Q)$. Then we denote
$$
\|u\|_{H^s(Q)}=\|u\|_{s,Q},\quad \|u\|_{W_p^s(Q)}=\|u\|_{s,p,Q},
$$
where $s\in\N$.\\
Additionally, we introduce
$$\eqal{
&\|u\|_{L_q(0,t;W_p^k(\Omega))}=\|u\|_{k,p,q,\Omega^t}, \cr
&\|u\|_{L_p(0,t;W_p^k(\Omega))}=\|u\|_{k,p,\Omega^t},\quad k\in\N,\ \ p,q\in[1,\infty].\cr}
$$
Moreover,
$$
\|u\|_{W_r^{1,1}(Q^t)}=\|u\|_{L_r(Q^t)}+\|u_{x}\|_{L_r(Q^t)}+\|u_{t}\|_{L_r(Q^t)},
$$
where $r\in[1,\infty]$ and $Q$ is equal either $\Omega$ or $S$ and
$$
\|u\|_{W_{r,s}^{1,1}(Q^t)}=\|u\|_{L_s(0,t;L_r(Q))}+\|u_{x}\|_{L_s(0,t;L_r(Q))}+ \|u_{t}\|_{L_s(0,t;L_r(Q))},
$$
where $s\in[1,\infty]$.

By $V_p^{2+{s}}(\Omega^t)$ we denote an anisotropic energy space in the $p$-approach with the norm
\begin{equation}
\eqal{
&\|u\|_{V_p^{2+{s}}(\Omega^t)}=\sup_{t'\le t}\|u(t')\|_{W_p^{2+{s}-2/p}(\Omega)}+\|u\|_{W_p^{2+{s},1+{s}/2}(\Omega^t)},\cr
&p\in[1,\infty],\ \ {s}\in(0,1).\cr}
\label{4.19}
\end{equation}
\end{definition}
By $C^\alpha(\Omega^T)$, $\alpha\in(0,1)$ we denote the H\"older space with the norm
$$\eqal{
\|u\|_{C^{\alpha}(\Omega^T)}&=\sup_{x,x',t} \frac{|u(x,t)-u(x',t)|}{|x-x'|^\alpha}\cr
&\quad+\sup_{x,t,t'}\frac{|u(x,t)-u(x,t')|}{|t-t'|^\alpha}+|u|_{\infty,\Omega^T}\cr}
$$
and the homogenous H\"older space with the norm
$$
\|u\|_{\dot C^\alpha(\Omega^T)}=\sup_{x,x',t}\frac{|u(x,t)-u(x',t)|}{|x-x'|^\alpha}+\sup_{x,t,t'} \frac{|u(x,t)-u(x,t')|}{|t-t'|^\alpha}.
$$
We also need a homogeneous Sobolev space
$$
\|u\|_{\dot W_{r,s}^{1,1}(\Omega^t)}=|u_{x}|_{s,r,\Omega^t}+|u_{t}|_{s,r,\Omega^t}.
$$
In this paper $\phi$ always denotes an increasing positive function. Moreover, $\phi$ is a generic function because it changes its form from formula to formula. The exponent $\bar{a}$ in the function $t^{\bar{a}}$ is assumed to be always positive.

\begin{definition}[Besov spaces (see \cite{BIN}, \cite{N}, \cite{G}, \cite{Tr1})]\label{d2.2}
There are many different but equivalent definitions of Besov spaces (see \cite{G} and \cite[Th.18.2 from Ch. 4, Sect. 18]{BIN}). In this paper it is assumed the following norm of anisotropic Besov space $B_{p,q}^{{\sigma},{\sigma}/2}(\Omega\times(\tau,T))$,
\begin{equation}\eqal{
&\|u\|_{B_{p,q}^{{\sigma},{\sigma}/2}(\Omega\times(\tau,T))}= \|u\|_{L_p(\Omega\times(\tau,T))}\cr
&\quad+\sum_{i=1}^n\bigg(\intop_0^{h_0} \frac{\|\Delta_i(h,\Omega)\partial_{x_i}^{[{\sigma}]}u\|_{L_p(\Omega\times(\tau,T))}^q} {h^{1+q({\sigma} -[{\sigma}])}}dh\bigg)^{1/q}\cr
&\quad+\bigg(\intop_0^{h_0} \frac{\|\Delta_0(h,(\tau,T))\partial_t^{[{\sigma}/2]}u\|_{L_p(\Omega\times(\tau,T))}^q}{ h^{1+q({\sigma}/2-[{\sigma} /2])}}dh\bigg)^{1/q},\cr}
\label{2.1}
\end{equation}
where
$$\eqal{
\Delta_i(h,\Omega)f(x,t)&=\begin{cases}f(x+h_i)-f(x)\ &{\rm if\ }[x,x+h_i]\subset\Omega,\cr
0&\rm otherwise,\cr\end{cases}\cr
\Delta_0(h,(\tau,T))f(x,t)&=\begin{cases}f(x,t+h)-f(x,t)&{\rm if\ }[t,t+h]\subset(\tau,T)\cr 0&\rm otherwise,\cr\end{cases}
\cr}
$$
where $h_i$ is the $i$-th coordinate of $h\in\R^n$, $\Omega\subset\R^n$, $\tau\in(-\infty,T)$, ${\sigma}\in\R_+$, $p,q\in[1,\infty]$.
\end{definition}

Following Lemma 7.44 from \cite{Ad} and for $p=q$ the norm (\ref{2.1}) is equivalent to the following
\begin{equation}\eqal{
&\|u\|_{\bar B_{p,p}^{{\sigma},{\sigma}/2}(\Omega\times(\tau,T))}= \|u\|_{L_p(\Omega\times(\tau,T))}\cr
&\quad+\bigg(\intop_\tau^Tdt\intop_\Omega dx'\intop_\Omega dx'' \frac{|D_{x'}^{[{\sigma}]}u(x',t)-D_{x''}^{[{\sigma}]}u(x'',t)|^p}{ |x'-x''|^{n+p({\sigma}-[{\sigma}])}}\bigg)^{1/p}\cr
&\quad+\bigg(\intop_\Omega dx\intop_\tau^Tdt'\intop_\tau^Tdt'' \frac{|\partial_{t'}^{[{\sigma}/2]}u(x,t')-\partial_{t''}^{[{\sigma}/2]}u(x,t'')|^p}{ |t'-t''|^{1+p({\sigma}/2-[{\sigma}/2])}}\bigg)^{1/p}.\cr}
\label{2.2}
\end{equation}
\begin{remark}[Sobolev-Slobodetskii spaces]  $\bar B_{p,p}^{{\sigma},{\sigma}/2}(\Omega\times(\tau,T))$ described in (\ref{2.2}) is denoted by $W_p^{{\sigma},{\sigma}/2}(\Omega\times(\tau,T))$ and the space is called the Sobolev-Slobodetskii space.

We also use isotropic Sobolev-Slobodetskii spaces $W_p^{\sigma}(\Omega)$ with the following norm
$$
\|u\|_{W_p^{\sigma}(\Omega)}=\|u\|_{L_p(\Omega)}+\bigg(\intop_\Omega\intop_\Omega \frac{|D_{x'}^{[{\sigma}]}u(x')-D_{x''}^{[{\sigma}]}u(x'')|^p}{|x'-x''|^{n+p({\sigma}-[{\sigma}])}}dx'dx''\bigg)^{1/p},
$$
where $\Omega\subset\R^n$, $D_x^k=\partial_{x_1}^{k_1}\dots\partial_{x_n}^{k_n}$, $k=\sum_{i=1}^nk_i$, $k_i,k\in\N\cup\{0\}$, $i=1,2$, $[{\sigma}]$ is the integer part of ${\sigma}\in\R_+$.
\end{remark}

\begin{definition}(Partition of unity, see \cite[Ch. 4, Sect. 4]{LSU}) \label{d2.3}
Let $\Omega\subset\R^3$ be the domain defined in Section \ref{s1}. We define two collections of open subsets $\{\omega^{(k)}\}$ and $\{\Omega^{(k)}\}$, $k\in{\frak M}\cup{\frak N}$, such that $\bar\omega^{(k)}\subset\Omega^{(k)}\subset\Omega$, $\bigcup_k\omega^{(k)}=\bigcup_k\Omega^{(k)}=\Omega$, $\bar\Omega^{(k)}\cap S=\phi$ for $k\in{\frak M}$, $\bar\Omega^{(k)}\cap S_1\not=\phi$ for $k\in{\frak N}_1$, $\bar\Omega^{(k)}\cap S_2\not=\phi$ for $k\in{\frak N}_2$ and $\bar\Omega^{(k)}$ is a neighborhood of a point $\xi\in L(a_i)=\partial S_2(a_i)$, $i=1,2$, for $k\in{\frak N}_3$.\\
Hence ${\frak N}={\frak N}_1\cup{\frak N}_2\cup{\frak N}_3$ and $\bar\Omega^{(k)}\cap S_i$, $k\in{\frak N}_i$, $i=1,2$, are located in a positive distance from the edge $L=\partial S_2$. We assume that at most $N_0$ of the $\Omega^{(k)}$ have nonempty intersections, $\sup_k{\rm diam}\Omega^{(k)}\le 2\lambda$, $\sup_k{\rm diam}\omega^{(k)}\le\lambda$ for some $\lambda>0$.

Let $\zeta^{(k)}(x)$ be a smooth function such that
$$
0\le\zeta^{(k)}(x)\le 1,\ \ \zeta^{(k)}(x)=1\quad {\rm for}\ \ x\in\omega^{(k)},\ \ \zeta^{(k)}(x)=0
$$
for $x\in\Omega\setminus\Omega^{(k)}$ and $|D_x^\nu\zeta^{(k)}(x)|\le c/\lambda^{|\nu|}$. Then $1\le\sum_k(\zeta^{(k)}(x))^2\le N_0$. Introducing the function
$$
\eta^{(k)}(x)=\frac{\zeta^{(k)}(x)}{\sum_l(\zeta^{(l)}(x))^2}
$$
we have that $\eta^{(k)}(x)=0$ for $x\in\Omega\setminus\Omega^{(k)}$, $\sum_k\eta^{(k)}(x)\zeta^{(k)}(x)=1$, $|D_x^\nu\eta^{(k)}(x)|\le c/\lambda^{|\nu|}$, where $D_x^\nu=\partial_{x_1}^{\nu_1}\partial_{x_2}^{\nu_2}\partial_{x_3}^{\nu_3}$, $|\nu|=\nu_1+\nu_2+\nu_3$, $\nu_i\in\N_0=\N\cup\{0\}$.
\end{definition}

Consider the problem
\begin{equation}\eqal{
&\varrho_t+v\cdot\nabla\varrho=0,\ \ \divv v=0\ \ \quad &{\rm in}\ \ \Omega^T,\cr
&\varrho|_{t=0}=\varrho(0) \quad &{\rm in}\ \ \Omega, \cr
&v\cdot\bar n=-d_1,\ \ \varrho=\varrho_1,\ \ d_1>0\quad &{\rm on}\ \ S_2(-a).\cr}
\label{2.3}
\end{equation}

\begin{lemma}\label{l2.4}
Assume that $\varrho(0)\in L_\infty(\Omega)$, $\varrho_1\in L_\infty(S_2^t(-a))$, $\frac{1}{\varrho(0)}\in L_\infty(\Omega)$, $\frac{1}{\varrho_1}\in L_\infty(S_2^t(-a))$. Then
\begin{equation}
|\varrho(t)|_{\infty,\Omega}\le|\varrho(0)|_{\infty,\Omega}+ |\varrho_1|_{\infty,S_2^t(-a)}\equiv\varrho^*.
\label{2.4}
\end{equation}
Moreover,
\begin{equation}
\varrho_*=\frac{\inf\varrho_1\cdot\inf\varrho(0)}{\inf\varrho_1+\inf\varrho(0)}\le \varrho.
\label{2.5}
\end{equation}
\end{lemma}

\begin{proof}
Multiply $(\ref{2.3})_1$ by $\varrho|\varrho|^{p-2}$, $p\in\R_+$, and integrate over $\Omega$. Then, we derive
$$
\frac{d}{dt}|\varrho|_{p,\Omega}^p+\intop_\Omega v\cdot\nabla|\varrho|^pdx=0.
$$
In view of the boundary conditions $(\ref{2.3})_2$ we conclude
$$
\frac{d}{dt}|\varrho|_{p,\Omega}^p\le\intop_{S_2(-a)}d_1|\varrho_1|^pdS_2.
$$
Next, integrating with respect to time yields
$$
|\varrho(t)|_{p,\Omega}^p\le\intop_{S_2^t(-a)}d_1|\varrho_1|^pdS_2dt'+ |\varrho(0)|_{p,\Omega}^p.
$$
Hence
$$
|\varrho(t)|_{p,\Omega}\le|d_1|_{\infty,S_2^t(-a)}^{1/p}|\varrho_1|_{p,S_2^t(-a)}+ |\varrho(0)|_{p,\Omega}.
$$
Passing with $p\to\infty$ implies (\ref{2.4}).

In order to achieve the second thesis of the lemma, we multiply $(\ref{2.3})_1$ by $\varrho|\varrho|^{-p-2}$, $p\in\R_+$ and integrate over $\Omega$. Then we obtain
$$
\frac{d}{dt}\bigg|\frac{1}{\varrho}\bigg|_{p,\Omega}^p\le\intop_{S_2(-a)}d_1 \bigg|\frac{1}{\varrho_1}\bigg|^pdS_2.
$$
Integrating with respect to time gives
$$
\left|\frac{1}{\varrho}\right|_{p,\Omega}^p\le\intop_{S_2^t(-a)}d_1 \bigg|\frac{1}{\varrho_1}\bigg|^pdS_2dt'+\bigg|\frac{1}{\varrho(0)}\bigg|_{p,\Omega}^p.
$$
Hence,
$$
\bigg|\frac{1}{\varrho}\bigg|_{p,\Omega}\le|d_1|_{\infty,S_2^t(-a)}^{1/p} \bigg|\frac{1}{\varrho_1}\bigg|_{p,S_2^t(-a)}+ \bigg|\frac{1}{\varrho(0)}\bigg|_{p,\Omega}.
$$
Passing with $p\to\infty$ yields
$$
\bigg|\frac{1}{\varrho}\bigg|_{\infty,\Omega}\le \bigg|\frac{1}{\varrho_1}\bigg|_{\infty,S_2^t(-a)}+ \bigg|\frac{1}{\varrho(0)}\bigg|_{\infty,\Omega}.
$$
This means that
$$
\inf\varrho\ge\frac{1}{\frac{1}{\inf\varrho_1}+\frac{1}{\inf\varrho(0)}}= \frac{\inf\varrho_1\cdot\inf\varrho(0)}{\inf\varrho_1+\inf\varrho(0)}.
$$
The above inequality implies (\ref{2.5}) and concludes the proof.
\end{proof}

\begin{lemma}[The Korn inequality, see \cite{SS}]\label{l2.5}
Assume that a function $w$ satisfies the following conditions: $E_\Omega(w)=|\D(w)|_{2,\Omega}^2<\infty$, $\divv w=0$, \\
$w\cdot\bar n|_S=0$ and $\Omega$ is not axially symmetric. Then there exists a constant $c$ independent of $w$ such that
\begin{equation}
\|w\|_{H^1(\Omega)}^2\le cE_\Omega(w).
\label{2.6}
\end{equation}
\end{lemma}

\begin{remark}\label{r2.6}
We recall some interpolations. Let ${\sigma}'<2+{\sigma}$, ${\sigma}\in(0,1)$. Then
\begin{equation}
\|u\|_{W_r^{{\sigma}',{\sigma}'/2}(\Omega^t)}\le\varepsilon \|u\|_{W_r^{2+{\sigma},1+{\sigma}/2}(\Omega^t)}+c(1/\varepsilon)|u|_{2,\Omega^t}.
\label{2.7}
\end{equation}
\end{remark}
Consider the Neumann problem
\begin{equation}\eqal{
&\Delta u=g\quad &{\rm in}\ \ \Omega,\cr
&\bar n\cdot\nabla u=f\quad &{\rm on}\ \ S,\cr
&\intop_\Omega udx=0,\cr}
\label{2.8}
\end{equation}
where $\Omega$ is a bounded domain described by (\ref{1.4}) and the following compa\-tibility condition holds:
$$ \int_{\Omega} g dx = \int_{S} f dS. $$

Let $G$ be the Green function to the Neumann problem (\ref{2.8}). It is a solution to the problem
\begin{equation}\eqal{
&\Delta_xG(x,y)=\delta(x,y)- \frac{1}{|\Omega|} \quad &{\rm in}\ \ \Omega,\cr
&\bar n\cdot\nabla_xG(x,y)=0\quad &{\rm on}\ \ S,\cr}
\label{2.9}
\end{equation}
where $\delta(x,y)$ is the delta function such that $\int_{\Omega} \delta\, dx = 1$ and $|\Omega| = {\rm meas}\ \Omega.$

\noindent Then any solution to problem (\ref{2.8}) can be expressed in the form
\begin{equation}
u(x)=\intop_\Omega G(x,y)g(y)dy-\intop_SG(x,y')f(y')dy',
\label{2.10}
\end{equation}
where $S=\partial\Omega$, $\Omega\subset\R^3$, and $G$ has a compact support with respect to $x$.

We are looking for the Green function in the form
$$ G(x,y) = \frac{1}{|x-y|} \eta(x,y) + g_0(x,y), $$
where a smooth function $\eta(x,y)$ is such that $\eta(x,y) =0$ for $|x-y|> \delta$ and $\eta(x,y)=1$ for $|x-y| < \delta/2$ with sufficiently small $\delta.$ Moreover, $g_0$ is a nonsingular function. The proof of existence of such a function can be found in \cite{S}.

Using an appropriate partition of unity and local flattening of the boun\-dary we can use in our considerations the Green function for problem (\ref{2.8}) in the half-space $\R_+^3=\{x\in\R^3\colon x_3>0\}$. In this case the Green function to problem (\ref{2.8}) has the form
\begin{equation}\eqal{
G(x,y)&=\frac{1}{4\pi}\frac{1}{\sqrt{(x_1-y_1)^2+(x_2-y_2)^2+(x_3-y_3)^2}}\cr
&\quad+\frac{1}{4\pi}\frac{1}{\sqrt{(x_1-y_1)^2+(x_2-y_2)^2+(x_3+y_3)^2}}.\cr}
\label{2.11}
\end{equation}
Hence
$$
\partial_{y_3}G(x,y)|_{y_3=0}=0.
$$
Then (\ref{2.10}) takes the form
\begin{equation}\eqal{
u(x)&=\frac{1}{4\pi}\intop_{\R_+^3}\bigg(\frac{1}{\sqrt{(x_\alpha-y_\alpha)^2+(x_3-y_3)^2}}\cr
&\quad+ \frac{1}{\sqrt{(x_\alpha-y_\alpha)^2+(x_3+y_3)^2}}\bigg)f(y)dy\cr
&\quad-\frac{1}{2\pi}\intop_{\R^2}\frac{1}{\sqrt{(x_\alpha-y_\alpha)^2+x_3^2}} g(x,y')dy',\cr}
\label{2.12}
\end{equation}
To guarantee that $(\ref{2.8})_2$ holds we recall the result
\begin{equation}
\frac{1}{2\pi}\intop_\R\intop_\R \frac{x_3}{(\sqrt{(x_\alpha-y_\alpha)^2+x_3})^3} dy_1dy_2=1,
\label{2.13}
\end{equation}
where $\alpha=1,2$ and the summation over repeated indices is assumed.

More results on the Green function to the Neumann problem in the half-space can be found in \cite{CP}.

To construct the Green function to problem (\ref{2.8}) in a bounded domain with sufficiently regular boundary we have to use an appropriate partition of unity. Then localizing problem (\ref{2.8}) and flattening locally the boundary we obtain problem (\ref{2.8}) in the half-space $x_3>0$. It can be described by formula (\ref{2.12}) modulo some non-singular term.

The existence of such Green function can be proved by the technique of regularizer. We have to emphasize, that singularity of the Green function does not change after local transformations.

A construction of the Green function in a bounded domain can be found in \cite{S}.

\begin{lemma}\label{l2.7}
Consider the problem (\ref{2.8}). Let $g=0$. Assume that the force $f$ has a compact support, $\int_S f dS = 0 $ and $f\in L_p(S)$. Then the solution $u\in L_r(\Omega)$ and
\begin{equation}
|u|_{r,\Omega}\le c|f|_{p,S},
\label{2.14}
\end{equation}
where $p>\frac{2r}{r+3}$, $r \ge p \ge 1.$
\end{lemma}

\begin{proof}
Using an appropriate partition of unity we can write (\ref{2.10}) in a local form
\begin{equation}
u(x)=\intop_{\R_+^3}dx\intop_{\R^2}\frac{1}{\sqrt{(x'-y')^2+x_3^2}}f(y')dy',
\label{2.15}
\end{equation}
where $x'=(x_1,x_2)$, $y'=(y_1,y_2)$ and $u$, $f$ have compact supports. Estima\-ting, we have
\begin{equation}
|u(x)|_{r,\R_+^3}\le\bigg(\intop_{\R_+^3}dx\bigg|\intop_{\R^2} \frac{1}{\sqrt{(x'-y')^2+x_3^2}}f(y')dy'\bigg|^r\bigg)^{1/r}\equiv I_1.
\label{2.16}
\end{equation}
By the Minkowski inequality, we obtain
$$
I_1\le\bigg(\intop_{\R^2}dx'\bigg|\intop_{\R^2}\bigg(\intop_{\R_+}dx_3\bigg| \frac{1}{\sqrt{(x'-y')^2+x_3^2}}\bigg|^r\bigg)^{1/r}f(y')dy'\bigg|^r\bigg)^{1/r}\equiv I_2.
$$
Let $a=|x'-y'|$. Introduce a new variable $z$ such that $x_3=az$. Then
$$
\bigg(\intop_{\R_+}dx_3\bigg|\frac{1}{\sqrt{a^2+x_3^2}}\bigg|^r\bigg)^{1/r}= \left(\frac{a}{ a^r}\right)^{1/r}\bigg(\intop_{\R_+}dz\bigg|\frac{1}{\sqrt{1+z^2}}\bigg|^r\bigg)^{1/r} \le \frac{c}{a^{1-1/r}}.
$$
Consequently,
$$
I_2\le c\bigg(\intop_{\R^2\cap\supp u}dx'\bigg|\intop_{\R^2\cap\supp f}\bigg(\frac{1}{\sqrt{x'-y')^2}}\bigg)^{1-1/r}f(y')dy'\bigg|^r\bigg)^{1/r}\equiv I_3.
$$
By the Young inequality (see \cite[Ch. 1, Sect. 2.14]{BIN}) we have
$$
I_3\le c|f|_{p,\R^2},
$$
where $1-1/p+1/r=1/{\sigma}$ and $K(x'-y')=\big(\frac{1}{\sqrt{(x'-y')^2}}\big)^{1-1/r}\in L_{\sigma}(\R^2)$ for $x'$, $y'$ belonging to some compact set, so ${\sigma}<\frac{2}{1-1/r}$.

Then
$$
1-\frac{1}{p}+\frac{1}{r}>\frac{r-1}{2r}\quad {\rm so}\quad p>\frac{2r}{r+3}.
$$
Using properties of the partition of unity, we conclude the proof.
\end{proof}
\begin{remark} 
In Lemma~\ref{l2.7} we have a low regularity of function $f$ such that $f \in L_p(S), p\ge 1.$ Therefore, to find estimates for solutions to 
(\ref{2.8}) we need to apply the potential method, thus the existence and properties of the Green function must be used. The restriction $p > \frac{2r}{r+3}$ follows from the Young inequality (see \cite{BIN}, Sect. 2.14), so any less restrictive interplay between $p$ and $r$ can not be expected. 
\end{remark}

\section{Method of successive approximations}\label{s3}

To prove the existence of local solutions to problem (\ref{1.1}) we use the method of successive approximations that we describe here and proceed in Sections~\ref{s5}-\ref{s6}.

Namely, we assume that $v_n(x,t)\subset W_{\sigma}^{2+{s},1+{s}/2}(\Omega^T)$ is given, $v_n$ is divergence free and there exists a positive constant $d_*=d_*(n)$ such that
\begin{equation}
v_{n;3}\ge d_*.
\label{3.1}
\end{equation}
For given $v_n$, we consider the function $\varrho_n$ - a solution to the problem
\begin{equation}\eqal{
&\varrho_{n,t}+v_n\cdot\nabla\varrho_n=0\quad &{\rm in}\ \ \Omega^T,\cr
&\divv v_n=0\quad &{\rm in}\ \ \Omega^T,\cr
&\varrho_n|_{S_2^T(-a)}=\varrho_1,\cr
&v_n\cdot\bar n|_{S_2^T(-a)}=-d_1,\ \ v_n\cdot\bar n|_{S_2^T(a)}=d_2,\ \ d_i>0,\ \ i=1,2,\cr
&\varrho_n|_{t=0}=\varrho(0).\cr}
\label{3.2}
\end{equation}
Next, for given $\varrho_n$ and $v_n$, we calculate $v_{n+1}$ and $p_{n+1}$ as solutions to the linear problem
\begin{equation}\eqal{
&\varrho_nv_{n+1,t}+\varrho_nv_n\cdot\nabla v_{n+1}-\nu\divv\D(v_{n+1})+\nabla p_{n+1}=\varrho_nf\quad &{\rm in}\ \ \Omega^T,\cr
&\divv v_{n+1}=0\quad &{\rm in}\ \ \Omega^T,\cr
&v_{n+1}\cdot\bar n=0,\ \ \nu\bar n\cdot\D(v_{n+1})\cdot\bar\tau_\alpha+\gamma v_{n+1}\cdot\bar\tau_\alpha=0,\ \ \alpha=1,2\quad &{\rm on}\ \ S_1^T,\cr
&v_{n+1}\cdot\bar n=d,\ \ \bar n\cdot\D(v_{n+1})\cdot\bar\tau_\alpha=0,\ \ \alpha=1,2\quad &{\rm on}\ \ S_2^T, \cr
& \int_{\Omega} p_{n+1} dx =0, \cr
&v_{n+1}|_{t=0}=v(0)\quad &{\rm in}\ \ \Omega.\cr}
\label{3.3}
\end{equation}
We are going to analyze solutions to (\ref{3.3}) and conclude some estimates. The first step is to derive the energy type estimates and this requires the homogeneous Dirichlet boundary conditions. For this purpose we consider a 'correction' function $\varphi$ as a solution to the Neumann problem

\begin{equation}\eqal{
&\Delta\varphi=0, \cr
&\bar n\cdot\nabla\varphi|_{S_1}=0,\ \ \bar n\cdot\nabla\varphi|_{S_2(-a)}= -d_1,\ \ \bar n\cdot\nabla\varphi|_{S_2(a)}= d_2\cr}
\label{3.4}
\end{equation}
and the following compatibility condition holds $$ - \int_{S_2(-a)} d_1 dS_2 + \int_{S_2(a)} d_2 dS_2=0. $$
Introducing the function
\begin{equation}
w_{n+1}=v_{n+1}-\nabla\varphi
\label{3.5}
\end{equation}
we see that $w_{n+1}$ and $p_{n+1}$ satisfy the system with homogeneous boundary conditions:
\begin{equation}\eqal{
&\varrho_nw_{n+1,t}+\varrho_nv_n\cdot\nabla w_{n+1}-\nu\divv\D(w_{n+1})+\nabla p_{n+1}\cr
&\quad =-\varrho_n\nabla\varphi_t- \varrho_nv_n\cdot\nabla\nabla\varphi+\nu\divv\D(\nabla\varphi)+\varrho_nf \quad &{\rm in} \ \Omega^T, \cr
&\divv w_{n+1}=0 \quad &{\rm in} \ \Omega^T, \cr
&w_{n+1}\cdot\bar n=0 \quad &{\rm on}\ \ S^T,\cr
&\nu\bar n\cdot\D(w_{n+1})\cdot\bar\tau_\alpha+\gamma w_{n+1}\cdot\bar\tau_\alpha= -\nu\bar n\cdot\D(\nabla\varphi)\cdot\bar\tau_\alpha\cr
&\quad-\gamma\nabla\varphi\cdot\bar\tau_\alpha,\ \ \alpha=1,2 \quad &{\rm on}\ \ S_1^T,\cr
&\bar n\cdot\D(w_{n+1})\cdot\bar\tau_\alpha=-\bar n\cdot\D(\nabla\varphi)\cdot\bar\tau_\alpha,\ \ \alpha=1,2\quad &{\rm on}\ \ S_2^T,\cr
&w_{n+1}|_{t=0}=v(0)-\nabla\varphi|_{t=0}\equiv w(0) \quad &{\rm in} \ \Omega .\cr}
\label{3.6}
\end{equation}

\begin{lemma}\label{l3.1}
For $\varphi$, a solution to (\ref{3.4}), the following estimate holds
\begin{equation}
\|\nabla\varphi\|_{W_{\sigma}^{2+{s},1+{s}/2}(\Omega^t)}\le c\sum_{i=1}^2\|d_i\|_{W_{\sigma}^{2+{s}-1/{\sigma},1+{s}/2}(S_2^t(a_i))}.
\label{3.7}
\end{equation}
\end{lemma}
\begin{proof} 
Problem (\ref{3.4}) is the elliptic problem so the time regularity follows from the set of problems 
\begin{equation} \eqal{
& \Delta \partial_t^l \varphi = 0, \ \bar{n} \cdot \Delta \partial_t^l \varphi|_{S_1} = 0, \cr & \bar{n} \cdot \Delta \partial_t^l \varphi|_{S_2(-a)} = - \partial_t^l d_1, \ \bar{n} \cdot \Delta \partial_t^l \varphi|_{S_2(a)} = \partial_t^l d_2,
 \cr}
\label{3.8}
\end{equation}
where $l \le 1+ \frac{s}{2}.$ The derivative $\partial_t^{1+s/2}$ must be treated as a fractional derivative. Hence, for any $l \le  1+ \frac{s}{2}$ we have the elliptic problem (\ref{3.8}). The existence of solutions to problem  (\ref{3.8}) and appropriate estimates in Besov spaces can be found in \cite{Tr2}, Ch.5, Th. 5.7.1. Moreover, the estimate (\ref{3.7}) also follows from \cite{ZZ}, where the heat equation and the Dirichlet problem are replaced by the Neumann problem for the Laplace equation.   
\end{proof}

\section{Estimates for first derivatives of $\varrho_n$}\label{s4}

In this Section we consider the density problem (\ref{3.2}) where $v_n$ is a given function belonging to $W_{\sigma}^{2+{s},1+{s}/2}(\Omega^t)$. Moreover, we assume that there exists a positive constant $d_*$ such that
\begin{equation}
v_{n;3}\ge d_*.
\label{4.1}
\end{equation}
To simplify presentation we drop the index $n$ in (\ref{3.2}). However, in the final result of this Section, the index $n$ will be recalled. Hence (\ref{3.2}) takes the form
\begin{equation}\eqal{
&\varrho_t+v\cdot\nabla\varrho=0, \cr
&\divv v=0, \cr
&\varrho|_{t=0}=\varrho(0),\ \ \varrho|_{S_2(-a)}=\varrho_1, \cr
&v\cdot\bar n|_{S_2(-a)}=-d_1,\ \ v\cdot\bar n|_{S_2(a)}=d_2,\cr}
\label{4.2}
\end{equation}
where $v$ is a given function.

Let
\begin{equation}\eqal{
&x'=(x_1,x_2),\ \  \varrho_{x'}=(\varrho_{x_1},\varrho_{x_2}),\ \ |\varrho_{x'}|=|\varrho_{x_1}|+|\varrho_{x_2}|,\cr
&v'=(v_1,v_2),\ \ v'_{x'}=(v_{1,x_1},v_{1,x_2},v_{2,x_1},v_{2,x_2}),\cr
&|v'_{x'}|=|v_{1,x_1}|+|v_{1,x_2}|+|v_{2,x_1}|+|v_{2,x_2}|\cr}
\label{4.3}
\end{equation}
and the summation convention of repeated indices is assumed.

\begin{lemma}\label{l4.1}
Assume that $v_{x'},v_t\in L_1(0,t;L_\infty(\Omega))$, $v\in L_\infty(\Omega^t)$, $d_1\in L_\infty(S_2^t(-a))$, $\varrho_{1,x'},\varrho_{1,t}\in L_r(S_2^t(-a))$, $\varrho_x(0),\varrho_t(0)\in L_r(\Omega)$, $r\ge 2$. Assume that $v_3\ge d_*>0$, where $d_*$ is a constant. Then $\varrho$ - a solution to the problem (\ref{4.2}), satisfies
\begin{equation}\eqal{
&|\varrho_{x'}|_{r,\Omega}+|\varrho_t|_{r,\Omega}\le\exp\bigg[\frac{c}{r}\bigg(1+\frac{1}{d_*^r}\bigg)
\intop_0^t(|v_{x'}|_{\infty,\Omega}+|v_t|_{\infty,\Omega})\cdot\cr
&\quad\cdot(1+|v'|_{\infty,\Omega}^r)dt'\bigg]\cdot\cr
&\quad\cdot\bigg[|d_1|_{\infty,S_2^t(-a)}^{1/r} \bigg(\intop_{S_2^t(-a)}(|\varrho_{1,x'}|^r+|\varrho_{1,t}|^r)dS_2dt'\bigg)^{1/r}\cr
&\quad+|\varrho_{x'}(0)|_{r,\Omega}+|\varrho_t|_{r,\Omega}\bigg].\cr}
\label{4.4}
\end{equation}
\end{lemma}

\begin{proof}
We differentiate $(\ref{2.4})_1$ with respect to $x_\alpha$, $\alpha=1,2$, then multiply by $\varrho_{x_\alpha}|\varrho_{x_\alpha}|^{r-2}$, $r\ge 2$, and integrate over $\Omega$. As the result, we obtain
\begin{equation}
\frac{1}{r}\frac{d}{dt}|\varrho_{x_\alpha}|_{r,\Omega}^r+\frac{1}{r}\intop_\Omega v\cdot\nabla|\varrho_{x_\alpha}|^rdx+\intop_\Omega v_{x_\alpha}\cdot\nabla\varrho\varrho_{x_\alpha}|\varrho_{x_\alpha}|^{r-2}dx=0.
\label{4.5}
\end{equation}
Using that $v$ is divergence free, the second term in (\ref{4.5}) equals
$$
\frac{1}{r}\intop_\Omega\divv(v|\varrho_{x_\alpha}|^r)dx=-\frac{1}{r}\intop_{S_2(-a)}d_1|\varrho_{1,x_\alpha}|^rdS_2+\frac{1}{r}
\intop_{S_2(a)}d_2|\varrho_{x_\alpha}|^rdS_2.
$$
The last integral in (\ref{4.5}) has the form
$$
\sum_{\beta=1}^2\intop_\Omega v_{\beta,x_\alpha}\varrho_{x_\beta}\varrho_{x_\alpha}|\varrho_{x_\alpha}|^{r-2}dx+ \intop_\Omega v_{3,x_\alpha}\varrho_{x_3}\varrho_{x_\alpha}|\varrho_{x_\alpha}|^{r-2}dx.
$$
In view of the above expressions we derive from (\ref{4.5}) the inequality
\begin{equation}\eqal{
&\frac{1}{r}\frac{d}{dt}|\varrho_{x_\alpha}|_{r,\Omega}^r\le\frac{1}{r}\intop_{S_2(-a)}d_1|\varrho_{1,x_\alpha}|^rdS_2\cr
&\quad+\sum_{\beta=1}^2\intop_\Omega|v_{\beta,x_\alpha}|\,|\varrho_{x_\beta}|\, |\varrho_{x_\alpha}|^{r-1}dx+\intop_\Omega|v_{3,x_\alpha}|\,|\varrho_{x_3}|\, |\varrho_{x_\alpha}|^{r-1}dx,\cr}
\label{4.6}
\end{equation}
where $\alpha=1,2$.

Using the notation (\ref{4.3}) we can write (\ref{4.6}) as follows
\begin{equation}\eqal{
&\frac{1}{r}\frac{d}{dt}|\varrho_{x'}|_{r,\Omega}^r\le\frac{1}{r}\intop_{S_2(-a)}d_1|\varrho_{1,x'}|^rdS_2\cr
&\quad+|v'_{x'}|_{\infty,\Omega}|\varrho_{x'}|_{r,\Omega}^r+ |v_{3,x'}|_{\infty,\Omega}|\varrho_{x_3}|_{r,\Omega}|\varrho_{x'}|_{r,\Omega}^{r-1}.\cr}
\label{4.7}
\end{equation}
Differentiate $(\ref{4.2})_1$ with respect to $t$,  multiply by $\varrho_t|\varrho_t|^{r-2}$ and integrate over $\Omega$. Then we have
\begin{equation}\eqal{
&\frac{1}{r}\frac{d}{dt}|\varrho_t|_{r,\Omega}^r\le\frac{1}{r}\intop_{S_2(-a)}d_1|\varrho_{1,t}|^rdS_2-\frac{1}{r}\intop_{S_2(a)}d_2|\varrho_t|^rdS_2\cr
&\quad+\intop_\Omega v_t\cdot\nabla\varrho\varrho_t|\varrho_t|^{r-2}dx=0.\cr}
\label{4.8}
\end{equation}
Simplifying, (\ref{4.8}) implies the inequality
\begin{equation}\eqal{
&\frac{1}{r}\frac{d}{dt}|\varrho_t|_{r,\Omega}^r\le\frac{1}{r}\intop_{S_2(-a)}d_1|\varrho_t|^rdS_2+|v'_t|_{\infty,\Omega} |\varrho_{x'}|_{r,\Omega}|\varrho_t|_{r,\Omega}^{r-1}\cr
&\quad+|v_{3,t}|_{\infty,\Omega}|\varrho_{x_3}|_{r,\Omega} |\varrho_t|_{r,\Omega}^{r-1}.\cr}
\label{4.9}
\end{equation}
From (\ref{4.7}) and (\ref{4.9}) we have
\begin{equation}\eqal{
&\frac{1}{r}\frac{d}{dt}(|\varrho_{x'}|_{r,\Omega}^r+|\varrho_t|_{r,\Omega}^r)\le\frac{1}{r}\intop_{S_2(-a)}d_1(|\varrho_{1,x'}|^r+|\varrho_{1,t}|^r)dS_2\cr
&\quad+|v'_{x'}|_{\infty,\Omega}|\varrho_{x'}|_{r,\Omega}^r+ |v_{3,x'}|_{\infty,\Omega}|\varrho_{x_3}|_{r,\Omega}|\varrho_{x'}|_{r,\Omega}^{r-1}\cr
&\quad+|v'_t|_{\infty,\Omega}|\varrho_{x'}|_{r,\Omega} |\varrho_t|_{r,\Omega}^{r-1}+|v_{3,t}|_{\infty,\Omega}|\varrho_{x_3}|_{r,\Omega} |\varrho_t|_{r,\Omega}^{r-1}.\cr}
\label{4.10}
\end{equation}
We transform this into the inequality
$$\eqal{
&\frac{1}{r}\frac{d}{dt}(|\varrho_{x'}|_{r,\Omega}^r+|\varrho_t|_{r,\Omega}^r) \le\frac{1}{r}\intop_{S_2(-a)}d_1(|\varrho_{1,x'}|^r+|\varrho_{1,t}|^r)dS_2\cr
&\quad+|v_{x'}|_{\infty,\Omega}\bigg(\frac{2r-1}{r}|\varrho_{x'}|_{r,\Omega}^r+\frac{1}{r}|\varrho_{x_3}|_{r,\Omega}^r\bigg)\cr
&\quad+|v_t|_{\infty,\Omega}\bigg(\frac{1}{r}|\varrho_{x'}|_{r,\Omega}^r+\frac{1}{r}|\varrho_{x_3}|_{r,\Omega}^r+\frac{2(r-1)}{r}
|\varrho_t|_{r,\Omega}^r\bigg).\cr}
$$
Finally, rearranging, we infer the formula
\begin{equation}\eqal{
&\frac{1}{r}\frac{d}{dt}(|\varrho_{x'}|_{r,\Omega}^r+|\varrho_t|_{r,\Omega}^r)\le\frac{1}{r}\intop_{S_2(-a)}d_1(|\varrho_{1,x'}|^r+|\varrho_{1,t}|^r)dS_2\cr
&\quad+(|v_{x'}|_{\infty,\Omega}+ |v_t|_{\infty,\Omega})\bigg(2|\varrho_{x'}|_{r,\Omega}^r+\frac{2}{ r}|\varrho_{x_3}|_{r,\Omega}^r+\frac{2(r-1)}{r}|\varrho_t|_{r,\Omega}^r\bigg).\cr}
\label{4.11}
\end{equation}
We express $(\ref{4.2})_1$ as the following equation
\begin{equation}
\varrho_t+v_\alpha\varrho_{x_\alpha}+v_3\varrho_{x_3}=0.
\label{4.12}
\end{equation}
Since $v_3\ge d_*>0$, we calculate
\begin{equation}
\varrho_{x_3}=-\frac{1}{v_3}(\varrho_t+v_\alpha\varrho_{x_\alpha}).
\label{4.13}
\end{equation}
Using (\ref{4.13}) in (\ref{4.11}) yields
\begin{equation}\eqal{
&\frac{1}{r}\frac{d}{dt}(|\varrho_{x'}|_{r,\Omega}^r+|\varrho_t|_{r,\Omega}^r) \le\frac{1}{r}\intop_{S_2(-a)}d_1(|\varrho_{1,x'}|^r+|\varrho_{1,t}|^r)dS_2\cr
&\quad+2(|v_{x'}|_{\infty,\Omega}+|v_t|_{\infty,\Omega})\bigg[ |\varrho_x|_{r,\Omega}^r+\frac{1}{ rd_*^r}|\varrho_t+v_\alpha\varrho_{x_\alpha}|^r\cr
&\quad+\frac{r-1}{r}|\varrho_t|_{r,\Omega}^r\bigg],\cr}
\label{4.14}
\end{equation}
where the expression under the square bracket is bounded by
$$
|\varrho_{x'}|_{r,\Omega}^r+\frac{2^r}{rd_*^r}(|\varrho_t|_{r,\Omega}^r+ |v'|_{\infty,\Omega}^r|\varrho_{x'}|_{r,\Omega}^r)+\frac{r-1}{r}|\varrho_t|_{r,\Omega}^r.
$$
Hence, (\ref{4.14}) takes the form
\begin{equation}\eqal{
&\frac{1}{r}\frac{d}{dt}(|\varrho_{x'}|_{r,\Omega}^r+|\varrho_t|_{r,\Omega}^r) \le\frac{1}{r}\intop_{S_2(-a)}d_1(|\varrho_{1,x'}|^r+|\varrho_{1,t}|^r)dS_2\cr
&\quad+2(|v_{x'}|_{\infty,\Omega}+|v_t|_{\infty,\Omega})\bigg[\bigg(1+\frac{2^r}{rd_*^r}|v'|_{\infty,\Omega}^r\bigg)|\varrho_{x'}|_{r,\Omega}^r\cr
&\quad+\bigg(\frac{2^r}{rd_*^r}+\frac{r-1}{r}\bigg)|\varrho_t|_{r,\Omega}^r\bigg].\cr}
\label{4.15}
\end{equation}
Let
\begin{equation}
X_r(t)=(|\varrho_{x'}|_{r,\Omega}^r+|\varrho_t|_{r,\Omega}^r)^{1/r}.
\label{4.16}
\end{equation}
Then for any finite $r$, (\ref{4.15}) implies the inequality
\begin{equation}\eqal{
&\frac{d}{dt}X_r^r(t)\le\intop_{S_2(-a)}d_1(|\varrho_{1,x'}|^r+|\varrho_{1,t}|^r)dS_2\cr
&\quad+c(r)\bigg(1+\frac{1}{d_*^r}\bigg)(|v_{x'}|_{\infty,\Omega}+|v_t|_{\infty,\Omega}) (1+|v'|_{\infty,\Omega}^r)X_r^r(t).\cr}
\label{4.17}
\end{equation}
Integrating with respect to time yields
$$\eqal{
X_r^r(t)&\le\exp\bigg(c(r)\bigg(1+\frac{1}{d_*^r}\bigg)\intop_0^t(|v_{x'}|_{\infty,\Omega}+|v_t|_{\infty,\Omega}) (1+|v'|_{\infty,\Omega}^r)dt'\bigg)\cdot\cr
&\quad\cdot\bigg[\intop_0^tdt'\intop_{S_2(-a)}d_1(|\varrho_{1,x'}|^r+ |\varrho_{1,t}|^r)dS_2+X_r^r(0)\bigg].\cr}
$$
The above inequality implies (\ref{4.4}) and concludes the proof of Lemma~\ref{l4.1}.
\end{proof}

\begin{remark}\label{r4.2}
From (\ref{4.13}) we have
\begin{equation}
|\varrho_{x_3}|_{r,\infty,\Omega^t}\le\frac{1}{d_*}(1+|v'|_{\infty,\Omega^t})X_r(t).
\label{4.18}
\end{equation}
\end{remark}
\noindent We introduce the space $V_{\sigma}^{2+{s}}(\Omega^t)$ through the norm given in (\ref{4.19}):
$$
\|u\|_{V_{\sigma}^{2+{s}}(\Omega^t)}=\sup_{t'\le t}\|u(t')\|_{W_{\sigma}^{2+{s}-2/{\sigma}}(\Omega)}+ \|u\|_{W_{\sigma}^{2+{s},1+{s}/2}(\Omega^t)},
$$
where ${s}\in(0,1)$, ${\sigma}\in[1,\infty]$. Recalling the index $n$ we have

\begin{theorem}\label{t4.3}
Assume that $v_n\in V_{\sigma}^{2+{s}}(\Omega^t)$, $3/{\sigma}<{s}$, ${\sigma}>3$, $d_1\in\break L_\infty(S_2^t(-a))$, $\varrho_1\in W_r^1(S_2^t(-a))$, $\varrho(0)\in W_r^1(\Omega)$, $v(0)\in W_r^1(\Omega)$, $r\ge 2$.
Let $X_{r;n} \equiv X_r(\varrho_n)$ be defined by (\ref{4.16}) for $\varrho_n$ - solutions to (\ref{4.2}). Let $v_{3;n}\ge d_*>0$, $d_*=d_*(n)$.\\
Then there exists an increasing positive function $\phi_1$ such that
\begin{equation}\eqal{
X_{r;n}(t)&\le\phi_1(t^{\bar{a}}\|v_n\|_{V_{\sigma}^{2+{s}}(\Omega^t)},\cr
&\quad t^{\bar{a}}\|p_{n,x_3}\|_{W_{\sigma}^{{s},{s}/2}(\Omega^t)}, |f_3|_{\infty,{\sigma},\Omega^t})\bar d_1\cr
&=\phi\bigg(\exp\bigg[\frac{t^{\bar{a}}}{\varrho_*}(\|p_{n,x_3}\|_{W_{\sigma}^{{s},{s}/2}(\Omega^t)}\cr
&\quad+\|f_3\|_{L_{\sigma}(0,t;L_\infty(\Omega))})\bigg]\bigg) \phi(t^{\bar{a}}\|v_n\|_{V_{\sigma}^{2+{s}}(\Omega^t)})\bar d_1, \cr}
\label{4.20}
\end{equation}
(see Lemma \ref{l8.1} for $n$-th step) where $\bar{a}>0$ and
\begin{equation}\eqal{
\bar d_1&=|d_1|_{\infty,S_2^t(-a)}(|\varrho_{1,x}|_{r,S_2^t(-a)}+ |\varrho_{1,t}|_{r,S_2^t(-a)})\cr
&\quad+\|\varrho_x(0)\|_{L_r(\Omega)}(1+\|v(0)\|_{W_r^1(\Omega)}).\cr}
\label{4.21}
\end{equation}
Moreover,
\begin{equation}
|\varrho_{n,x_3}|_{r,\infty,\Omega^t}\le \frac{1}{d_*(n)}(1+\|v_n\|_{V_{\sigma}^{2+{s}}(\Omega^t)})X_r(t).
\label{4.22}
\end{equation}
\end{theorem}

\begin{proof}
First, we examine the argument of exp function in (\ref{4.4}). In the proof the index $n$ for $\varrho$ and $v$ is dropped.
$$\eqal{
&\intop_0^t(|v_x|_{\infty,\Omega}+|v_t|_{\infty,\Omega}) (1+|v|_{\infty,\Omega}^r)dt'\cr
&\le(1+\sup_t|v|_{\infty,\Omega}^r)t^{\frac{\sigma-1}{\sigma}}\bigg(\intop_0^t (|v_x|_{\infty,\Omega}^{\sigma}+|v_t|_{\infty,\Omega}^{\sigma} )dt'\bigg)^{1/{\sigma}}\cr
&\le c(1+\sup_t\|v\|_{W_{\sigma}^{2+{s}-2/\sigma}(\Omega)}^r)t^\frac{\sigma-1}{\sigma}\bigg(\intop_0^t \|v\|_{W_{\sigma} ^{2+{s}}(\Omega)}^{\sigma}dt'\bigg)^{1/{\sigma}}\cr
&\equiv\phi_1(t^{\bar{a}}\|v\|_{V_{\sigma}^{2+{s}}(\Omega^t)}),\cr}
$$
where we used that $3/{\sigma}<{s}$.

From (\ref{4.18}) we have
$$
|\varrho_{x_3}|_{r,\infty,\Omega^t}\le\frac{1}{d_*}(1+\sup_t|v|_{\infty,\Omega})X_r(t).
$$
Hence (\ref{4.22}) follows. This ends the proof.
\end{proof}

\begin{corollary} Estimates (\ref{4.20}) and (\ref{4.22}) imply
\begin{equation}\eqal{
&|\varrho_{n,x}|_{r,\Omega}+|\varrho_{n,t}|_{r,\Omega}\le\bigg[\frac{1}{d_*}(1+\|v_n\|_{V_{\sigma}^{2+{s}}(\Omega^t)})
+1\bigg]\cdot\cr
&\quad\cdot(\phi_1(t^{\bar{a}}\|v_n\|_{V_{\sigma}^{2+{s}}(\Omega^t)},t^{\bar{a}} \|p_{n,x_3}\|_{W_{\sigma}^{{s},{s}/2}(\Omega^t)},|f_3|_{\infty,{\sigma},\Omega^t}))\bar d_1,\cr}
\label{4.23}
\end{equation}
where $d_*=d_*(n)$ which is calculated in Lemma \ref{l8.1} for the $n$-th step.
\end{corollary}

\section{Estimates for solutions to (\ref{3.3})}\label{s5}

In this Section we consider solutions to problem~(\ref{3.3}). At first, we derive the energy estimate for solutions to problem (\ref{3.6}) and then our main issue is the proof of Lemma~\ref{l5.2} where the inequality for higher norms of solutions to (\ref{3.3}) in Sobolev and Sobolev-Slobodetskii spaces is established.

\begin{lemma}\label{l5.1}
Assume that $d_1\in L_\infty(S_2^t(-a))$, $\varrho_1\in L_6(0,t;L_3(S_2(-a))$, $v_n\in L_4(0,t;L_3(\Omega))$, $\nabla\varphi\in W_2^{2,1}(\Omega^t)$, $f\in L_2(0,t;L_{6/5}(\Omega))$, $w(0)\in L_2(\Omega)$.\\
Then, for solutions to (\ref{3.6}),
\begin{equation}\eqal{
&|w_{n+1}(t)|_{2,\Omega}^2+\|w_{n+1}\|_{1,2,\Omega^t}^2\le\exp(c |d_1|_{\infty,S_2^t(-a)}^6|\varrho_1|_{3,6,S_2^t(-a)}^6)\cdot\cr
&\quad\cdot[(1+|v_n|_{3,4,\Omega^t}^2)\|\nabla\varphi\|_{W_2^{2,1}(\Omega^t)}+ |f|_{6/5,2,\Omega^t}^2+|w(0)|_{2,\Omega}^2].\cr}
\label{5.1}
\end{equation}
\end{lemma}

\begin{proof}
Multiply $(\ref{3.6})_1$ by $w_{n+1}$, integrate over $\Omega$ and use boundary conditions $(\ref{3.6})_{3,4,5}$. Employing $(\ref{3.2})_1$ we obtain
\begin{equation}\eqal{
&\frac{1}{2}\frac{d}{dt}\intop_\Omega\varrho_n|w_{n+1}|^2dx+\nu |\D(w_{n+1})|_{2,\Omega}^2+\gamma\intop_{S_1}|w_{n+1}\cdot \bar\tau_\alpha|^2dS_1\cr
&=\intop_{S_2}\varrho_1d_1w_{n+1}^2dS_2-\intop_\Omega\varrho_n\nabla\varphi_t w_{n+1}dx+\intop_\Omega\varrho_n v_n\nabla^2\varphi w_{n+1}dx\cr
&\quad+\nu\intop_\Omega\D(\nabla\varphi)\cdot\nabla w_{n+1}dx-\gamma\intop_{S_1} \nabla\varphi\cdot\bar\tau_\alpha w_{n+1}\cdot\bar\tau_\alpha dS_1\cr
&\quad+\intop_\Omega\varrho_nf\cdot w_{n+1}dx.\cr}
\label{5.2}
\end{equation}
Using the Korn inequality (see Lemma \ref{l2.5}) and Lemma \ref{l2.4} applied to problem (\ref{3.2}), we have
\begin{equation}\eqal{
&\frac{d}{dt}\intop_\Omega\varrho_nw_{n+1}^2dx+\|w_{n+1}\|_{1,\Omega}^2+ |w_{n+1}\cdot\bar\tau_\alpha|_{2,S_1}^2\cr
&\le c\intop_{S_2}\varrho_1d_1w_{n+1}^2dS_2+c(|\nabla\varphi_t|_{6/5,\Omega}^2+ |v_n|_{3,\Omega}^2|\nabla^2\varphi|_{2,\Omega}^2\cr
&\quad+|\D(\nabla\varphi)|_{2,\Omega}^2+|\nabla\varphi\cdot\bar\tau_\alpha|_{2,S_1}^2+ |f|_{6/5,\Omega}^2)\cr
&\le c\intop_{S_2}\varrho_1d_1w_{n+1}^2dS_2+c(|\nabla\varphi_t|_{6/5,\Omega}^2+ |v_n|_{3,\Omega}^2\|\nabla\varphi\|_{1,\Omega}^2\cr
&\quad+\|\nabla\varphi\|_{1,\Omega}^2+|f|_{6/5,\Omega}^2).\cr}
\label{5.3}
\end{equation}
Using in (\ref{5.3}) the interpolation
$$\eqal{
&\intop_{S_2}\varrho_1d_1w_{n+1}^2dS_2\le|\varrho_1|_{3,S_2}|d_1|_{\infty,S_2} |w_{n+1}|_{3,S_2}^2\cr
&\le\varepsilon|\nabla w_{n+1}|_{2,\Omega}^2+c(1/\varepsilon) |\varrho_1|_{3,S_2}^6|d_1|_{\infty,S_2}^6|w_{n+1}|_{2,\Omega}^2\cr}
$$
 and integrating the result with respect to time, yield
\begin{equation}\eqal{
&|w_{n+1}(t)|_{2,\Omega}^2+\|w_{n+1}\|_{1,2,\Omega^t}^2+ |w_{n+1}\cdot\bar\tau_\alpha|_{2,S_1^t}^2\cr
&\le\exp(|d_1|_{\infty,S_2^t(-a)}^6|\varrho_1|_{3,6,S_2^t(-a)}^6)[(1+ |v_n|_{3,4,\Omega^t}^2)\cdot\|\nabla\varphi\|_{W_2^{2,1}(\Omega^t)}^2\cr
&\quad+|f|_{6/5,2,\Omega^t}^2+|w(0)|_{2,\Omega}^2].\cr}
\label{5.4}
\end{equation}
This inequality implies (\ref{5.1}) and concludes the proof of Lemma~\ref{5.1}.
\end{proof}

\begin{lemma}\label{l5.2}
Assume that $p_{n+1}\in W_{\sigma}^{{s},{s}/2}(\Omega^t)$, $\varrho_n\in W_{r,\infty}^{1,1}(\Omega^t)$, $v_{n+1}\in L_\infty(0,t;L_2(\Omega))$, $v_n\in V_{\sigma}^{2+{s}}(\Omega^t)$, $f\in W_{\sigma}^{{s},{s}/2}(\Omega^t)$, $r>{\sigma}$, $\frac{5}{r}<{s}$, $\frac{3}{\sigma}-\frac{3}{r}<s$, ${s}\in(0,1)$, $d\in W_{\sigma}^{2+{s}-1/{\sigma},1+{s}/2-1/2{\sigma}}(S_2^t)$, $v(0)\in W_{\sigma}^{2+{s}-2/{\sigma}}(\Omega)$. Then there exists $\bar{a}>0$ such that for solutions to (\ref{3.3})
\begin{equation}\eqal{
&\|v_{n+1}\|_{V_{\sigma}^{2+{s}}(\Omega^t)}+\|\nabla p_{n+1}\|_{W_{\sigma}^{{s},{s}/2}(\Omega^t)}\le c[\|p_{n+1}\|_{W_{\sigma} ^{{s},{s}/2}(\Omega^t)}\cr
&\quad+\phi_2(t^{\bar{a}}\|\varrho_n\|_{W_{r,\infty}^{1,1}(\Omega^t)},t^{\bar{a}} \|v_n\|_{V_{\sigma}^{2+{s}}(\Omega^t)})|v_{n+1}|_{2,\infty,\Omega^t}\cr
&\quad+|v_{n+1}|_{2,\Omega^t}+(1+t^{\bar{a}}\|\varrho_n\|_{W_{r,\infty}^{1,1}(\Omega^t)}) \|f\|_{W_{\sigma}^{{s},{s}/2}(\Omega^t)}\cr
&\quad+\|d\|_{W_{\sigma}^{2+{s}-1/{\sigma},1+{s}/2-1/2{\sigma}}(S_2^t)}+ \|v(0)\|_{W_{\sigma}^{2+{s}-2/{\sigma}}(\Omega)}].\cr}
\label{5.5}
\end{equation}
\end{lemma}

\begin{proof}
We consider the linear Stokes system (\ref{3.3}) and we localize the problem, following \cite{RZ3}. Hence, we introduce the partition of unity
$$
\zeta^{(k,l)}(x,t)=\zeta^{(k)}(x)\zeta_0^{(l)}(t),\quad k,l\in\N,
$$
such that
$$\eqal{
&\sup_k\diam\supp\zeta^{(k)}(x)\le\lambda,\cr
&\sup_l\diam\supp\zeta_0^{(l)}(t)\le\lambda,\cr}
$$
where $\lambda$ will be chosen later.

Let $\xi^{(k)}$, $\xi_0^{(l)}$ be interior points of $\supp\zeta^{(k)}$ and $\supp\zeta_0^{(l)}$, respectively. If $\overline{\supp\zeta^{(k)}}\cap S\not=\phi$ then $\xi^{(k)}$ is an interior point of $\overline{\supp\zeta^{(k)}}\cap S$.

Let
$$\eqal{
&\tilde v_{n+1}^{(k,l)}=v_{n+1}\zeta^{(k,l)}, \cr &\tilde p_{n+1}^{(k,l)}=p_{n+1}\zeta^{(k,l)},\cr
&\tilde f^{(k,l)}=f\zeta^{(k,l)}.\cr}
$$
Then problem (\ref{3.3}) takes the form
\begin{equation}\eqal{
&\varrho_n(\xi^{(k)},\xi_0^{(l)})\tilde v_{n+1,i,t}^{(k,l)}-(\divv\T(\tilde v_{n+1}^{(k,l)},\tilde p_{n+1}^{(k,l)}))_i\cr
&=[\varrho_n(\xi^{(k)},\xi_0^{(l)})-\varrho_n(x,t)]\tilde v_{n+1,i,t}^{(k,l)}+\varrho_nv_{n+1,i}\zeta_t^{(k,l)}\cr
&\quad+p_{n+1}\nabla_i \zeta^{(k,l)}-\nu\partial_{x_j}(v_{n+1,j}\zeta_{x_i}^{(k,l)}+ v_{n+1,i}\zeta_{x_j}^{(k,l)})+\varrho_n\tilde f_i^{(k,l)}\cr
&\quad-\nu(v_{n+1,i,x_j}+v_{n+1,j,x_i})\zeta_{x_j}^{(k,l)}+\varrho_nv_n\cdot \nabla v_{n+1,i}\zeta^{(k,l)} &\  {\rm in}\ \Omega^T,\cr
&\divv\tilde v_{n+1}^{(k,l)}=v_{n+1}\cdot\nabla\zeta^{(k,l)}\quad &{\rm in}\ \  \Omega, \cr
&\tilde v_{n+1}^{(k,l)}\cdot\bar n=0 \quad &{\rm on}\ \ S_1,\cr  & \nu\bar n\cdot\D(\tilde v_{n+1}^{(k,l)})\cdot\bar\tau_\alpha
=\nu n_i(v_{n+1,i}\zeta_{x_j}^{(k,l)}+v_{n+1,j}\zeta_{x_i}^{(k,l)})\cdot \tau_{\alpha j} \cr & \quad -\gamma\tilde v_{n+1}^{(k,l)}\cdot\bar\tau_\alpha\quad &{\rm on}\ \ S_1,\cr
&\tilde v_{n+1}^{(k,l)}\cdot\bar n=\tilde d^{(k,l)}\quad &{\rm on}\ \ S_2,\cr
&\bar n\cdot\D(\tilde v_{n+1}^{(k,l)})\cdot\bar\tau_\alpha=n_i(v_{n+1,i}\zeta_{x_j}^{(k,l)}+ v_{n+1,j}\zeta_{x_i}^{(k,l)})\cdot\tau_{\alpha j}\quad &{\rm on}\ \ S_2,\cr
&\tilde v_{n+1}^{(k,l)}|_{t=0}=\tilde v^{(k,l)}(0)\quad & {\rm in}\  \Omega,\cr}
\label{5.6}
\end{equation}
where $i=1,2,3$.

The r.h.s of this system needs some analysis. We examine the coefficient $[\varrho_n(\xi^{(k)},\xi_0^{(l)})-\varrho_n(x,t)]$ of the first term from the r.h.s. of (\ref{5.6}). It can be estimated in the following way.
$$\eqal{
&\sup_{k,l}\supp_{\zeta(k,l)}|\varrho_n(x,t)-\varrho_n(\xi^{(k)},\xi_0^{(l)})|\cr
&\le\sup_{k,l}\supp_{\zeta(k,l)}\bigg[\frac{|\varrho_n(x,t)-\varrho_n(\xi^{(k)},t)|}{|x-\xi^{(k)}|^\alpha}|x-\xi^{(k)}|^\alpha\cr
&\quad+\frac{|\varrho_n(\xi^{(k)},t)-\varrho_n(\xi^{(k)},\xi_0^{(l)})|}{|t-\xi_0^{(l)}|^\alpha}|t-\xi_0^{(l)}|^\alpha\bigg]\cr
&\le c\|\varrho_n\|_{\dot C^\alpha(\Omega^t)}\lambda^\alpha\le\delta,\cr}
$$
where in view of results from Section \ref{s4} we use the imbedding
$$
\|\varrho_n\|_{C^\alpha(\Omega^t)}\le c\|\varrho_n\|_{W_{r,\infty}^{1,1}(\Omega^t)}
$$
which holds under the condition $3/r+\alpha<1$.

Next, we consider the norm of the first term on the r.h.s. of (\ref{5.6}).
$$\eqal{
&\|(\varrho_n(x,t)-\varrho_n(\xi^{(k)},\xi_0^{(l)}))\tilde v_{n+1,t}^{(k,l)}\|_{W_{\sigma}^{{s},{s}/2}(\Omega^t)}\cr
&\le\|\varrho_n\|_{\dot C^\alpha(\Omega^t)}\lambda^\alpha\|\tilde v\|_{n+1,t}^{(k,l)}\|_{W_{\sigma}^{{s},{s}/2}(\Omega^t)}\cr
&\quad+\bigg(\intop_0^tdt\intop_\Omega\intop_{n+1,t} \frac{|\varrho_n(x',t)-\varrho_n(x'',t)|^{\sigma}}{ |x'-x''|^{3+{\sigma}{s}}}|\tilde v_{n+1,t}^{(k,l)}(x',t)|^{\sigma} dx'dx''\bigg)^{1/{\sigma}}\cr
&\quad+\bigg(\intop_\Omega dx\intop_0^t\intop_0^t\frac{|\varrho_n(x,t')-\varrho_n(x,t'')|^\sigma}{|t'-t''|^{1+{\sigma}{s}/2}} |\tilde v_{n+1,t}^{(k,l)}(x,t')|^{\sigma} dt'dt''\bigg)^{1/{\sigma}}\cr
&\equiv I_1+I_2+I_3,\cr}
$$
where
$$\eqal{
I_2&\le\bigg(\intop_0^tdt'\bigg(\intop_\Omega\intop_\Omega \frac{|\varrho_n(x',t')-\varrho_n(x'',t')|^{\sigma}\lambda_1}{ |x'-x''|^{3+{\sigma} \lambda_1[\frac{1}{ {\sigma}\lambda_1}(\frac{3}{2}\lambda_1-3)+s]}} dx'dx''\bigg)^{1/\lambda_1}\cr
&\quad\cdot\bigg(\intop_\Omega\intop_\Omega\frac{|\tilde v_{n+1,t}^{(k,l)}(x',t')|^{{\sigma} \lambda_2}}{|x'-x''|^{\frac{3}{2}\lambda_2}}dx'dx''\bigg)^{1/\lambda_2} \bigg)^{1/{\sigma}}\equiv I_2^1,\cr}
$$
where $1/\lambda_1+1/\lambda_2=1$, $\lambda_2<2$. Then we can perform integration with respect to $x''$ in the second integral.

Let ${s}'=\frac{1}{\sigma\lambda_1}\big(\frac{3}{2}\lambda_1-3\big)+{s}$. Then
$$\eqal{
I_2^1&\le\sup_t\|\varrho_n\|_{W_{{\sigma}\lambda_1}^{{s}'}(\Omega)}\bigg(\intop_0^tdt' \bigg(\intop_\Omega|\tilde v_{n+1,t}^{(k,l)}(x',t')|^{{\sigma} \lambda_2}dx'\bigg)^{1/\lambda_2}\bigg)^{1/{\sigma}}\cr
&\le ct^{\bar{a}}\sup_t\|\varrho_n\|_{W_{{\sigma}\lambda_1}^{{s}'}(\Omega)} \|\tilde v_{n+1,t}^{(k,l)}\|_{L_{{\sigma}\lambda_2}(\Omega^t)}\equiv I_2^2.\cr}
$$
Using the interpolation
\begin{equation}
\|\tilde v_{n+1,t}^{(k,l)}\|_{L_{{\sigma}\lambda_2}(\Omega^t)}\le c\|\tilde v_{n+1}^{(k,l)}\|_{W_{\sigma}^{2+{s},1+{s}/2}(\Omega^t)}^\theta|\tilde v_{n+1}^{(k,l)}|_{2,\Omega^t}^{1-\theta},
\label{5.7}
\end{equation}
where $\theta$ is a solution to the equation $\frac{5}{ {\sigma}\lambda_2}-2=(1-\theta)\frac{5}{2}+\theta\big(\frac{5}{ {\sigma}}-2-{s}\big)$ and $\theta<1$ implies the restriction $5/r<{s}$, where we used that ${\sigma}\lambda_1=r$. The last condition implies that $r>{\sigma}$. Then we have the bound
$$
I_2^2\le\varepsilon\|\tilde v_{n+1}^{(k,l)}\|_{W_{\sigma}^{2+{s},1+{s}/2}(\Omega^t)}+c(1/\varepsilon) (t^{\bar{a}}\sup_t\|\varrho_n\|_{W_r^1(\Omega)})^{1/(1-\theta)}\cdot|\tilde v_{n+1,t}^{(k,l)}|_{2,\Omega^t},
$$
where we used that ${s}'\le 1$.

Next, we consider $I_3$,
$$
I_3\le\bigg(\intop_\Omega dx\intop_0^tdt'\intop_0^tdt''\frac{|t'-t''|^{{\sigma}-1}\intop_0^t|\varrho_{n,t}|^{\sigma} dt}{ |t'-t''|^{1+{\sigma} {s}/2}}|\tilde v_{n+1,t}^{(k,l)}(x,t')|^{\sigma}\bigg)^{1/{\sigma}}\equiv I_3^1.
$$
Let ${\sigma}-1-1-{\sigma}{s}/2>-1$ so ${\sigma}(1-{s}/2)>1$. Then
$$
I_3^1\le t^{1-\frac{1}{{\sigma}}-\frac{{s}}{2}}\bigg(\intop dx\intop_0^tdt'\intop_0^t|\varrho_{n,t}(t'')|^{\sigma} dt''|\tilde v_{n+1,t}^{(k,l)}(x,t')|^{\sigma}\bigg)^{1/{\sigma}}\equiv I_3^2.
$$
By the H\"older inequality, we have
$$
I_3^2\le t^{1-{s}/2}\bigg(\intop_\Omega\bigg|\intop_0^t |\varrho_{n,t}|^{\sigma} dt''\bigg|^{\lambda_1}dx\bigg)^{1/{\sigma} \lambda_1}\bigg(\intop_\Omega dx\bigg|\intop_0^t|\tilde v_{n+1,t}^{(k,l)}|^{\sigma} dt\bigg|^{\lambda_2}\bigg)^{1/{\sigma}\lambda_2}\equiv I_3^3,
$$
where $1/\lambda_1+1/\lambda_2=1$. Using interpolation (\ref{5.7}),
$$
I_3^2\le\varepsilon\|\tilde v_{n+1}^{(k,l)}\|_{W_{\sigma}^{2+{s},1+{s}/2}(\Omega^t)}+c(1/\varepsilon) (t^{1-{s}/2}\sup_t|\varrho_{n,t}|_{r,\Omega})^{1/(1-\theta)}|\tilde v_{n+1}^{(k,l)}|_{2,\Omega^t},
$$
where $\theta$ is the same as in (\ref{5.7}).

Using \cite{RZ3} and that $\delta$ is sufficiently small we obtain from (\ref{5.6}), after summing up over all subdomains of the partition of unity, the inequality
\begin{equation}\eqal{
&\sup_t\|v_{n+1}\|_{W_{\sigma}^{2+{s}-2/{\sigma}}(\Omega)}+ \|v_{n+1}\|_{W_{\sigma}^{2+{s},1+{s}/2}(\Omega^t)}+\|\nabla p_{n+1}\|_{W_{\sigma} ^{{s},{s}/2}(\Omega^t)}\cr
&\le\phi(t^{\bar{a}}\|\varrho_n\|_{W_{r,\infty}^{1,1}(\Omega^t)}) [\|v_{n+1}\|_{W_{\sigma}^{{s},{s}/2}(\Omega^t)}+ \|p_{n+1}\|_{W_{\sigma} ^{{s},{s}/2}(\Omega^t)}\cr
&\quad+\|\nabla v_{n+1}\|_{W_{\sigma}^{{s},{s}/2}(\Omega^t)}+\|\varrho_nv_{n+1}\|_{W_{\sigma}^{{s},{s}/2}(\Omega^t)}+\|\!\varrho_nv_n\!\cdot\!\nabla v_{n+1}\!\|_{W_{\sigma}^{{s},{s}/2}(\Omega^t)} \cr
&\quad +\|\varrho_nf\|_{W_{\sigma}^{{s},{s}/2}(\Omega^t)}+ \|v_{n+1}\|_{W_{\sigma}^{1+{s},1/2+{s}/2}(\Omega^t)}\cr
&\quad+\|v_{n+1}\|_{W_{\sigma}^{1+{s}-1/{\sigma},1/2+{s}/2-1/2{\sigma}}(S^t)}+ \|d\|_{W_{\sigma}^{2+{s}-1/{\sigma},1+{s}/2-1/2{\sigma}}(S^t)}\cr
&\quad+\|v(0)\|_{W_{\sigma}^{2+{s}-2/{\sigma}}(\Omega)}]+\phi(t^{\bar{a}}\|\varrho_n\|_{\dot W_{r,\infty}^{1,1}(\Omega^t)})|v_{n+1}|_{2,\Omega^t},\cr}
\label{5.8}
\end{equation}
where $\bar{a}>0$ and
$$
\|\varrho_n\|_{\dot W_{r,\infty}^{1,1}(\Omega^t)}=\sup_t|\nabla\varrho_n|_{r,\Omega}+ \sup_t|\varrho_{n,t}|_{r,\Omega}.
$$
By interpolation (\ref{2.7}), from (\ref{5.8}) we infer
\begin{equation}\eqal{
&\sup_t\|v_{n+1}\|_{W_{\sigma}^{2+{s}-2/{\sigma}}(\Omega)}+ \|v_{n+1}\|_{W_{\sigma}^{2+{s},1+{s}/2}(\Omega^t)}+\|\nabla p_{n+1}\|_{W_{\sigma} ^{{s},{s}/2}(\Omega^t)}\cr
&\le\phi(t^{\bar{a}}\|\varrho_n\|_{W_{r,\infty}^{1,1}(\Omega^t)}) [\|p_{n+1}\|_{W_{\sigma}^{{s},{s}/2}(\Omega^t)}+ \|\varrho_nv_{n+1}\|_{W_{\sigma} ^{{s},{s}/2}(\Omega^t)}\cr
&\quad+\|\varrho_nv_n\cdot\nabla v_{n+1}\|_{W_{\sigma}^{{s},{s}/2}(\Omega^t)}\cr
&\quad+\|\varrho_nf\|_{W_{\sigma}^{{s},{s}/2}(\Omega^t)}+ \|d\|_{W_{\sigma}^{2+{s}-1/{\sigma},1+{s}/2-1/2{\sigma}}(S^t)}\cr
&\quad+\|v(0)\|_{W_{\sigma}^{2+{s}-2/{\sigma}}(\Omega)}+\phi(t^{\bar{a}} \|\varrho_n\|_{W_{r,\infty}^{1,1}(\Omega^t)})|v_{n+1}|_{2,\Omega^t}].\cr}
\label{5.9}
\end{equation}
In the following subsections, we analyze r.h.s. terms $\|\varrho_nv_{n+1}\|_{W_{\sigma}^{{s},{s}/2}(\Omega^t)},$ $\|\varrho_nv_n\cdot\nabla v_{n+1}\|_{W_{\sigma}^{{s},{s}/2}(\Omega^t)}$ and $\|\varrho_nf\|_{W_{\sigma}^{{s},{s}/2}(\Omega^t)}$ to come up with inequalities (\ref{5.19}), (\ref{5.30}) and  (\ref{5.37}).

\subsection{Analysis of inequality (\ref{5.9}): term $\|\varrho_nv_{n+1}\|_{W_{\sigma}^{{s},{s}/2}(\Omega^t)}$}

Consider the second term under the square bracket on the r.h.s. of (\ref{5.9}),
\begin{equation}\eqal{
&\|\varrho_nv_{n+1}\|_{W_{\sigma}^{{s},{s}/2}(\Omega^t)}=|\varrho_nv_{n+1}|_{{\sigma},\Omega^t}\cr
&\quad+\|\varrho_nv_{n+1}\|_{L_{\sigma}(0,t;W_{\sigma}^{s}(\Omega))}+ \|\varrho_nv_{n+1}\|_{L_{\sigma}(\Omega;W_{\sigma}^{{s}/2}(0,t))}\cr
&\equiv I_1+I_2+I_3,\cr}
\label{5.10}
\end{equation}
By interpolation (\ref{2.7}) and Lemma \ref{l2.4} we have
$$
I_1\le\varrho^*|v_{n+1}|_{{\sigma},\Omega^t}\le\varepsilon_1 \|v_{n+1}\|_{W_{\sigma}^{2+{s},1+{s}/2}(\Omega^t)}+c(1/\varepsilon_1) |v_{n+1}|_{2,\Omega^t}.
$$
Then we analyze $I_2$,
$$\eqal{
I_2&=\bigg(\intop_0^tdt'\intop_\Omega\intop_\Omega\frac{|\varrho_n(x',t')v_{n+1}(x',t')- \varrho_n(x'',t'')v_{n+1}(x'',t')|^{\sigma}}{|x'-x''|^{3+{\sigma}{s}}}dx'dx''\bigg)^{1/{\sigma}}\cr
&\le\bigg(\intop_0^tdt'\intop_\Omega\intop_\Omega \frac{|\varrho_n(x',t')-\varrho_n(x'',t')|^{\sigma}}{|x'-x''|^{3+{\sigma}{s}}} |v_{n+1}(x',t)|^{\sigma} dx'dx''\bigg)^{1/{\sigma}}\cr
&\quad+\bigg(\intop_0^tdt'\intop_\Omega \frac{|\varrho_n(x'',t')|^{\sigma}|v_{n+1}(x',t')-v_{n+1}(x'',t')|^{\sigma}}{|x'-x''|^{3+{\sigma}{s}}} dx'dx''\bigg)^{1/{\sigma}}\cr
&\equiv I_2^1+I_2^2.\cr}
$$
By the H\"older inequality,
$$\eqal{
I_2^1&\le\bigg(\intop_0^tdt'\bigg(\intop_\Omega\intop_\Omega \frac{|\varrho_n(x',t')-\varrho_n(x'',t')|^{{\sigma}\lambda_1}}{ |x'-x''|^{3+{\sigma} \lambda_1[\frac{1}{ {\sigma}\lambda_1}(\frac{3}{2}\lambda_1-3)+{s}]}} dx'dx''\bigg)^{1/\lambda_1}\cdot\cr
&\quad\cdot\bigg(\intop_\Omega\intop_\Omega \frac{|v_{n+1}(x',t')|^{{\sigma}\lambda_2}}{ |x'-x''|^{\frac{3}{2}\lambda_2}}dx'dx''\bigg)^{1/\lambda_2}\bigg)^{1/{\sigma}}\equiv I_2^{11},\cr}
$$
where $1/\lambda_1+1/\lambda_2=1$. For $\lambda_2<2$ we can integrate with respect to $x''$ in the second factor of $I_2^{11}$.

Let ${s}'=\frac{1}{\sigma\lambda_1}\big(\frac{3}{2}\lambda_1-3\big)+{s}$. Then
$$
I_2^{11}\le c\bigg(\intop_0^tdt'\|\varrho_n(t')\|^{\sigma}_{W_{{\sigma}\lambda_1}^{{s}'}(\Omega)} \|v_{n+1}(t')\|_{L_{{\sigma}\lambda_2}(\Omega)}^{\sigma} \bigg)^{1/{\sigma}}\equiv I_2^{12}.
$$
Assuming
\begin{equation}
\frac{3}{r}-\frac{3}{\sigma\lambda_1}+s'\le 1\quad {\rm or}\quad \frac{3}{r}-\frac{3}{{\sigma}\lambda_1}+\frac{3}{ 2{\sigma}}-\frac{3}{\sigma\lambda_1}+s\le 1,
\label{5.11}
\end{equation}
we can use the imbedding and next the estimate (\ref{4.23}),
\begin{equation}\eqal{
&\|\varrho_n\|_{W_{{\sigma}\lambda_1}^{{s}'}(\Omega)}\le c\|\varrho_n\|_{W_r^1(\Omega)} \cr
&\le c\bigg[\frac{1}{d_*(n)}(1+\|v_n\|_{V_{\sigma}^{2+{s}}(\Omega^t)})+1\bigg]\exp(\phi_1(t^{\bar{a}} \|v_n\|_{V_{\sigma}^{2+{s}}(\Omega^t)}))\bar d_1.\cr}
\label{5.12}
\end{equation}
We use the interpolation
\begin{equation}
\|v_{n+1}\|_{L_{{\sigma}\lambda_2}(\Omega)}\le\varepsilon \|v_{n+1}\|_{W_{\sigma}^{2+{s}-2/{\sigma}}(\Omega)}+c(1/\varepsilon)|v_{n+1}|_{2,\Omega},
\label{5.13}
\end{equation}
which holds for $\theta_1$ such that
\begin{equation}
\frac{3}{ {\sigma}\lambda_2}=(1-\theta_1)\frac{3}{2}+\theta_1\bigg(\frac{3}{\sigma}-\left(2+{s}-\frac{2}{\sigma}\right)\bigg)
\label{5.14}
\end{equation}
so
$$
\theta_1=\frac{\frac{3}{2}-\frac{3}{\sigma\lambda_2}}{2+s+\frac{3}{2}-\frac{5}{\sigma}}
$$
and the condition $\theta_1<1$ implies
\begin{equation}
\frac{5}{\sigma}-\frac{3}{\sigma\lambda_2}<2+{s}.
\label{5.15}
\end{equation}
Conditions (\ref{5.11}) and (\ref{5.15}) do not imply one restriction in the case ${\sigma}<r$, $3<{\sigma}$.

Then we obtain
$$\eqal{
I_2^{12}&\le\varepsilon\|v_{n+1}\|_{W_{\sigma}^{2+{s},1+{s}/2}(\Omega^t)}+c (1/\varepsilon)t^{1/2}\exp(c\phi_1(t^{\bar{a}}\|v_n\|_{V_{\sigma} ^{2+{s}}(\Omega^t)}))\cdot\cr
&\quad\cdot\bigg\{\bigg[\frac{1}{d_*(n)}(1+\|v_n\|_{V_{\sigma}^{2+{s}}(\Omega^t)})+1\bigg]\bar d_1\bigg\}^{a_1} |v_{n+1}|_{2,\infty,\Omega^t},\cr}
$$
where $a_1=\frac{1}{1-\theta_1}$ and $\bar{a}$ appears in (\ref{4.23}).

In view of Lemma \ref{l2.4} and interpolation (\ref{2.7}) we have
$$\eqal{
I_2^2&\le\varrho^*\|v_{n+1}\|_{L_{\sigma}(0,t;W_{\sigma}^{s}(\Omega))}\cr
&\le\varepsilon\|v_{n+1}\|_{W_{\sigma}^{2+{s},1+{s}/2}(\Omega^t)}+c (1/\varepsilon,\varrho^*)|v_{n+1}|_{2,\Omega^t}.\cr}
$$
Estimates of $I_2^1$ and $I_2^2$ imply
\begin{equation}\eqal{
&\|\varrho_nv_n\|_{L_{\sigma}(0,t;W_{\sigma}^{s}(\Omega))}\le\varepsilon \|v_{n+1}\|_{W_{\sigma}^{2+{s},1+{s}/2}(\Omega^t)}\cr
&\quad+c(1/\varepsilon)t^{1/2}\exp(c\phi_1(t^{\bar{a}}\|v_n\|_{V_{\sigma}^{2+{s}}(\Omega^t)}))\cdot\cr
&\quad\cdot\bigg\{\bigg[\frac{1}{d_*}(1+\|v_n\|_{V_{\sigma}^{2+{s}}(\Omega^t)})+1\bigg]\bar d_1\bigg\}^{a_1}|v_{n+1}|_{2,\infty,\Omega^t}\cr
&\quad+c(1/\varepsilon)|v_{n+1}|_{2,\Omega^t},\cr}
\label{5.16}
\end{equation}
where $\bar d_1$ is defined in (\ref{4.21}), $a_1=\frac{1}{1-\theta_1}$ and $\bar{a}>0$ appears in (\ref{4.23}).

Since $\lambda_1>2$ and $\lambda_2<2$ inequalities (\ref{5.11}) and (\ref{5.15}) imply the restrictions
\begin{equation}
\frac{3}{r}+\frac{11}{2\sigma}<5+{s}.
\label{5.17}
\end{equation}
Finally, we estimate the last term on the r.h.s. of (\ref{5.10}),
$$\eqal{
I_3&=\bigg(\intop_\Omega dx\intop_0^t\intop_0^t \frac{|\varrho_n(x,t')v_{n+1}(x,t')-\varrho_n(x,t'')v_{n+1}(x,t'')|^{\sigma} }{|t'-t''|^{1+{\sigma} {s}/2}}dt'dt''\bigg)^{1/{\sigma}}\cr
&\le\bigg(\intop_\Omega dx\intop_0^t\intop_0^t \frac{|\varrho_n(x,t')-\varrho_n(x,t'')|^{\sigma}}{ |t'-t''|^{1+{\sigma}{s}/2}}|v_{n+1}(x,t')|^{\sigma} dt'dt''\bigg)^{1/{\sigma}}\cr
&\quad+\bigg(\intop_\Omega dx\intop_0^t\intop_0^t|\varrho_n(x,t'')|^{\sigma} \frac{|v_{n+1}(x,t')-v_{n+1}(x,t'')|^{\sigma}}{|t'-t''|^{1+{\sigma} {s}/2}}dt'dt''\bigg)^{1/{\sigma}}\cr
&\equiv I_3^1+I_3^2.\cr}
$$
where
$$\eqal{
I_3^1&\le\bigg(\intop_\Omega dx\intop_0^t\intop_0^t\frac{|t'-t''|^{{\sigma}-1}}{|t'-t''|^{1+{\sigma}{s}/2}}\intop_{t'}^{t''} |\varrho_{n,t}(x,t)|^{\sigma} dt|v_{n+1}(x,t')|^{\sigma} dt'dt''\bigg)^{1/{\sigma}}\cr
&\le t^{1-1/{\sigma}-{s}/2}\bigg(\intop_\Omega dx\intop_0^t |\varrho_{n,t}(x,t')|^{\sigma} dt'\intop_0^t|v_{n+1}(x,t')|^{\sigma} dt'\bigg)^{1/{\sigma} }\cr
&\le ct^{1-\frac{1}{{\sigma}}-{s}/2}|\varrho_{n,t}|_{{\sigma}\lambda_1,\Omega^t} |v_{n+1}|_{{\sigma}\lambda_2,\Omega^t}\cr
&\le ct^{1-\frac{1}{\sigma}-{s}/2+1/r}\sup_t|\varrho_{n,t}|_{r,\Omega} |v_{n+1}|_{W_{\sigma}^{2+{s},1+{s}/2}(\Omega^t)}^{\theta_2} |v_{n+1}|_{2,\Omega^t}^{1-\theta_2}\equiv I_3^{11},\cr}
$$
where we assumed that ${\sigma}\lambda_1=r$ and $\theta_2$ satisfies the equation
\begin{equation}
\frac{5}{\sigma\lambda_2}=(1-\theta_2)\frac{5}{2}+\theta_2\bigg(\frac{5}{\sigma}-(2+{s})\bigg)\quad {\rm so}\quad \theta_2=\frac{\frac{5}{2}-\frac{5}{{\sigma}\lambda_2}}{2+{s}+\frac{5}{2}-\frac{5}{\sigma}}
\label{5.18}
\end{equation}
and ${\sigma}\lambda_2= \frac{r\sigma}{r-{\sigma}}$.

By the Young inequality, we obtain
$$
I_3^{11}\le\varepsilon|v_{n+1}|_{W_{\sigma}^{2+{s},1+{s}/2}(\Omega^t)}+ \bigg(\frac{1}{\varepsilon}t^{1-\frac{1}{\sigma}-\frac{s}{2}+\frac{1}{r}}\sup_t|\varrho_{n,t}|_{r,\Omega} \bigg)^{a_2}|v_{n+1}|_{2,\Omega^t},
$$
where $a_2=\frac{1}{1-\theta_2}$.

By interpolation (\ref{2.7}), the second term in $I_3$ is bounded by
$$\eqal{
I_3^2&\le\varrho^*\|v_{n+1}\|_{L_{\sigma}(\Omega;W_{\sigma}^{{s}/2}(0,t))}\cr
&\le\varepsilon\|v_{n+1}\|_{W_{\sigma}^{2+{s},1+{s}/2}(\Omega^t)}+c (1/\varepsilon,\varrho^*)|v_{n+1}|_{2,\Omega^t}.\cr}
$$
In view of bounds for $I_3^1$ and $I_3^2$ we have
\begin{equation}\eqal{
&\|\varrho_nv_{n+1}\|_{L_{\sigma}(\Omega;W_{\sigma}^{{s}/2}(0,t))}\le\varepsilon \|v_{n+1}\|_{W_{\sigma}^{2+{s},1+{s}/2}(\Omega^t)}\cr
&\quad+c(1/\varepsilon,\varrho^*)[(t^{1-1/{\sigma}-{s}/2+1/r}\sup_t |\varrho_{n,t}|_{r,\Omega})^{a_2}+1]|v_{n+1}|_{2,\Omega^t},\cr}
\label{5.19}
\end{equation}
where $a_2=\frac{1}{1-\theta_2}$, and $\theta_2$ is described by (\ref{5.18}). Since $\theta_2<1$, we have
\begin{equation}
r>{\sigma},\quad 1-\frac{1}{\sigma}-\frac{s}{2}+\frac{1}{r}>0.
\label{5.20}
\end{equation}

\subsection{Analysis of inequality (\ref{5.9}): term \\ $\|\varrho_nv_n\cdot\nabla v_{n+1}\|_{W_{\sigma}^{{s},{s}/2}(\Omega^t)}$}

Next we examine the third term on the r.h.s. of (\ref{5.9}). It is equal
$$\eqal{
&|\varrho_nv_n\cdot\nabla v_{n+1}|_{{\sigma},\Omega^t}+\|\varrho_nv_n\cdot\nabla v_{n+1}\|_{L_{\sigma}(0,t;W_{\sigma}^{s}(\Omega))}\cr
&\quad+\|\varrho_nv_n\cdot\nabla v_{n+1}\|_{L_{\sigma}(\Omega;W_{\sigma}^{{s}/2}(0,t))}\equiv I_1+I_2+I_3. \cr}
$$
Consider $I_1$. Then we have
$$\eqal{
|I_1|&\le\varrho^*|v_n\cdot\nabla v_{n+1}|_{{\sigma},\Omega^t}=\varrho^*\bigg(\intop_0^t|v_n\cdot\nabla v_{n+1}|_{{\sigma},\Omega}^{\sigma} dt'\bigg)^{1/{\sigma}}\cr
&\le\varrho^*\bigg(\intop_0^t|v_n|_{{\sigma}\lambda_1,\Omega}^{\sigma}|\nabla v_{n+1}|_{{\sigma}\lambda_2,\Omega}^{\sigma} dt'\bigg)^{1/{\sigma} }\equiv I_1^1,\cr}
$$
where $1/\lambda_1+1/\lambda_2=1$. Continuing,
$$
I_1^1\le\varrho^*\sup_t\|v_n\|_{W_{\sigma}^{2+{s}-2/{\sigma}}(\Omega)}\bigg(\intop_0^t |\nabla v_{n+1}|_{{\sigma}\lambda_2,\Omega}^{\sigma} dt'\bigg)^{1/{\sigma}}\equiv I_1^2,
$$
where we used the imbedding
$$
|v_n|_{{\sigma}\lambda_1,\Omega}\le c\|v_n\|_{W_{\sigma}^{2+{s}-2/{\sigma}}(\Omega)}\quad {\rm for}\quad \frac{3}{ {\sigma}}-\frac{3}{ {\sigma} \lambda_1}\le 2+{s}-\frac{2}{{\sigma}}.
$$
Next, we use the interpolation
$$
|\nabla v_{n+1}|_{{\sigma}\lambda_2,\Omega}\le c\|v_{n+1}\|_{W_{\sigma}^{2+{s}}(\Omega)}^{\theta_3}|v_{n+1}|_{2,\Omega}^{1-\theta_3},
$$
where $\theta_3$ is a solution to the equation
$$
\frac{3}{{\sigma}\lambda_2}-1=\left(1-\theta_3\right)\frac{3}{2}+\theta_3\bigg(\frac{3}{\sigma}-(2+{s})\bigg).
$$
To show the existence of such $\theta_3$ we calculate
$\theta_3 = \displaystyle{\frac{\frac{3}{2}+1-\frac{3}{{\sigma}\lambda_2}}{2+{s}+\frac{3}{2}-\frac{3}{{\sigma}}}}$. Then $\theta_3>0$ because $\frac{5}{2}-\frac{3}{{\sigma}\lambda_2}>0$ and $\theta_3<1$ implies $\frac{3}{2}- \frac{3}{{\sigma}\lambda_2}<1+{s}+\frac{3}{2}-\frac{3}{{\sigma} }$.

Hence
$$
\frac{3}{\sigma}-\frac{3}{\sigma\lambda_2}<1+{s}.
$$
Since we have two restrictions
$$\eqal{
&\frac{3}{ {\sigma}}-\frac{3}{{\sigma}\lambda_1}\le 2+{s}- \frac{2}{\sigma}, \cr
&\frac{3}{ {\sigma}}-\frac{3}{ {\sigma}\lambda_2}<1+{s}, \cr}
$$
so $\frac{6}{\sigma}+\frac{2}{\sigma}-\frac{3}{{\sigma}}<3+2{s}$. Therefore, we obtain the inequality $\frac{5}{{\sigma}}<3+2{s}$ which always holds.

Then
$$\eqal{
I_1^2&\le\varepsilon\|v_{n+1}\|_{W_{\sigma}^{2+{s},1+{s}/2}(\Omega^t)}\cr
&\quad+\bigg(c(1/\varepsilon) t^{1/{\sigma}}\sup_t\|v_n\|_{W_{\sigma}^{2+{s}-2/{\sigma}}(\Omega)}\bigg)^{\frac{1}{ 1-\theta_3}} \cdot|v_{n+1}|_{2,\infty,\Omega^t}.\cr}
$$
Next, we examine $I_2$,
$$\eqal{
&I_2=\bigg(\intop_0^t\intop_\Omega\intop_\Omega\bigg[ \frac{|\varrho_n(x',t')v_n(x',t')\cdot\nabla v_{n+1}(x',t')}{ |x'-x''|^{3+{\sigma}{s}}}\cr
&\hskip3,1cm-\frac{\varrho_n(x'',t')v_n(x'',t')\cdot\nabla v_{n+1}(x'',t')|^{\sigma}}{ |x'-x''|^{3+{\sigma}{s}}}\bigg]dx'dx''dt'\bigg)^{1/{\sigma}}\cr
&\le\bigg(\intop_0^t\intop_\Omega\!\!\intop_\Omega\!\!\frac{|\varrho_n(x',t')-\varrho_n(x'',t')|^{\sigma}}{ |x'-x''|^{3+{\sigma}{s}}}|v_n(x',t')|^{\sigma} |\nabla v_{n+1}(x',t')|^{\sigma} dx'dx''dt'\bigg)^{1/{\sigma}}\cr
&\quad+\varrho^*\bigg(\intop_0^t\intop_\Omega\intop_\Omega \frac{|v_n(x',t')-v_n(x'',t')|^{\sigma}|\nabla v_{n+1}(x',t')|^{\sigma}}{|x'-x''|^{3+{\sigma} {s}}}dx'dx''dt'\bigg)^{1/{\sigma}}\cr
&\quad+\varrho^*\bigg(\intop_0^t\intop_\Omega\intop_\Omega \frac{|v_n(x'',t')|^{\sigma}|\nabla v_{n+1}(x',t')-\nabla v_{n+1}(x'',t')|^{\sigma}}{|x'-x''|^{3+{\sigma}{s}}}dx'dx''dt'\bigg)^{1/{\sigma}}\cr
&\equiv I_2^1+I_2^2+I_2^3.\cr}
$$
We consider $I_2^1$,
$$\eqal{
& I_2^1\le\bigg[\intop_0^tdt'\bigg(\intop_\Omega\intop_\Omega \frac{|\varrho(x',t')-\varrho_n(x'',t')|^{\sigma\lambda_1}}{ |x'-x''|^{(3/2)\lambda_1+{\sigma}{s}\lambda_1}}dx'dx''\bigg)^{1/\lambda_1}\cr
&\bigg(\intop_\Omega\intop_\Omega \frac{|v_n(x',t')|^{\sigma\lambda_2}}{|x'-x''|^{(3/2)\lambda_2}} dx'dx''\bigg)^{1/\lambda_2}\bigg(\intop_\Omega\intop_\Omega|\nabla v_{n+1}|^{{\sigma}\lambda_3}dx'dx''\bigg)^{1/\lambda_3}\bigg]^{1/{\sigma}}\equiv I_2^{11},\cr}
$$
where $1/\lambda_1+1/\lambda_2+1/\lambda_3=1$. If $\lambda_2<2$ then
$$\eqal{
I_2^{11}&\le c\bigg(\intop_0^tdt\intop_\Omega\intop_\Omega \frac{|\varrho_n(x',t')-\varrho_n(x'',t')|^{{\sigma}\lambda_1}} {|x'-x''|^{3+{\sigma} \lambda_1[\frac{1}{ {\sigma}\lambda_1}(\frac{3}{2}\lambda_1-3)+{s}]}}dx'dx''\bigg)^{1/{\sigma}\lambda_1}\cdot\cr
&\quad\cdot\bigg(\intop_0^t\intop_\Omega|v_n(x,t')|^{{\sigma}\lambda_2}dxdt' \bigg)^{1/{\sigma}\lambda_2}\bigg(\intop_0^t\intop_\Omega|\nabla v_{n+1}(x,t')|^{{\sigma}\lambda_3}dxdt'\bigg)^{1/{\sigma}\lambda_3}\cr
&\equiv J_1J_2J_3.\cr}
$$
Assuming
$$
\frac{1}{ {\sigma}\lambda_1}\bigg(\frac{3}{2}\lambda_1-3\bigg)+{s}\le 1,
$$
we obtain
$$
J_1\le\sup_t\|\varrho_n\|_{W_{{\sigma}\lambda_1}^1(\Omega)}t^{1/{\sigma}\lambda_1}.
$$
Let $r={\sigma}\lambda_1$. Then the above condition reads
\begin{equation}
\frac{3}{ 2{\sigma}}-\frac{3}{r}+{s}\le 1.
\label{5.21}
\end{equation}
If
\begin{equation}
3/{\sigma}-3/{\sigma}\lambda_2\le 2+{s}-2/{\sigma},
\label{5.22}
\end{equation}
then
$$
J_2\le t^{1/{\sigma}\lambda_2}\sup_t\|v_n\|_{W_{\sigma}^{2+{s}-2/{\sigma}}(\Omega)}.
$$
Finally, we use the interpolation
$$
|\nabla v_{n+1}|_{{\sigma}\lambda_3,\Omega}\le c\|v_{n+1}\|_{W_{\sigma}^{2+{s}-2/{\sigma}}(\Omega)}^{\theta_4}|v_{n+1}|_{2,\Omega}^{1-\theta_4},
$$
where
$$
\frac{3}{{\sigma}\lambda_3}-1=(1-\theta_4)\frac{3}{2}+\theta_4\bigg(\frac{3}{\sigma}-\bigg(2+{s}-\frac{2}{ {\sigma}}\bigg)\bigg),\quad {\rm so}
$$
\begin{equation}
\theta_4=\frac{\frac{5}{2}-\frac{3}{{\sigma}\lambda_3}}{\frac{3}{2}-\frac{5}{{\sigma}}+2+{s}}.
\label{5.23}
\end{equation}
Hence $\theta_4<1$ implies the restriction
\begin{equation}
\frac{5}{\sigma}-\frac{3}{\sigma\lambda_3}<1+{s}.
\label{5.24}
\end{equation}
Taking into account (\ref{5.21}), (\ref{5.22}) and (\ref{5.24}) implies
$$
\frac{17}{ 2(4+s)}<\sigma,
$$
which is always satisfied because ${\sigma}>3$. Summarizing, we have
\begin{equation}\eqal{
&\|\varrho_nv_n\cdot\nabla v_{n+1}\|_{L_{\sigma}(0,t;W_{\sigma}^{s}(\Omega))}\le \varepsilon\sup_t\|v_{n+1}\|_{W_{\sigma}^{2+{s}-2/{\sigma} }(\Omega)}\cr
&\quad+(c(1/\varepsilon)t^{\bar{a}}\sup_t\|\varrho_n\|_{W_r^1(\Omega)}\sup_t \|v_n\|_{W_{\sigma}^{2+{s}-2/{\sigma} }(\Omega)})^{a_4}|v_{n+1}|_{2,\infty,\Omega^t},\cr}
\label{5.25}
\end{equation}
where $a_4=\frac{1}{1-\theta_4}$ and $\theta_4$ is defined by (\ref{5.23}).

Next, we calculate
$$\eqal{
I_3&=\bigg(\intop_\Omega dx\intop_0^t\intop_0^t\frac{|\varrho(x,t')v_n(x,t')\nabla v_{n+1}(x,t')}
{|t'-t''|^{1+{\sigma}{s}/2}}\cr
&\hskip5cm\frac{-\varrho_n(x,t'')v_n(x,t'')\nabla v_{n+1}(x,t'')|^{\sigma}}{ |t'-t''|^{1+{\sigma}{s}/2}}dt'dt''\bigg)^{1/{\sigma}}\cr
&\le\bigg(\intop_\Omega dx\intop_0^t\intop_0^t\frac{|\varrho_n(x,t')-\varrho_n(x,t'')|^{\sigma}|v_n(x,t')|^{\sigma}|\nabla v_{n+1}(x,t')|^{\sigma}}{|t'-t''|^{1+{\sigma}{s}/2}}dt'dt''\bigg)^{1/{\sigma}}\cr
&\quad+\varrho^*\bigg(\intop_\Omega dx\intop_0^t\intop_0^t \frac{|v_n(x,t')-v_n(x,t'')|^{\sigma}}{|t'-t''|^{1+{\sigma}{s}/2}} |\nabla v_{n+1}(x,t')|^{\sigma} dt'dt''\bigg)^{1/{\sigma}}\cr
&\quad+\varrho^*\bigg(\intop_\Omega dx\intop_0^t\intop_0^t\frac{|v_n(x,t'')|^{\sigma}|\nabla v_{n+1}(x,t')-\nabla v_{n+1}(x,t'')|^\sigma}{|t'-t''|^{1+{\sigma}{s}/2}}dt'dt''\bigg)^{1/{\sigma}} \cr
&\equiv I_3^1+I_3^2+I_3^3.\cr}
$$
First, we estimate
$$\eqal{
I_3^1\!&\!\le\bigg(\!\intop_\Omega\!\! dx\intop_0^t\!\intop_0^t\!{(t'-t'')^{{\sigma}-1}\!\!\intop_0^t|\varrho_{n,\tau}|^{\sigma}d\tau}{ |t'-t''|^{1+{\sigma s}/2}}|v_n(x,t')|^{\sigma} \cr & |\nabla v_{n+1}(x,t')|^{\sigma} dt'dt''\bigg)^{1/{\sigma}} 
\equiv I_3^{11}.\cr}
$$ 
If ${\sigma}-1-1-{\sigma}{s}/2>-1,$ so ${\sigma}(1-{s}/2)>1,$ we obtain
$$\eqal{
I_3^{11}&\le ct^{1-\frac{1}{\sigma}-\frac{{s}}{2}}\bigg(\intop_\Omega dx\intop_0^t|\varrho_{n,t}|^{\sigma} dt\intop_0^t|v_n(x,t')|^{\sigma} |\nabla v_{n+1}(x,t')|^{\sigma} dt'\bigg)^{1/{\sigma}}\cr
& \le ct^{1-{s}/2}\sup_t\|v_n\|_{W_{\sigma}^{2+{s}-2/{\sigma}}(\Omega)}\bigg|\intop_0^t |\varrho_{n,t}|^{\sigma} dt\bigg|_{\lambda_1,\Omega}^{1/{\sigma} }\cdot\bigg|\intop_0^t|\nabla v_{n+1}|^{\sigma} dt\bigg|_{\lambda_2,\Omega}^{1/{\sigma}}\cr
&\equiv I_3^{12},\cr}
$$
where $1/\lambda_1+1/\lambda_2=1$ and we need $\frac{5}{\sigma}<2+{s}$. Set ${\sigma}\lambda_1=r$ and apply the interpolation
\begin{equation}
|\nabla v_{n+1}|_{{\sigma}\lambda_2,{\sigma},\Omega^t}\le c\|v_{n+1}\|_{W_{\sigma}^{2+{s},1+{s}/2}(\Omega^t)}^{\theta_5} |v_{n+1}|_{2,\Omega^t}^{1-\theta_5},
\label{5.26}
\end{equation}
where $\theta_5$ satisfies
\begin{equation}
\frac{3}{\sigma\lambda_2}+\frac{2}{\sigma}-1=(1-\theta_5)\frac{5}{2}+\theta_5\bigg(\frac{5}{\sigma}-(2+{s})\bigg),
\label{5.27}
\end{equation}
so
$$
\theta_5=\frac{\frac{5}{2}+1-\frac{2}{\sigma}-\frac{3}{{\sigma}\lambda_2}}{ 2+{s}+\frac{5}{2}-\frac{5}{\sigma}}.
$$
The condition $\theta_5<1$ implies $\frac{3}{r}<1+{s}$.

Summarizing the above estimates yields
$$\eqal{
I_3^{12}&\le ct^{1+1/\sigma-{s}/2}\sup_t\|v_n\|_{W_{\sigma}^{2+{s}-2/{\sigma}}(\Omega)}\sup_t|\varrho_{n,t}|_{r,\Omega} \|v_{n+1}\|_{W_{\sigma}^{2,1}(\Omega^t)}^{\theta_5}\cdot|v_{n+1}|_{2,\Omega^t}^{1-\theta_5}\cr
&\le\varepsilon\|v_{n+1}\|_{W_{\sigma}^{2+{s},1+{s}/2}(\Omega^t)}\cr
&\quad+(c(1/\varepsilon) t^{1+1/{\sigma}-{s}/2}\sup_t\|v_n\|_{W_{\sigma}^{2+{s}-2/{\sigma}}(\Omega))} |\varrho_{n,t}|_{r,\infty,\Omega^t})^{a_5}|v_{n+1}|_{2,\Omega^t},\cr}
$$
where $a_5=1/(1-\theta_5)$.

Next, we examine
$$\eqal{
I_3^2&\le\varrho^*\bigg(\intop_\Omega dx\intop_0^t\intop_0^t \frac{|t'-t''|^{{\sigma}-1}\intop_{t'}^{t''}|v_{n,t}|^{\sigma}d\tau}{ |t'-t''|^{1+{\sigma} {s}/2}}|\nabla v_{n+1}(x,t')|^{\sigma} dt'dt''\bigg)^{1/{\sigma}}\cr
&\le\varrho^*t^{1-\frac{1}{\sigma}-\frac{s}{ 2}}\bigg(\intop_\Omega dx\intop_0^t|v_{n,t}|^{\sigma} dt\intop_0^t|\nabla v_{n+1}(x,t')|^{\sigma} dt'\bigg)^{1/{\sigma}}\cr
&\le c\varrho^*t^{1-\frac{1}{\sigma}-\frac{s}{2}}|v_{n,t}|_{{\sigma}\lambda_1,{\sigma},\Omega^t}|\nabla v_{n+1}|_{{\sigma}\lambda_2,{\sigma} ,\Omega^t}\equiv I_3^{21},\cr}
$$
where the last inequality follows from the Minkowski inequality.

By imbedding
$$
\|v_{n,t}\|_{L_{{\sigma}\lambda_1,{\sigma}}(\Omega^t)}\le c\|v_n\|_{W_{\sigma}^{2+{s},1+{s}/2}(\Omega^t)},
$$
which holds if the relation is satisfied
\begin{equation}
\frac{5}{\sigma}-\frac{3}{ {\sigma}\lambda_1}-\frac{2}{\sigma}\le{s},\quad {\rm so}\quad \frac{3}{ {\sigma}\lambda_2}\le{s}.
\label{5.28}
\end{equation}
In this case we use interpolation (\ref{5.26}) with $\theta_5$ defined by (\ref{5.27}).

Since (\ref{5.28}) holds the restriction $\theta_5<1$ is satisfied if $\frac{3}{ {\sigma}\lambda_1}=\frac{3}{ {\sigma}}-\frac{3}{ {\sigma} \lambda_2}\le1+{s}$. Combining the above restrictions yields $3/{\sigma}<1+2{s}$. Then
$$\eqal{
I_3^{21}&\le c\varrho^*t^{1-1/{\sigma}-{s}/2}\|v_n\|_{W_{\sigma}^{2+{s},1+{s}/2}(\Omega^t)} \|v_{n+1}\|_{W_{\sigma} ^{2+{s},1+{s}/2}(\Omega^t)}^{\theta_5} |v_{n+1}|_{2,\Omega^t}^{1-\theta_5}\cr
&\le\varepsilon\|v_{n+1}\|_{W_{\sigma}^{2+{s},1+{s}/2}(\Omega^t)}\cr
&\quad+(c(1/\varepsilon) t^{1-1/{\sigma}-{s}/2}\|v_n\|_{W_{\sigma}^{2+{s},1+{s}/2}(\Omega^t)})^{1/(1-\theta_5)} |v_{n+1}|_{2,\Omega^t}.\cr}
$$
Finally,
$$\eqal{
I_3^3&\le\varrho^*\sup_t|v_n|_{\infty,\Omega}\|\nabla v_{n+1}\|_{L_{\sigma}(\Omega;W_{\sigma}^{{s}/2}(0,t))}\cr
&\le\varrho^*\sup_t|v_n|_{\infty,\Omega}\|\nabla v_{n+1}\|_{W_{\sigma}^{{s},{s}/2}(\Omega^t)}\equiv I_3^{31}. \cr}
$$
We use the imbedding
$$
|v_n|_{\infty,\Omega}\le c\|v_n\|_{W_{\sigma}^{2+{s}-2/{\sigma}}(\Omega)},
$$
which holds for $3/{\sigma}<2+{s}-2/{\sigma}$ and the interpolation
$$
\|\nabla v_{n+1}\|_{W_{\sigma}^{{s},{s}/2}(\Omega^t)}\le c\|v_{n+1}\|_{W_{\sigma}^{2+{s},1+{s}/2}(\Omega^t)}^{\theta_6} |v_{n+1}|_{2,\Omega^t}^{1-\theta_6},
$$
where $\theta_6$ satisfies the relation
$$
\frac{5}{ {\sigma}}-1-{s}=(1-\theta_6)\frac{5}{2}+\theta_6\bigg(\frac{5}{ {\sigma}}-2-{s}\bigg),
$$
so $\theta_6$ is equal to
\begin{equation}
\theta_6=\frac{1+{s}+5/2-5/{\sigma}}{ 2+{s}+5/2-5/{\sigma}},
\label{5.29}
\end{equation}
where $\theta_6<1$ always holds.

Then, we have
$$\eqal{
I_3^{31}&\le c\varrho^*\sup_t\|v_n\|_{W_{\sigma}^{2+{s}-2/{\sigma}}(\Omega)} \|v_{n+1}\|_{W_{\sigma}^{2+{s},1+{s}/2}(\Omega^t)}^{\theta_6} |v_{n+1}|_{2,\Omega^t}^{1-\theta_6}\cr
&\le\varepsilon\|v_{n+1}\|_{W_{\sigma}^{2+{s},1+{s}/2}(\Omega^t)}\cr
&\quad+c(1/\varepsilon) (\varrho^*\sup_t\|v_n\|_{W_{\sigma}^{2+{s}-2/{\sigma}}(\Omega)})^{a_6}t^{1/2} |v_{n+1}|_{2,\infty,\Omega^t},\cr}
$$
where $a_6=1/(1-\theta_6)$.

Summarizing the above calculations yields
\begin{equation}\eqal{
&\|\varrho_nv_n\cdot\nabla v_{n+1}\|_{L_{\sigma}(\Omega;W_{\sigma}^{{s}/2}(0,t))}\le \varepsilon\|v_{n+1}\|_{W_{\sigma} ^{2+{s},1+{s}/2}(\Omega^t)}\cr
&\quad+c(1/\varepsilon[(t^{\bar{a}}\|v_n\|_{V_{\sigma}^{2+{s}}(\Omega^t)} |\varrho_{n,t}|_{r,\infty,\Omega^t})^{a_5}|v_{n+1}|_{2,\Omega^t}\cr
&\quad+(t^{\bar{a}}\|v_n\|_{V_{\sigma}^{2+{s}}(\Omega^t)})^{a_5}|v_{n+1}|_{2,\Omega^t}\cr
&\quad+(\varrho^*\|v_n\|_{V_{\sigma}^{2+{s}}(\Omega^t)})^{a_6}t^{1/2} |v_{n+1}|_{2,\infty,\Omega^t},\cr}
\label{5.30}
\end{equation}
where $a_5=1/(1-\theta_5)$ and $a_6=1/(1-\theta_6)$.

Moreover, $\theta_5$ is defined in (\ref{5.27}) and $\theta_6$ in (\ref{5.29}).

\subsection{Analysis of inequality (\ref{5.9}): term $\|\varrho_nf\|_{W_{\sigma}^{{s},{s}/2}(\Omega^t)}$}\label{s5.3}

Finally, we calculate
\begin{equation}\eqal{
&\|\varrho_nf\|_{W_{\sigma}^{{s},{s}/2}(\Omega^t)}\le\varrho^* \|f\|_{W_{\sigma}^{{s},{s}/2}(\Omega^t)}\cr
&\quad+\bigg(\intop_0^t\intop_\Omega\intop_\Omega \frac{|f(x',t)|^{\sigma}|\varrho_n(x',t)-\varrho_n(x'',t)|^{\sigma}}{|x'-x''|^{3+{s} {\sigma}}} dx'dx''dt\bigg)^{1/{\sigma}}\cr
&\quad+\bigg(\intop_\Omega\intop_0^t\intop_0^t\frac{|f(x,t')|^{\sigma}|\varrho_n(x,t')- \varrho_n(x,t'')|^{\sigma}}{|t'-t''|^{1+{\sigma} {s}/2}}dt'dt''dx\bigg)^{1/{\sigma}}\cr
&\equiv\varrho^*\|f\|_{W_{\sigma}^{{s},{s}/2}(\Omega^t)}+I_1+I_2.\cr}
\label{5.31}
\end{equation}
Consider $I_1$. Applying the H\"older inequality with respect to space variables we get
$$\eqal{
I_1&\le\bigg[\intop_0^t\bigg(\intop_\Omega\intop_\Omega \frac{|f(x',t)|^{{\sigma}\lambda_2}dx'dx''}{|x'-x''|^{(3/2)\lambda_2}}\bigg)^{1/\lambda_2} \cdot\cr
&\quad\cdot\bigg(\intop_\Omega\intop_\Omega \frac{|\varrho_n(x',t)-\varrho_n(x'',t)|^{{\sigma}\lambda_1}}{|x'-x''|^{3+{\sigma} \lambda_1{s}'}}dx'dx''\bigg)^{1/\lambda_1}dt\bigg]^{1/{\sigma}} \equiv I_1^1,\cr}
$$
where $1/\lambda_1+1/\lambda_2=1$, ${s}'=\frac{1}{ {\sigma}\lambda_1}\big(\frac{3}{2}\lambda_1-3+{s}\big)$. Assuming that $\lambda_2<2$ we can integrate with respect to $x''$ in the first factor under the time integral.

Then we obtain
$$
I_1^1\le\sup_t\|\varrho_n\|_{W_{{\sigma}\lambda_1}^{{s}'}(\Omega)} |f|_{{\sigma}\lambda_2,{\sigma},\Omega^t}\equiv I_1^2.
$$
Using the imbeddings
$$
\|\varrho_n\|_{W_{{\sigma}\lambda_1}^{{s}'}(\Omega)}\le c\|\varrho_n\|_{W_r^1(\Omega)}
$$
for
\begin{equation}
\frac{3}{r}-\frac{3}{\sigma\lambda_1}+{s}'\le 1\quad {\rm so}\quad \frac{3}{r}-\frac{3}{\sigma\lambda_1}+\frac{3}{ 2{\sigma}}-\frac{3}{\sigma\lambda_1}+{s}\le 1
\label{5.32}
\end{equation}
and
$$
|f|_{{\sigma}\lambda_2,{\sigma},\Omega^t}\le t^{\bar{a}}|f|_{{\sigma}\lambda_2,{\sigma}',\Omega^t}\le ct^{\bar{a}}\|f\|_{W_{\sigma}^{{s},{s}/2}(\Omega^t)}
$$
for ${\sigma}'>{\sigma}$ and
\begin{equation}
\frac{5}{\sigma}-\frac{3}{\sigma\lambda_2}-\frac{2}{\sigma}<{s}.
\label{5.33}
\end{equation}
Restrictions (\ref{5.32}) and (\ref{5.33}) imply
\begin{equation}
\frac{3}{r}+\frac{3}{2{\sigma}}<1+{s}.
\label{5.34}
\end{equation}
Hence
\begin{equation}
I_1\le I_1^2\le ct^{\bar{a}}\sup_t\|\varrho_n\|_{W_r^1(\Omega)}\|f\|_{W_{\sigma}^{{s},{s}/2}(\Omega^t)}.
\label{5.35}
\end{equation}
Consider $I_2$,
$$\eqal{
I_2&\le\bigg(\intop_\Omega\intop_0^t\intop_0^t \frac{|f(x,t')|^{\sigma}|t'-t''|^{{\sigma}-1}\intop_{t'}^{t''}|\varrho_{n,\tau}(\cdot,\tau)|^{\sigma} d\tau}{ |t'-t''|^{1+{\sigma}{s}/2}}dt'dt''dx\bigg)^{1/{\sigma}}\cr
&\le t^{1-{s}/2}\bigg(\intop_\Omega\intop_0^t|\varrho_{n,\tau}(\cdot,\tau)|^{\sigma}d\tau \intop_0^t|f(x,t)|^{\sigma} dtdx\bigg)^{1/{\sigma}}\cr
&\le ct^{1-{s}/2}\bigg(\intop_0^t|\varrho_{n,\tau}|_{{\sigma}\lambda_1,\Omega}^{\sigma}d\tau\bigg)^{1/{\sigma}} \bigg(\intop_0^t|f(\cdot,t)|_{{\sigma} ,\lambda_2,\Omega}^{\sigma} dt\bigg)^{1/{\sigma}}\equiv I_2^1.\cr}
$$
Setting ${\sigma}\lambda_1=r$, and using the imbedding
$$
|f|_{{\sigma}\lambda_2,{\sigma},\Omega^t}\le c\|f\|_{W_{\sigma}^{{s},{s}/2}(\Omega^t)},
$$
which holds for $\frac{3}{r}\le{s}$. Therefore
\begin{equation}
I_2\le t^{1+\frac{1}{ {\sigma}}-\frac{s}{2}}|\varrho_{n,\tau}|_{r,\infty,\Omega^t} \|f\|_{W_{\sigma}^{{s},{s}/2}(\Omega^t)}.
\label{5.36}
\end{equation}
Hence, applying (\ref{5.34}) and (\ref{5.35}) in (\ref{5.31}) yields
\begin{equation}
\|\varrho_nf\|_{W_{\sigma}^{{s},{s}/2}(\Omega^t)}\le c(\varrho^*+t^{\bar{a}}\|\varrho_n\|_{W_{r,\infty}^{1,1}(\Omega^t)}) \|f\|_{W_{\sigma} ^{{s},{s}/2}(\Omega^t)}.
\label{5.37}
\end{equation}
Using estimates (\ref{5.19}), (\ref{5.30}) and  (\ref{5.37}) in (\ref{5.9}) yields
\begin{equation}\eqal{
&\|v_{n+1}\|_{V_{\sigma}^{2+{s}}(\Omega^t)}+\|\nabla p_{n+1}\|_{W_{\sigma}^{{s},{s}/2}(\Omega^t)}\cr
&\le c[\|p_{n+1}\|_{W_{\sigma}^{{s},{s}/2}(\Omega^t)}+\phi_2(t^{\bar{a}} \|\varrho_n\|_{\dot W_{r,\infty}^{1,1}(\Omega^t)},\cr
&t^{\bar{a}}\|v_n\|_{V_{\sigma}^{2+{s}}(\Omega^t)})|v_{n+1}|_{2,\infty,\Omega^t}+ |v_{n+1}|_{2,\Omega^t}\cr
&\quad+\|f\|_{W_{\sigma}^{{s},{s}/2}(\Omega^t)}(\varrho^*+t^{\bar{a}} \|\varrho_n\|_{\dot W_{r,\infty}^{1,1}(\Omega^t)})]+\bar d_2,\cr}
\label{5.38}
\end{equation}
where $\bar{a}>0$, $\frac{1}{\sigma}-\frac{1}{r}<1-\frac{s}{ 2}$, $r>{\sigma}$ and
\begin{equation}
\bar d_2=\|d\|_{W_{\sigma}^{2+{s}-\frac{1}{{\sigma}},1+\frac{{s}}{2}-\frac{1}{2{\sigma}}}(S_2^t)}+ \|v(0)\|_{W_{\sigma}^{2+{s}-\frac{2}{{\sigma} }}(\Omega)}.
\label{5.39}
\end{equation}
The above inequality implies (\ref{5.5}) and ends the proof of Lemma~\ref{l5.2}.
\end{proof}

\noindent
 There is yet some problematic term on the r.h.s. of (\ref{5.5}): $\|p_{n+1}\|_{W_{\sigma}^{{s},{s}/2}(\Omega^t)}$. This will be analyzed in the next section.

\section{Estimates for $v_{n+1}$ and $p_{n+1}$}\label{s6}

\begin{remark}\label{r7.1}
We need the imbeddings
$$\eqal{
&|\varrho_1|_{\infty,S_2^t(-a)}\le c\|\varrho_1\|_{W_r^{1,1}(S_2^t(-a))}\quad &{\rm for}\ \ \frac{4}{r}<1,\cr
&\|\varrho_1\|_{L_4(S_2(-a);W_{\sigma}^{{s}/2}(0,t))}\le c\|\varrho_1\|_{W_r^{1,1}(S_2^t(-a))}\quad &{\rm for}\ \ \frac{4}{r}+{s}\le\frac{3}{2}+\frac{2}{\sigma},\cr
&|d_t|_{p,{\sigma},S_2^t}\le c\|d\|_{W_{\sigma}^{2+{s}-\frac{1}{{\sigma}},1+\frac{{s}}{2}-\frac{1}{2{\sigma}}}(S_2^t)}\quad &{\rm for}\ \ \frac{3}{{\sigma}}- \frac{2}{p} \le{s},\cr
&|d_t|_{4,\infty,S_2^t}\le c\|d\|_{W_{\sigma}^{2+{s}-\frac{1}{{\sigma}},1+\frac{{s}}{2}-\frac{1}{2{\sigma}}}(S_2^t)}\quad &{\rm for}\ \ \frac{5}{\sigma}-\frac{1}{2} < {s},\cr
&|d_t|_{L_2(S_2;W_{\sigma}^{{s}/2}(0,t))}\le c\|d\|_{W_{\sigma}^{2+{s}-\frac{1}{{\sigma}},1+\frac{{s}}{2}-\frac{1}{2{\sigma}}}(S_2^t)}\quad &{\rm for}\ \ \frac{3}{\sigma} \le 1,\cr
&|f|_{\frac{p}{r-p},\infty,\Omega^t}\le c\|f\|_{W_{\sigma}^{{s},{s}/2}(\Omega^t)}\quad &{\rm for}\ \ \frac{5}{\sigma}-\frac{2}{ r}+\frac{3}{p} < {s},\cr
&&\phantom{for\ } \frac{1}{ p}=\frac{1}{3}+\frac{1}{\sigma},\cr
&|d_1|_{\infty,S_2^t(-a)}\le c\|d_1\|_{W_{\sigma}^{2+{s}-\frac{1}{{\sigma}},1+\frac{{s}}{2}-\frac{1}{2{\sigma}}}(S_2^t(-a))}\quad &{\rm for}\ \ \frac{5}{\sigma} <  2+{s},\cr
&|\varrho_1|_{3,6,S_2^t(-a)}\le c\|\varrho_1\|_{W_{\sigma}^{1,1}(S_2^t(-a))}\quad &{\rm for}\ \ \frac{4}{{\sigma}} \le 2. \cr}
$$
\end{remark}

In this section we estimate some norms of velocity $v_{n+1}$ and the pressure term $\|p_{n+1}\|_{W_{\sigma}^{{s},{s}/2}(\Omega^t)}$, where $p_{n+1}$ is a solution to the problem
\begin{equation}\eqal{
&\Delta p_{n+1}=\divv(-\varrho_nv_{n+1,t}-\varrho_nv_n\cdot\nabla v_{n+1}+\nu\Delta v_{n+1}+\varrho_nf) \quad {\rm in} \ \Omega,\cr
&\frac{\partial p_{n+1}}{\partial n}\bigg|_{S_2(-a)}\!=\!\bar n\cdot(-\varrho_nv_{n+1,t}-\varrho_nv_n\cdot\nabla v_{n+1}+\nu\Delta v_{n+1}+\varrho_nf)|_{S_2(-a)}\cr
&=-\varrho_1d_{1,t}-\varrho_1\bar n\cdot v_n\cdot\nabla v_{n+1}+\nu\partial_{\tau_\alpha}^2d_1+\nu\bar n\cdot\partial_n^2v_{n+1}+\varrho_1f\cdot\bar n,\cr
&\frac{\partial p_{n+1}}{\partial n}\bigg|_{S_2(a)}=\bar n(-\varrho_nv_{n+1,t}-\varrho_nv_n\cdot\nabla v_{n+1}+\nu\Delta v_{n+1}+\varrho_nf)|_{S_2(a)}\cr
&=-\varrho_nd_{2,t}-\varrho_n\bar n\cdot v_n\cdot\nabla v_{n+1}+\nu\partial_{\tau_\alpha}^2d_2+\nu\bar
n\partial_n^2v_{n+1}+\varrho_nf\cdot\bar n),\cr
&\frac{\partial p_{n+1}}{\partial n}\bigg|_{S_1}=\bar n\cdot(-\varrho_nv_{n+1,t}-\varrho_nv_n\cdot\nabla v_{n+1}+\nu\Delta v_{n+1}+\varrho_nf)|_{S_1}\cr
&=\varrho_nv_nv_{n+1}\cdot\nabla\bar n+\nu\bar n\partial_n^2v_{n+1}+\varrho_nf\cdot\bar n,\cr
& \int_{\Omega}p_{n+1} dx = 0, \cr}
\label{6.1}
\end{equation}
where we introduced a local curvilinear coordinate system on the boundary such that $\Delta=\partial_{\tau_1}^2+\partial_{\tau_2}^2+\partial_n^2$, $\tau_1$, $\tau_2$ are tangent parameters and $n$ is the normal.\\
We can write (\ref{6.1}) in the short form
\begin{equation}\eqal{
&\Delta p_{n+1}=\divv h_n\quad &{\rm in}\ \ \Omega,\cr
&\frac{\partial p_{n+1}}{\partial n}=h_n\cdot\bar n\quad &{\rm on}\ \ S,\cr
&\int_{\Omega}p_{n+1} dx= 0, \cr}
\label{6.2}
\end{equation}
where 
$$
h_n=-\varrho_nv_{n+1,t}-\varrho_nv_n\cdot\nabla v_{n+1}+\nu\Delta v_{n+1}+\varrho_nf.
$$

\begin{lemma}\label{l6.1}
Assume that $v_{n+1},v_n\in V_{\sigma}^{2+{s}}(\Omega^t)$, $\varrho_n\in W_{r,\infty}^{1,1}(\Omega^t)$,\\ $\varrho_1\in W_r^{1,1}(S_2^t(-a))$, $d_i\in W_{\sigma}^{2+{s}-1/{\sigma},1+{s}/2-1/2{\sigma}}(S_2^t(a_i))$, $i=1,2$, $f\in L_{\sigma}(\Omega^t)\cap L_{\frac{pr}{5-p},\infty}(\Omega^t)$, $\partial_t^{{s}/2}f\in L_{\sigma}(0,t;L_p(\Omega))$, $\frac{1}{ p}=\frac{1}{3}+\frac{1}{\sigma}$, $r>{\sigma}>\frac{3}{{s}}$.
Then for solutions to (\ref{6.2}) there exists $\underline{\theta}\in(0,1)$ such that holds the inequality
\begin{equation}\eqal{
&\|p_{n+1}\|_{W_{\sigma}^{{s},{s}/2}(\Omega^t)}\le\varepsilon \|v_{n+1}\|_{V_{\sigma}^{2+{s}}(\Omega^t)}\cr
&\quad+(c(1/\varepsilon)t^{\bar{a}}\|\varrho_n\|_{W_{r,\infty}^{1,1}(\Omega^t)}+ \|\varrho_1\|_{W_r^{1,1}(S_2^t(-a))}+\varrho^*)\cr
&\quad \cdot\sum_{i=1}^2\|d_i\|_{W_{\sigma}^{2+{s}-1/{\sigma},1+{s}/2-1/2{\sigma}}(S_2^t(a_i))}\cr
&\quad+c(1/\varepsilon,\varrho^*)[\|\varrho_n\|_{W_{r,\infty}^{1,1}(\Omega^t)}^\frac{1}{ 1-\underline{\theta}}+\|v_n\|_{V_{\sigma}^{2+{s}}(\Omega^t)}^\frac{1}{1-\underline{\theta}}\cr
&\quad+\|\varrho_n\|_{W_{r,\infty}^{1,1}(\Omega^t)}^\frac{1}{ 1-\underline{\theta}} \|v_n\|_{V_{\sigma}^{2+{s}}(\Omega^t)}^\frac{1}{ 1-\underline{\theta}}+1]t^{\bar{a}} |v_{n+1}|_{2,\infty,\Omega^t}\cr
&\quad+c\varrho^*|f|_{{\sigma},\Omega^t}+ct^{1/{\sigma}}\|\varrho_n\|_{W_{r,\infty}^{1,1}(\Omega^t)} |f|_{\frac{pr}{ r-p},\infty,\Omega^t}\cr
&\quad+c\varrho^*|\partial_t^{{s}/2}f|_{p,{\sigma},\Omega^t},\cr}
\label{6.3}
\end{equation}
where $\frac{1}{p}=\frac{1}{3}+\frac{1}{\sigma}$.
\end{lemma}

\begin{proof}
Let $G=G(x,y)$ be the Green function to the Neumann problem (\ref{6.2}). From the properties of the Green function, we have
\begin{equation}\eqal{
p_{n+1}&= \int_{\Omega} G \divv h_n dx + \int_S \bar{n}\cdot \nabla G p_{n+1} dS - \int_S G\bar{n} \cdot \nabla p_{n+1} dS \cr &  =  \int_{\Omega} G \divv h_n dx - \int_S G\bar{n} \cdot \nabla p_{n+1} dS = \int_{\Omega} \divv(G h_n) dx \cr & - \int_{\Omega} \nabla G \cdot h_n dx  - \int_S G \bar{n} \cdot \nabla p_{n+1} dS   =  \intop_\Omega h_n\cdot\nabla Gdx \cr & =\intop_\Omega(\varrho_nv_{n+1,t}+\varrho_nv_n\cdot\nabla v_{n+1}-\nu\Delta v_{n+1}-\varrho_nf)\nabla Gdx \cr & = \intop_\Omega\varrho_nv_{n+1,t}\cdot\nabla Gdx+ \intop_\Omega(\varrho_nv_n\cdot\nabla v_{n+1} \cr & -\nu\Delta v_{n+1}-\varrho_nf)\cdot\nabla Gdx   \equiv J_1+J_2.\cr}
\label{6.4}
\end{equation}
Consider $J_1$. Integrating by parts yields
\begin{equation}
J_1=-\intop_\Omega\nabla\varrho_n\cdot v_{n+1,t}Gdx+\intop_S\varrho_nv_{n+1,t}\cdot\bar nGdS\equiv J_{11}+J_{12},
\label{6.5}
\end{equation}
where
$$
J_{12}=\intop_{S_2(-a)}\varrho_nd_{1,t}GdS_2+\intop_{S_2(a)}\varrho_nd_{2,t}GdS_2 \equiv J_{12}^1+J_{12}^2,
$$
and we used that $v_{n+1}\cdot\bar n|_{S_1}=0$.

Continuing, we have
$$
\|J_{12}^1\|_{W_{\sigma}^{{s},{s}/2}(\Omega^t)}= \|J_{12}^1\|_{L_{\sigma}(0,t;W_{\sigma}^{s}(\Omega))}+ \|J_{12}^1\|_{L_{\sigma}(\Omega;W_{\sigma} ^{{s}/2}(0,t))}\equiv L_1^1+L_2^1.
$$
To examine $L_1^1$ we use the estimate
$$ L_1^1 \le c\|J^1_{12}\|_{L_\sigma(0,t;W^1_\sigma(\Omega))}. $$
Then, we use the proof of Lemma \ref{l2.7} to estimate the singular part of $J_{12}^1$ as follows
$$
\bigg(\intop_{\R^2}dx'\bigg|\intop_{\R^2}\bigg(\intop_{\R_+}dx_3 \bigg|\frac{x_3}{(\sqrt{|x'-y'|^2+x_3^2})^3}\bigg|^{\sigma}\bigg)^{1/{\sigma} }\varrho_1 d_{1,t}dy'\bigg|^{\sigma}\bigg)^{1/{\sigma}}.
$$
Next, the integration with respect to $x_3$, i.e. applying the Minkowski inequality (as in the proof of Lemma \ref{l2.7}) gives
$$
\bigg(\intop_{\R^2}dx'\bigg|\intop_{\R^2}\bigg|\frac{1}{\sqrt{|x'-y'|^2}}\bigg|^{2-1/{\sigma}} \varrho_1d_{1,t}dy'\bigg|^{\sigma}\bigg)^{1/{\sigma} }\equiv I'_3.
$$
Since
$$
K(x'-y')=\bigg(\frac{1}{\sqrt{|x'-y'|^2}}\bigg)^{2-1/{\sigma}}\in L_{s'}(\R^2)\quad {\rm for}\ \ {s'}<\frac{2}{ 2-1/{\sigma}}
$$
and $x'$, $y'$ belong to some compact set, we obtain by the Young inequality that
$$
I'_3\le c\|\varrho_1d_{1,t}\|_{L_p(\R^2)},
$$
where $1+1/{\sigma}-1/p=1/{s'}>\frac{2-1/{\sigma}}{2}$ so $p>\frac{2}{3}{\sigma}$.

Then
$$
L_1^1\le c|\varrho_1|_{\infty,S_2^t(-a)}\|d_{1,t}\|_{L_p(S_2^t(-a))}.
$$
Similarly as above, using Lemma \ref{l2.7},
$$\eqal{
L_2^1&\le c\|\,\|\varrho_1d_{1,t}\|_{W_{\sigma}^{{s}/2}(0,t)}\|_{L_2(S_2(-a))}\cr
&\le c\|\varrho_1\|_{L_4(S_2(-a);W_{\sigma}^{{s}/2}(0,t))} \|d_{1,t}\|_{L_\infty(0,t;L_4(S_2(-a)))}\cr
&\quad+c|\varrho_1|_{\infty,S_2^t(-a)}\|d_{1,t}\|_{L_2(S_2(-a);W_{\sigma}^{{s}/2}(0,t))}\cr
&\le c\|\varrho_1\|_{W_r^{1,1}(S_2^t(-a))}\cdot \|d_1\|_{W_{\sigma}^{2+{s}-1/{\sigma},1+\frac{s}{ 2}-1/2{\sigma}}(S_2^t(-a))},\cr}
$$
where we used Remark \ref{r7.1}.

For $J_{12}^2$ we obtain the similar estimate, where $\varrho_1d_{1,t}$ is replaced by $\varrho_nd_{2,t}$, respectively. Hence
$$\eqal{
&\|\,\|\varrho_nd_{2,t}\|_{W_{\sigma}^{{s}/2}(0,t)}\|_{L_2(S_2)}\le\varepsilon t^{\bar{a}}\|\varrho_n\|_{L_4(S_2;W_{\sigma}^{{s}/2}(0,t))} \|d_{2,t}\|_{L_\infty(0,t;L_4(S_2(a)))}\cr
&\quad+c\varrho^*\|d_{2,t}\|_{L_2(S_2(a);W_{\sigma}^{{s}/2}(0,t))}\cr
&\le(\varepsilon t^{\bar{a}}\|\varrho_n\|_{W_{r,\infty}^{1,1}(\Omega^t)}+c\varrho^*)\cdot \|d_2\|_{W_{\sigma}^{2+{s}-1/{\sigma},1+{s}/2-1/2{\sigma} }(S_2^t(a))},\cr}
$$
where we used Remark \ref{r7.1} and the imbedding
$$
\|\varrho_n\|_{L_4(S_2;W_{\sigma}^{{s}/2}(0,t))}\le c\|\varrho_n\|_{W_{r,\infty}^{1,1}(\Omega^t)}.
$$
Exploiting the above estimates for $J_{12}$ we obtain
\begin{equation}\eqal{
\|J_{12}\|_{W_{\sigma}^{{s},{s}/2}(\Omega^t)}&\le(\varepsilon t^{\bar{a}}\|\varrho_n\|_{W_{r,\infty}^{1,1}(\Omega^t)}+c\varrho^*+ \|\varrho_1\|_{W_r^{1,1}(S_2^t(-a))})\cdot\cr
&\quad\cdot\sum_{i=1}^2\|d_i\|_{W_{\sigma}^{2+{s}-1/{\sigma},1+{s}/2-1/2{\sigma}}(S_2^t(a_i))}.\cr}
\label{6.6}
\end{equation}
Next, we calculate
\begin{equation}
\|J_{11}\|_{W_{\sigma}^{{s},{s}/2}(\Omega^t)}=\|J_{11}\|_{L_{\sigma}(0,t;W_{\sigma}^{s}(\Omega))} +\|J_{11}\|_{L_{\sigma}(\Omega;W_{\sigma} ^{{s}/2}(0,t))}\equiv K_1+K_2.
\label{6.7}
\end{equation}
By the Young inequality (see \cite[Ch. 1, Sect. 2.14]{BIN}), we have
$$\eqal{
K_1&\le\bigg(\intop_0^t\bigg\|\intop_\Omega G\nabla\varrho_nv_{n+1,t}\bigg\|^{\sigma}_{W_{\sigma}^{s}(\Omega)}dt\bigg)^{1/{\sigma}} \cr & \le c \bigg(\intop_0^t\bigg\|\intop_\Omega G\nabla\varrho_nv_{n+1,t}\bigg\|^{\sigma}_{W_{\sigma}^{1}(\Omega)}dt\bigg)^{1/{\sigma}}\cr
&\le c\bigg(\intop_0^t|\nabla\varrho_nv_{n+1,t}|_{p,\Omega}^{\sigma} dt\bigg)^{1/{\sigma}}\le c\bigg(\intop_0^t|\nabla\varrho_n|_{p'',\Omega}^{\sigma} |v_{n+1,t}|_{{\sigma},\Omega}^{\sigma} dt\bigg)^{1/{\sigma}}\cr
&\le c|\nabla\varrho_n|_{p'',\infty,\Omega^t}|v_{n+1,t}|_{{\sigma},\Omega^t}\equiv K_1^1,\cr}
$$
where $1-\frac{1}{p}+\frac{1}{\sigma}>\frac{2}{3}$, so $\frac{1}{3}+\frac{1}{\sigma}>\frac{1}{p}=\frac{1}{\sigma}+\frac{1}{p''}$. Thus $p''>3$. We use the interpolation
$$
|v_{n+1,t}|_{{\sigma},\Omega^t}\le c\|v_{n+1}\|_{W_{\sigma}^{2+{s},1+{s}/2}(\Omega^t)}^{\theta_1} |v_{n+1}|_{2,\Omega^t}^{1-\theta_1},
$$
where $\theta_1$ satisfies
$$
\frac{5}{\sigma}-2=(1-\theta_1)\frac{5}{2}+\theta_1\bigg(\frac{5}{\sigma}-2-{s}\bigg),
$$
so
$$
\theta_1=\frac{2+5/2-5/{\sigma}}{2+{s}+5/2-5/{\sigma}}.
$$
Then
\begin{equation}
K_1\le K_1^1\le\varepsilon\|v_{n+1}\|_{W_{\sigma}^{2+{s},1+{s}/2}(\Omega^t)}+ (c(1/\varepsilon)|\nabla\varrho_n|_{p'',\infty,\Omega^t}^{1/(1-\theta_1)}) |v_{n+1}|_{2,\Omega^t},
\label{6.8}
\end{equation}
where $p''\le r$ so $|\nabla\varrho_n|_{p'',\infty,\Omega^t}\le c\|\varrho_n\|_{W_{r,\infty}^{1,1}(\Omega^t)}$.

Finally, we calculate
$$\eqal{
K_2&\le\bigg(\intop_\Omega\intop_0^t\bigg|\intop_\Omega G\partial_t^{{s}/2}\nabla\varrho_nv_{n+1,t}dx\bigg|^{\sigma} dxdt\bigg)^{1/{\sigma}}\cr
&\quad+\bigg(\intop_\Omega\intop_0^t\bigg|\intop_\Omega G\nabla\varrho_n\partial_t^{{s}/2}v_{n+1,t}dx\bigg|^{\sigma} dxdt\bigg)^{1/{\sigma}}\equiv I_1+I_2.\cr}
$$
To examine $I_2$ we use the Young inequality from \cite[Ch. 1, Sect. 2.14]{BIN}. Then we get
$$\eqal{
I_2&\le c\bigg(\intop_0^t|\nabla\varrho_n\partial_t^{{s}/2}v_{n+1,t}|_{p,\Omega}^{\sigma} dt \bigg)^{1/{\sigma}}\cr
&\le c\bigg(\intop_0^t|\nabla\varrho_n|_{p',\Omega}^{\sigma} |\partial_t^{{s}/2}v_{n+1,t}|_{2,\Omega}^{\sigma} dt\bigg)^{1/{\sigma}}\cr
&\le c|\nabla\varrho_n|_{p',\infty,\Omega^t}\bigg(\intop_0^t |\partial_t^{{s}/2}v_{n+1,t}|_{2,\Omega}^{\sigma} dt\bigg)^{1/{\sigma}}\equiv I_2^1,\cr}
$$
where
$$
1-\frac{1}{p}+\frac{1}{\sigma}>\frac{1}{3},\quad {\rm so}\quad \frac{2}{3}+\frac{1}{\sigma}>\frac{1}{p}=\frac{1}{2}+\frac{1}{ p'}.
$$
Hence $\frac{1}{6}+\frac{1}{\sigma}>\frac{1}{p'}$ and $p'>\frac{6{\sigma}}{\sigma+6}$.

Now, we use the interpolation
$$
|\partial_t^{{s}/2}v_{n+1,t}|_{2,{\sigma},\Omega^t}\le c\|v_{n+1}\|_{W_{\sigma}^{2+{s},1+{s}/2}(\Omega^t)}^{\theta_2} |v_{n+1}|_{2,\Omega^t}^{1-\theta_2},
$$
where $\theta_2$ is a solution to the equation
$$
\frac{3}{2}+\frac{2}{\sigma}-(2+{s})=(1-\theta_2)\frac{5}{2}+\theta_2\bigg(\frac{5}{\sigma}-(2+{s})\bigg),
$$
so
$$
\theta_2=\frac{2+{s}+\frac{5}{2}-\frac{3}{2}-\frac{2}{\sigma}}{2+{s}+\frac{5}{2}-\frac{5}{\sigma}}
$$
and $\theta_2<1$ because $3/{\sigma}<3/2$.

Therefore,
\begin{equation}
I_2^1\le\varepsilon\|v_{n+1}\|_{W_{\sigma}^{2+{s},1+{s}/2}(\Omega^t)}+ (c(1/\varepsilon)|\nabla\varrho_n|_{p',\infty,\Omega^t})^{1/(1-\theta_2)} |v_{n+1}|_{2,\Omega^t}.
\label{6.9}
\end{equation}
Next, we analyze
$$\eqal{
I_1&\le\bigg(\intop_{\Omega^t}\bigg|\intop_\Omega G\partial_t^{{s}/2}\nabla\varrho_nv_{n+1,t}dx\bigg|^{\sigma} dxdt\bigg)^{1/{\sigma}}\cr
&\le\bigg(\intop_{\Omega^t}\bigg|\intop_\Omega\nabla G\partial_t^{{s}/2}\varrho_nv_{n+1,t}dx\bigg|^{\sigma} dxdt\bigg)^{1/{\sigma}}\cr
&\quad+\bigg(\intop_{\Omega^t}\bigg|\intop_SG\partial_t^{{s}/2}\varrho_n v_{n+1,t}\cdot\bar ndS\bigg|^{\sigma} dxdt\bigg)^{1/{\sigma}}\equiv I_1^1+I_1^2,\cr}
$$
where
$$
I_1^1\le c|\partial_t^{{s}/2}\varrho_n|_{p',\infty,\Omega^t} |v_{n+1,t}|_{{\sigma},\Omega^t}\equiv I_1^{11}
$$
and
$$
1-\frac{1}{p}+\frac{1}{\sigma}>\frac{2}{3},\quad {\rm so}\quad \frac{1}{3}+\frac{1}{\sigma}>\frac{1}{p}=\frac{1}{\sigma}+\frac{1}{p'},\quad {\rm thus}\ \  p'>3.
$$
In view of the interpolation
$$
|v_{n+1,t}|_{{\sigma},\Omega^t}\le c\|v_{n+1}\|_{W_{\sigma}^{2+{s},1+{s}/2}(\Omega^t)}^{\theta_3} |v_{n+1}|_{2,\Omega^t}^{1-\theta_3},
$$
where $\theta_3$ is a solution to the equation
$$
\frac{5}{\sigma}-2=(1-\theta_3)\frac{5}{2}+\theta_3\bigg(\frac{5}{\sigma}-(2+{s})\bigg),
$$
thus
$$
\theta_3=\frac{2+5/2-5/{\sigma}}{2+{s}+5/2-5/{\sigma}}.
$$
Therefore, we obtain
\begin{equation}
I_1^{11}\le\varepsilon\|v_{n+1]}|_{W_{\sigma}^{2+{s},1+{s}/2}(\Omega^t)}+ (c(1/\varepsilon)|\partial_t^{{s}/2}\varrho_n|_{p',\infty,\Omega^t})^{1/(1-\theta_3)} |v_{n+1}|_{2,\Omega^t},
\label{6.10}
\end{equation}
where
$$
|\partial_t^{{s}/2}\varrho_n|_{p',\infty,\Omega^t}\le c\|\varrho_n\|_{W_{r,\infty}^{1,1}(\Omega^t)}\quad {\rm for}\ \ r\ge p'.
$$
In view of boundary conditions and Lemma \ref{l2.7}, we have
\begin{equation}\eqal{
I_1^2&\le\bigg(\intop_{\Omega^t}\bigg|\intop_{S_2(-a)}G\partial_t^{{s}/2} \varrho_1d_{1,t}dS_2\bigg|^{\sigma} dxdt\bigg)^{1/{\sigma}}\cr
&\quad+\bigg(\intop_{\Omega^t}\bigg|\intop_{S_2(a)}G\partial_t^{{s}/2} \varrho_nd_{2,t}dS_2\bigg|^{\sigma} dxdt\bigg)^{1/{\sigma}}\cr
&\le\bigg(\intop_0^t|\partial_t^{{s}/2}\varrho_1d_{1,t}|_{p,S_2(-a)}^{\sigma} dt'\bigg)^{1/{\sigma}} +\bigg(\intop_0^t|\partial_t^{{s}/2}\varrho_nd_{2,t}|_{p,S_2(a)}^{\sigma} dt'\bigg)^{1/{\sigma}}\cr
&\le c\sup_t|\partial_t^{{s}/2}\varrho_1|_{p',S_2(-a)}|d_{1,t}|_{{\sigma},S_2^t(-a)}\cr
&\quad+c\sup_t|\partial_t^{{s}/2}\varrho_n|_{p',S_2(a)} |d_{2,t}|_{{\sigma},S_2^t(a)}\equiv I_1^{21},\cr}
\label{6.11}
\end{equation}
where $p> \frac{3{\sigma}}{ {\sigma}+3}$ so $\frac{1}{2}+\frac{3}{ 2{\sigma}}>\frac{1}{p}=\frac{1}{ p'}+\frac{1}{\sigma}$. Hence $p'>\frac{2{\sigma}}{ {\sigma} +1}$.

We imply the imbedding
$$
\sup_t|\partial_t^{{s}/2}\varrho_n|_{p',S_2}\le c\|\varrho_n\|_{W_{{\sigma},\infty}^{1,1}(\Omega^t)},
$$
which holds for
$$
\frac{3}{{\sigma}}-\frac{2}{p'}+\frac{s}{2}\le 1,\quad {\rm so}\quad \frac{3}{\sigma}+\frac{s}{2}\le 1+\frac{2}{p'},\ \ p'>\frac{2\sigma}{\sigma +1}.
$$
The above restrictions always hold for ${\sigma}>3$, ${s}<1$.
\begin{equation}\eqal{
I_1^{12}&\le c(t^{\bar{a}}\|\varrho_n\|_{W_{{\sigma},\infty}^{1,1}(\Omega^t)}\cr
&\quad+ \|\varrho_1\|_{W_r^{1,1}(S_2^t(-a))})\sum_{i=1}^2 \|d_i\|_{W_{\sigma}^{2+{s}-1/{\sigma},1+{s}/2-1/2{\sigma}}(S_2^t(a_i))}.\cr}
\label{6.12}
\end{equation}
Now we find estimate for $\|J_{11}\|_{W_{\sigma}^{{s},{s}/2}(\Omega^t)}$ that appeared in (\ref{6.7}). From the above estimates we have
\begin{equation}
\|J_{11}\|_{W_{\sigma}^{{s},{s}/2}(\Omega^t)}\le K_1+K_2,
\label{6.13}
\end{equation}
where $K_1\le K_1^1$ and $K_1^1$ is bounded in (\ref{6.8}) and $K_2\le I_1+I_2\le I_2^1+I_1+I_1^2$.

Continuing, we have that $I_2^1$ is estimated in (\ref{6.9}), $I_1\le I_1^{11}$ is estimated in (\ref{6.10}).

Finally, $I_1^2$ is bounded in (\ref{6.11}) and $I_1^2\le I_1^{21}$, which is estimated in (\ref{6.12}).

Exploiting estimates (\ref{6.8})--(\ref{6.12}) in (\ref{6.13}) yields
\begin{equation}\eqal{
&\|J_{11}\|_{W_{\sigma}^{{s},{s}/2}(\Omega^t)}
\le\varepsilon\|v_{n+1}\|_{W_{\sigma}^{2+{s},1+{s}/2}(\Omega^t)}\cr
&\quad+(c(1/\varepsilon)\|\varrho_n\|_{W_{r,\infty}^{1,1}(\Omega^t)})^{1/(1-\theta_*)} |v_{n+1}|_{2,\Omega^t}+c(t^{\bar{a}}\|\varrho_n\|_{W_{r,\infty}^{1,1}(\Omega^t)}\cr
&\quad+ \|\varrho_1\|_{W_r^{1,1}(S_2^t(-a))})\sum_{i=1}^2 \|d_i\|_{W_{\sigma}^{2+{s}-1/{\sigma},1+{s}/2-1/2{\sigma}}(S_2^t(-a))},\cr}
\label{6.14}
\end{equation}
where
$$
\frac{1}{1-\theta_*}=\max_{i\in\{1,2,3\}}\frac{1}{1-\theta_i}.
$$
Next, we consider
\begin{equation}
\|J_2\|_{\dot W_{\sigma}^{{s},{s}/2}(\Omega^t)}= \|J_2\|_{L_{\sigma}(0,t;W_{\sigma}^{s}(\Omega))}+ \|J_2\|_{L_{\sigma}(\Omega;W_{\sigma} ^{{s}/2}(0,t))}\equiv J_2^1+J_2^2,
\label{6.15}
\end{equation}
where
\begin{equation}\eqal{
J_2^1&\le\|J_2\|_{L_{\sigma}(0,t;W_{\sigma}^1(\Omega))}\le\bigg(\intop_0^tdt \bigg\|\intop_\Omega(\varrho_nv_n\cdot\nabla v_{n+1}\cr
&\quad-\nu\Delta v_{n+1}-\varrho_nf)\nabla Gdx\bigg\|_{W_{\sigma}^1(\Omega)}^{\sigma} \bigg)^{1/{\sigma}}\cr
&\le c\varrho^*\bigg(\intop_{\Omega^t}|v_n\cdot\nabla v_{n+1}|^{\sigma} dxdt\bigg)^{1/{\sigma}}+c\bigg(\intop_{\Omega^t}|\Delta v_{n+1}|^{\sigma} dxdt\bigg)^{1/{\sigma}}\cr
&\quad+c\bigg(\intop_{\Omega^t}|\varrho_nf|^{\sigma} dxdt\bigg)^{1/{\sigma}}\equiv L_1+L_2+L_3.\cr}
\label{6.16}
\end{equation}
First, we estimate
\begin{equation}
L_1\le c\varrho^*\bigg(\intop_0^t|v_n|_{{\sigma}\lambda_1,\Omega}^{\sigma}|\nabla v_{n+1}|_{{\sigma}\lambda_2,\Omega}^{\sigma} dt'\bigg)^{1/{\sigma} }\equiv L_1^1,
\label{6.17}
\end{equation}
where $1/\lambda_1+1/\lambda_2=1$.

Using the imbedding
$$
|v_n|_{{\sigma}\lambda_1,\Omega}\le c\|v_n\|_{W_{\sigma}^{2+{s}-2/{\sigma}}(\Omega)},\quad \frac{3}{\sigma}-\frac{3}{\sigma\lambda_1}\le 2+{s}-\frac{2}{{\sigma}}
$$
and the interpolation
$$
|\nabla v_{n+1}|_{{\sigma}\lambda_2,\Omega}\le c\|v_{n+1}\|_{W_{\sigma}^{2+{s}-2/{\sigma}}(\Omega)}^{\theta_4}|v_{n+1}|_{2,\Omega}^{1-\theta_4},
$$
where $\theta_4$ satisfies
$$
\frac{3}{{\sigma}\lambda_2}-1=(1-\theta_4)\frac{3}{2}+\theta_4\bigg(\frac{3}{\sigma}-(2+{s}-2/{\sigma})\bigg),
$$
so
$$
\theta_4=\frac{1+3/2-3/{\sigma}\lambda_2}{2+{s}-2/{\sigma}+3/2-3/{\sigma}},
$$
where the condition $\theta_4<1$ and the above imbedding imply
$$
\frac{7}{\sigma}\le 3+2{s},
$$
which is always satisfied.

Then,
\begin{equation}\eqal{
L_1^1&\le c\varrho^*\|v_n\|_{V_{\sigma}^{2+{s}}(\Omega^t)}\bigg(\intop_0^t \|v_{n+1}\|_{W_{\sigma}^{2+{s}-2/{\sigma}}(\Omega)}^{\theta_4{\sigma}} |v_{n+1}|_{2,\Omega}^{(1-\theta_4){\sigma}}dt\bigg)^{1/{\sigma}}\cr
&\le c\varrho^*\|v_n\|_{V_{\sigma}^{2+{s}}(\Omega^t)}\|v_{n+1} \|_{V_{\sigma}^{2+{s}}(\Omega^t)}^{\theta_4} t^{1/{\sigma} }|v_{n+1}|_{2,\infty,\Omega^t}^{1-\theta_4}\cr
& \le\varepsilon\|v_{n+1}\|_{V_{\sigma}^{2+{s}}(\Omega^t)}\! \cr
& \quad +\! \bigg(\!c(1/\varepsilon)\varrho^* \|v_n\|_{V_{\sigma}^{2+{s}}(\Omega^t)}t^{1/{\sigma} }\bigg)^{1/(1-\theta_4)}\!t^{\frac{1}{2}} |v_{n+1}|_{2,\infty,\Omega^t}.\cr}
\label{6.18}
\end{equation}
Next,
\begin{equation}
L_2=|\Delta v_{n+1}|_{{\sigma},\Omega^t}\le c\|v_{n+1}\|_{W_{\sigma}^{2+{s},1+{s}/2}(\Omega^t)}^{\theta _5}|v_{n+1}|_{2,\Omega^t}^{1-\theta_5}\equiv L_2^1,
\label{6.19}
\end{equation}
where $\theta_5$ satisfies
$$
\frac{5}{\sigma}-2=(1-\theta_5)\frac{5}{2}+\theta_5\bigg(\frac{5}{\sigma}-(2+{s})\bigg)\quad {\rm so}\quad \theta_5=\frac{2+5/2-5/{\sigma}}{2+{s}+5/2-5/{\sigma}}.
$$
Then
\begin{equation}
L_2^1\le\varepsilon\|v_{n+1}\|_{W_{\sigma}^{2+{s},1+{s}/2}(\Omega^t)}+ c(1/\varepsilon)|v_{n+1}|_{2,\Omega^t}.
\label{6.20}
\end{equation}
Finally,
\begin{equation}
L_3\le\varrho^*|f|_{{\sigma},\Omega^t}.
\label{6.21}
\end{equation}
Using (\ref{6.17})--(\ref{6.21}) in (\ref{6.16}) yields
\begin{equation}\eqal{
J_2^1&\le\varepsilon\|v_{n+1}\|_{V_{\sigma}^{2+{s}}(\Omega^t)}+(c(1/\varepsilon) \varrho^*\|v_n\|_{V_{\sigma} ^{2+{s}}(\Omega^t)})^{1/(1-\theta_4)}t^{1/{\sigma}} |v_{n+1}|_{2,\infty,\Omega^t}\cr
&\quad+c(1/\varepsilon)|v_{n+1}|_{2,\Omega^t}+\varrho^*|f|_{{\sigma},\Omega^t}.\cr}
\label{6.22}
\end{equation}
At the end, we consider $J_2^2$,
$$\eqal{
J_2^2&=\bigg(\intop_\Omega dx\bigg|\intop_\Omega\|\varrho_nv_n\cdot\nabla v_{n+1}-\nu\Delta v_{n+1}-\varrho_nf\|_{W_{\sigma}^{{s}/2}(0,t)}\nabla Gdx\bigg|^{\sigma}\bigg)^{1/{\sigma}}\cr
&=\bigg(\intop_\Omega dx\bigg|\intop_\Omega\intop_0^t(\partial_t^{{s}/2}(\varrho_nv_n\cdot\nabla v_{n+1})-\partial_t^{{s}/2}\nu\Delta v_{n+1}\cr
&\quad-\partial_t^{{s}/2}(\varrho_nf))dt'\nabla Gdx'\bigg|^{\sigma}\bigg)^{1/{\sigma}}\cr
&\le\bigg(\intop_\Omega dx\bigg|\intop_\Omega\intop_0^t\partial_t^{{s}/2} (\varrho_nv_n\cdot\nabla v_{n+1})\nabla Gdx'dt'\bigg|^{\sigma} \bigg)^{1/{\sigma}}\cr
&\quad+\bigg(\intop_\Omega dx\bigg|\intop_\Omega\intop_0^t\partial_t^{{s}/2} \Delta v_{n+1}\nabla Gdx'dt'\bigg|^{\sigma}\bigg)^{1/{\sigma}}\cr
&\quad+\bigg(\intop_\Omega dx\bigg|\intop_\Omega\intop_0^t\partial_t^{{s}/2} (\varrho_nf)\nabla Gdx'dt'\bigg|^{\sigma}\bigg)^{1/{\sigma}}\equiv K_1+K_2+K_3.\cr}
$$
Consider $K_2$. By the Minkowski and Young inequalities we have
\begin{equation}
K_2\le c\bigg(\intop_0^t|\partial_t^{{s}/2}\Delta v_{n+1}|_{p,\Omega}^{\sigma} dt\bigg)^{1/{\sigma}}\equiv K_2^1,
\label{6.23}
\end{equation}
where the Young theorem (see \cite[Ch. 1, Sect. 2.14]{BIN}) is used with
$$
1-\frac{1}{p}+\frac{1}{\sigma}=\frac{1}{\delta}>\frac{2}{3}\quad {\rm so}\quad \frac{1}{3}+\frac{1}{\sigma}>\frac{1}{p}.
$$
Hence,
$$
K_2^1=c|\partial_t^{{s}/2}\Delta v_{n+1}|_{p,{\sigma},\Omega^t}\le c\|v_{n+1}\|_{W_{\sigma}^{2+{s},1+{s}/2}(\Omega^t)}^{\theta_5} |v_{n+1}|_{2,\Omega^t}^{1-\theta_5},
$$
where
$$
\frac{3}{p}+\frac{2}{\sigma}-(2+{s})=(1-\theta_5)\frac{5}{2}+\theta_5\left(\frac{5}{\sigma}-(2+{s})\right),
$$
so
$$
\theta_5=\frac{2+{s}+5/2-(3/p+2/{\sigma})}{2+{s}+5/2-5/{\sigma}}.
$$
The condition $\theta_5<1$ implies that $\frac{3}{p}>\frac{3}{\sigma}$.

Hence
\begin{equation}
K_2\le K_2^1\le\varepsilon\|v_{n+1}\|_{W_{\sigma}^{2+{s},1+{s}/2}(\Omega^t)}+ c(1/\varepsilon)|v_{n+1}|_{2,\Omega^t}.
\label{6.24}
\end{equation}
Comparing this with restriction $\frac{1}{3}+\frac{1}{\sigma}>\frac{1}{p}$ implies that $\frac{1}{3}+\frac{1}{\sigma}>\frac{1}{\sigma},$ so there is no restriction. Now, we estimate $K_1$.

Consider $K_1$. Then we have
\begin{equation}\eqal{
K_1&\le\bigg(\intop_\Omega dx\bigg|\intop_\Omega\intop_0^t\partial_t^{{s}/2}\varrho_nv_n\nabla v_{n+1}\nabla Gdx'dt'\bigg|^{\sigma}\bigg)^{1/{\sigma} }\cr
&\quad+\bigg(\intop_\Omega dx\bigg|\intop_\Omega\intop_0^t\varrho_n\partial_t^{{s}/2}v_n\nabla v_{n+1}\nabla Gdx'dt'\bigg|^{\sigma}\bigg)^{1/{\sigma} }\cr
&\quad+\bigg(\intop_\Omega dx\bigg|\intop_\Omega\intop_0^t\varrho_nv_n\partial_t^{{s}/2}\nabla v_{n+1}\nabla Gdx'dt'\bigg|^{\sigma}\bigg)^{1/{\sigma} }\cr
&\equiv K_1^1+K_1^2+K_1^3.\cr}
\label{6.25}
\end{equation}
Applying the Minkowski and Young inequalities
$$\eqal{
K_1^1&\le\sup_t|v_n|_{\infty,\Omega}\bigg(\intop_0^t|\partial_t^{{s}/2}\varrho_n\nabla v_{n+1}|_{p,\Omega}^{\sigma} dt'\bigg)^{1/{\sigma}}\cr
&\le\sup_t\|v_n\|_{W_{\sigma}^{2+{s}-2/{\sigma}}(\Omega)}\bigg(\intop_0^t |\partial_t^{{s}/2}\varrho_n|_{p',\Omega}^{\sigma}|\nabla v_{n+1}|_{p'',\Omega}^{\sigma} dt'\bigg)^{1/{\sigma}}\cr
&\equiv K_1^{11},\cr}
$$
where $\frac{1}{3}+\frac{1}{\sigma}>\frac{1}{p}$, $\frac{1}{p'}+\frac{1}{p''}=\frac{1}{p}$.

Setting $p'=r$, we obtain
$$\eqal{
K_1^{11}&\le\sup_t\|v_n\|_{W_{\sigma}^{2+{s}-2/{\sigma}}(\Omega)}\sup_t |\partial_t\varrho_n|_{r,\Omega}\bigg(\intop_0^t|\nabla v_{n+1}|_{p'',\Omega}^{\sigma} dt'\bigg)^{1/{\sigma}}\cr
&\equiv K_1^{12},\cr}
$$
where $\frac{1}{p''}=\frac{1}{p}-\frac{1}{r}$.

We apply the interpolation
$$
|\nabla v_{n+1}|_{p'',\Omega}\le c\|v_{n+1}\|_{W_{\sigma}^{2+{s}}(\Omega)}^{\theta_6}|v_{n+1}|_{2,\Omega}^{1-\theta_6},
$$
where $\theta_6$ is a solution to the equation
$$
\frac{3}{p''}-1=(1-\theta_6)\frac{3}{2}+\theta_6\bigg(\frac{3}{\sigma}-(2+{s})\bigg),
$$
so
$$
\theta_6=\frac{\frac{3}{2}+1-\frac{3}{p''}}{2+{s}+\frac{3}{2}-\frac{3}{\sigma}}.
$$
Hence the condition $\theta_6<1$ yields $\frac{3}{{\sigma}}<1+\frac{3}{p''}=1+\frac{3}{p}-\frac{3}{r}$.

Using that $\frac{1}{p}<\frac{1}{3}+\frac{1}{{\sigma}}$ we finally derive the restriction
$$
\frac{3}{r}<2.
$$
Hence, there is no restriction because $r>5$.

By the Young inequality
\begin{equation}\eqal{
K_1^1&\le K_1^{12}\le\varepsilon\|v_{n+1}\|_{W_{\sigma}^{2+{s},1+{s}/2}(\Omega^t)}\cr
&\quad+ (c(1/\varepsilon)|\partial_t\varrho_n|_{r,\infty,\Omega^t} \|v_n\|_{V_{\sigma}^{2+{s}}(\Omega^t)})^{1/1-\theta_6}\cdot |v_{n+1}|_{2,{\sigma} ,\Omega^t}.\cr}
\label{6.26}
\end{equation}
Next, we examine
$$\eqal{
K_1^2&\le\varrho^*\bigg(\intop_\Omega dx\bigg|\intop_\Omega\intop_0^t \partial_t^{{s}/2}v_n\nabla v_{n+1}\nabla Gdx'dt'\bigg|^{\sigma} \bigg)^{1/{\sigma}}\cr
&\le\varrho^*\bigg(\intop_0^t|\partial_t^{{s}/2}v_n\nabla v_{n+1}|_{p,\Omega}^{\sigma} dt'\bigg)^{1/{\sigma}}\equiv K_1^{21},\cr}
$$
where $\frac{1}{3}+\frac{1}{\sigma}>\frac{1}{p}$ and the Minkowski and Young inequalities were used.

By the H\"older inequality we have
$$
K_1^{21}\le\varrho^*\bigg(\intop_0^t|\partial_t^{{s}/2}v_n|_{p',\Omega}^{\sigma} |\nabla v_{n+1}|_{p'',\Omega}^{\sigma} dt'\bigg)^{1/{\sigma}}\equiv K_1^{22},
$$
where $1/p'+1/p''=1/p$.

In view of the imbedding
$$\eqal{
&\sup_t|\partial_t^{{s}/2}v_n|_{p',\Omega}\le\sup_t \|\partial_t^{{s}/2}v_n\|_{W_{\sigma}^{2-2/{\sigma}}(\Omega)}\le c\|\partial_t^{{s}/2}v_n\|_{W_{\sigma}^{2,1}(\Omega^t)}\cr
&\le c\|v_n\|_{W_{\sigma}^{2+{s},1+{s}/2}(\Omega^t)},\cr}
$$
which holds for $\frac{3}{\sigma}-\frac{3}{p'}<2-\frac{2}{\sigma},$ so $\frac{5}{\sigma}-\frac{3}{p'}<2$ and the interpolation
$$
|\nabla v_{n+1}|_{p'',\Omega}\le c\|v_{n+1}\|_{W_{\sigma}^{2+{s}}(\Omega)}^{\theta_7}|v_{n+1}|_{2,\Omega}^{1-\theta_7},
$$
with $\theta_7$ satisfying the equation
$$
\frac{3}{p''}-1=(1-\theta_7)\frac{3}{2}+\theta_7\bigg(\frac{3}{\sigma}-(2+{s})\bigg), \quad {\rm so}\quad \theta_7=\frac{1+\frac{3}{2}-\frac{3}{p''} }{ 2+{s}+\frac{3}{2}-\frac{3}{\sigma}}
$$
and $\theta_7<1$ implies $\frac{3}{\sigma}-\frac{3}{p''}<1+{s}$. We obtain the bound
$$
K_1^{22}\le\varrho^*\|v_n\|_{W_{\sigma}^{2+{s},1+{s}/2}(\Omega^t)}\bigg(\intop_0^t \|v_{n+1}\|_{W_{\sigma}^{2+{s}}(\Omega)}^{{\sigma} \theta_7}|v_{n+1}|_{2,\Omega}^{{\sigma}(1-\theta_7)} dt'\bigg)^{1/{\sigma}}\equiv K_1^{23}.
$$
Restrictions $\frac{5}{\sigma}-\frac{3}{p'}<2$, $\frac{3}{\sigma}-\frac{3}{p''}<1+{s}$ imply $\frac{8}{{\sigma}}-\frac{3}{ p}<3+{s}$. Since $\frac{1}{p}<\frac{1}{3}+\frac{1}{\sigma}$, we obtain $\frac{5}{\sigma}<4+{s}$. Hence, there is no restriction.

Hence, by the Young inequality, we get
\begin{equation}\eqal{
K_1^2&\le K_1^{23}\le\varepsilon\|v_{n+1}\|_{W_{\sigma}^{2+{s},1+{s}/2}(\Omega^t)}\cr
&\quad+ (c(1/\varepsilon)\varrho^*\|v_n\|_{W_{\sigma}^{2+{s},1+{s}/2}(\Omega^t)} )^{1/1-\theta_7}\cdot|v_{n+1}|_{2,{\sigma},\Omega^t}.\cr}
\label{6.27}
\end{equation}
Finally,
$$
K_1^3\le\varrho^*\sup_t\|v_n\|_{L_\infty(\Omega)}\bigg(\intop_0^t |\partial_t^{{s}/2}\nabla v_{n+1}|_{p,\Omega}^{\sigma} dt'\bigg)^{1/{\sigma}}\equiv K_1^{31},
$$
where $\frac{1}{3}+\frac{1}{\sigma}>\frac{1}{p}$.

Using the interpolation
$$
|\partial_t^{{s}/2}\nabla v_{n+1}|_{p,{\sigma},\Omega^t}\le c\|v_{n+1}\|_{W_{\sigma}^{2+{s},1+{s}/2}(\Omega^t)}^{\theta_8} |v_{n+1}|_{2,\Omega^t}^{1-\theta_8}
$$
with $\theta_8$ satisfying the equation
$$
\frac{3}{p}+\frac{2}{\sigma}-{s}-1=(1-\theta_8)\frac{5}{2}+\theta_8\bigg(\frac{5}{\sigma}-(2+{s})\bigg),
$$
so
$$
\theta_8=\frac{\frac{5}{2}+1+{s}-\frac{3}{p}-\frac{2}{\sigma}}{2+{s}+\frac{5}{2}-\frac{5}{\sigma}}
$$
and $\theta_8<1$ implies $\frac{3}{\sigma}-\frac{3}{p}<1$, we obtain
$$
K_1^{31}\le c\varrho^*\sup_t\|v_n\|_{L_\infty(\Omega)} \|v_{n+1}\|_{W_{\sigma}^{2+{s},1+{s}/2}(\Omega^t)}^{\theta_8} |v_{n+1}|_{2,\Omega^t}^{1-\theta_8}\equiv K_1^{32}.
$$
Collecting all restrictions
$$
\frac{3}{\sigma}<1+\frac{3}{p}<2+\frac{3}{\sigma}\quad {\rm and}\quad \frac{3}{\sigma}<2+{s}-\frac{2}{\sigma},\quad {\rm so}\quad \frac{5}{\sigma}<2+{s}
$$
and applying the Young inequality, we obtain
\begin{equation}\eqal{
K_1^3&\le K_1^{32}\le\varepsilon\|v_{n+1}\|_{W_{\sigma}^{2+{s},1+{s}/2}(\Omega^t)} \cr
&\quad+(c(1/\varepsilon)\varrho^*\|v_n\|_{V_{\sigma}^{2+{s}}(\Omega^t)})^{1/(1-\theta_8)} \cdot|v_{n+1}|_{2,\Omega^t}.\cr}
\label{6.28}
\end{equation}
From (\ref{6.25})--(\ref{6.28}) we obtain
\begin{equation}\eqal{
K_1&\le\varepsilon\|v_{n+1}\|_{W_{\sigma}^{2+{s},1+{s}/2}(\Omega^t)}\cr
&\quad+(c(1/\varepsilon)\|\varrho_n\|_{W_{r,\infty}^{1,1}(\Omega^t)} \|v_n\|_{V_{\sigma}^{2+{s}}(\Omega^t)})^\frac{1}{ 1-\theta_0} |v_{n+1}|_{2,{\sigma} ,\Omega^t},\cr}
\label{6.29}
\end{equation}
where
$$
\frac{1}{1-\theta_0}=\max_{i\in\{6,7,8\}}\frac{1}{1-\theta_i}.
$$
In the end,
$$
K_3\le c\bigg(\intop_0^t(|f\partial_t^{{s}/2}\varrho_n|_{p,\Omega}^{\sigma}+ |\varrho_n\partial_t^{{s}/2}f|_{p,\Omega}^{\sigma})dt'\bigg)^{1/{\sigma} }\equiv K_3^1,
$$
where $1/3+1/{\sigma}=1/p$. By the H\"older inequality we have
$$
K_3^1\le\bigg(\intop_0^t|f|_{p\lambda_1,\Omega}^{\sigma} |\partial_t^{{s}/2}\varrho_n|_{p\lambda_2,\Omega}^{\sigma} dt'\bigg)^{1/{\sigma}}+\varrho^* \bigg(\intop_0^t|\partial_t^{{s}/2}f|_{p,\Omega}^{\sigma} dt'\bigg)^{1/{\sigma}}\equiv K_3^2,
$$
where $1/\lambda_1+1/\lambda_2=1$. Setting $p\lambda_2=r$ we get $p\lambda_1=\frac{pr}{r-p}$. Then we obtain
\begin{equation}
K_3\le K_3^1\le K_3^2\le\sup_t|\partial_t^{{s}/2}\varrho_n|_{r,\Omega}t^{1/{\sigma}}|f|_{\frac{pr}{ r-p},\infty,\Omega^t}+c\varrho^*|\partial_t^{{s}/2}f|_{p,{\sigma},\Omega^t}.
\label{6.30}
\end{equation}
Using estimates (\ref{6.24}), (\ref{6.29}) and (\ref{6.30}) in the r.h.s. of $J_2^2$ yields
\begin{equation}\eqal{
J_2^2&\le\varepsilon\|v_{n+1}\|_{W_{\sigma}^{2+{s},1+{s}/2}(\Omega^t)}\cr
&\quad+[c(1/\varepsilon)(\|\varrho_n\|_{W_{r,\infty}^{1,1}(\Omega^t)} \|v_n\|_{V_{\sigma} ^{2+{s}}(\Omega^t)}+1)]^{1/(1-\theta_0)}|v_{n+1}|_{2,{\sigma} ,\Omega^t}\cr
&\quad+ct^{1/{\sigma}}\|\varrho_n\|_{W_{r,\infty}^{1,1}(\Omega^t)}|f|_{{\frac{pr}{r-p}},\infty,\Omega^t}
+c\varrho^*|\partial_t^{{s}/2}f|_{p,{\sigma},\Omega^t},\cr
&{\rm where}\  \frac{1}{1-\theta_0}=\max_{i\in\{6,7,8\}}\frac{1}{1-\theta_i}.\cr}
\label{6.31}
\end{equation}
Now we are ready to collect all estimates for $\|p_{n+1}\|_{W_{\sigma}^{{s},{s}/2}(\Omega^t)}$.

From (\ref{6.4}),
\begin{equation}
\|p_{n+1}\|_{W_{\sigma}^{{s},{s}/2}(\Omega^t)}\le \|J_1\|_{W_{\sigma}^{{s},{s}/2}(\Omega^t)}+\|J_2\|_{W_{\sigma}^{{s},{s}/2}(\Omega^t)}.
\label{6.32}
\end{equation}
Next, (\ref{6.6}) and (\ref{6.14}) imply
\begin{equation}\eqal{
&\|J_1\|_{W_{\sigma}^{{s},{s}/2}(\Omega^t)}\le \|J_{11}\|_{W_{\sigma}^{{s},{s}/2}(\Omega^t)}+ \|J_{12}\|_{W_{\sigma}^{{s},{s}/2}(\Omega^t)}\cr
&\le c(\varepsilon t^{\bar{a}}\|\varrho_n\|_{W_{r,\infty}^{1,1}(\Omega^t)}+\varrho^*+ \|\varrho_1\|_{W_r^{1,1}(S_2^t(-a))})\cdot\cr
&\quad\cdot\sum_{i=1}^2\|d_i\|_{W_{\sigma}^{2+{s}-1/{\sigma},1+{s}/2-1/2{\sigma}}(S_2^t(a_i))}\cr
&\quad+\varepsilon\|v_{n+1}\|_{W_{\sigma}^{2+{s},1+{s}/2}(\Omega^t)}+ (c(1/\varepsilon)\|\varrho_n\|_{W_{r,\infty}^{1,1}(\Omega^t)})^{1/(1-\theta_*)} |v_{n+1}|_{2,\Omega^t},\cr}
\label{6.33}
\end{equation}
where
$$
\frac{1}{1-\theta_*}=\max_{i\in\{1,2,3\}}\frac{1}{1-\theta_i}.
$$
Finally,
\begin{equation}
\|J_2\|_{W_{\sigma}^{{s},{s}/2}(\Omega^t)}\le J_2^1+J_2^2,
\label{6.34}
\end{equation}
where (\ref{6.22}) implies
\begin{equation} \eqal{
J_2^1\le\varepsilon\|v_{n+1}\|_{V_{\sigma}^{2+{s}}(\Omega^t)} \cr +[c(1/\varepsilon) (\varrho^*\|v_n\|_{V_{\sigma}^{2+{s}}(\Omega^t)}+1)]^{1/(1-\theta_4)} |v_{n+1}|_{2,{\sigma},\Omega^t}+\varrho^*|f|_{{\sigma},\Omega^t}
\label{6.35} \cr}
\end{equation}
and (\ref{6.31}) yields
\begin{equation}\eqal{
J_2^2&\le\varepsilon\|v_{n+1}\|_{V_{\sigma}^{2+{s}}(\Omega^t)}+(c(1/\varepsilon) (\|\varrho_n\|_{W_{r,\infty}^{1,1}(\Omega^t)}\|v_n\|_{V_{\sigma} ^{2+{s}}(\Omega^t)}\cr
&\quad +1))^{1/(1-\theta_0)}|v_{n+1}|_{2,{\sigma},\Omega^t}+ct^{1/{\sigma}} \|\varrho_n\|_{W_{r,\infty}^{1,1}(\Omega^t)}|f|_{\frac{pr}{r-p},\infty,\Omega^t}\cr
&\quad+c\varrho^*|\partial_t^{{s}/2}f|_{p,{\sigma},\Omega^t},\cr}
\label{6.36}
\end{equation}
where $\frac{1}{1-\theta_0}=\max_{i\in\{4,6,7,8\}}\frac{1}{1-\theta_i}$ and $\frac{1}{3}+\frac{1}{\sigma}=\frac{1}{p}$.

Estimates (\ref{6.32})--(\ref{6.36}) imply (\ref{6.3}) for $\underline{\theta}=\max_{i\in \{1,\ldots,8\}} \theta_i.$ This ends the proof.

\end{proof}


Using the imbeddings from Remark \ref{r7.1} we obtain

\begin{lemma}\label{l7.2}
Assume that $f\in W_{\sigma}^{{s},{s}/2}(\Omega^t)$, $v(0)\in W_{\sigma}^{2+{s}-2/{\sigma}}(\Omega)$,\break $\nabla\varphi\in W_2^{1,1}(\Omega^t)$, $\varrho_n\in W_{r,\infty}^{1,1}(\Omega^t)$, $v_n\in V_{\sigma}^{2+{s}}(\Omega^t)$, $\varrho_1\in W_r^{1,1}(S_2^t(-a))$,\break $d\in W_{\sigma} ^{2+{s}-\frac{1}{{\sigma}},1+\frac{{s}}{2}-\frac{1}{2{\sigma}}}(S_2^t)$. Let $\bar{a}>0$ and $5/r<{s}$, $3/{\sigma}<{s}$, $r>{\sigma}$.\\
Then for solutions to problem (\ref{3.3}) holds
\begin{equation}\eqal{
&\|v_{n+1}\|_{V_{\sigma}^{2+{s}}(\Omega^t)}+\|\nabla p_{n+1}\|_{{s},{\sigma},\Omega^t}\cr
&\le\phi(\|\varrho_1\|_{W_r^{1,1}(S_2^t(-a))} \|d_1\|_{W_{\sigma}^{2+{s}+\frac{1}{\sigma},1+\frac{{s}}{2}-\frac{1}{2{\sigma}}}(S_2^t(-a))})[\phi(t^{\bar{a}} \|\varrho_n\|_{W_{r,\infty}^{1,1}(\Omega^t)},\cr
&\quad t^{\bar{a}}\|v_n\|_{V_{\sigma}^{2+{s}}(\Omega^t)}) (\|\nabla\varphi\|_{W_2^{2,1}(\Omega^t)}+\|f\|_{W_{\sigma}^{{s},{s}/2}(\Omega^t)}\cr
&\quad+\|d\|_{W_{\sigma}^{2+{s}-\frac{1}{{\sigma}},1+\frac{{s}}{2}-\frac{1}{2{\sigma}}}(S_2^t)}+|v(0)|_{2,\Omega})\cr
&\quad+\phi(\|\varrho_1\|_{W_r^{1,1}(S_2^t(-a))} \|d\|_{W_{\sigma}^{2+{s}-\frac{1}{{\sigma}},1+\frac{{s}}{2}-\frac{1}{2{\sigma}}}(S_2^t)}+ \|f\|_{W_{\sigma}^{{s},{s}/2}(\Omega^t)}\cr
&\quad+\|\nabla\varphi\|_{W_2^{2,1}(\Omega^t)} \|v(0)\|_{W_{\sigma}^{2+{s}-2/{\sigma}}(\Omega))}).\cr}
\label{7.1}
\end{equation}
\end{lemma}

\begin{proof}
In view of Remark \ref{r7.1} inequality (\ref{5.1}) implies
\begin{equation}\eqal{
\|v_{n+1}\|_{V(\Omega^t)}&\le\phi(\|d_1\|_{W_{\sigma}^{2+{s}-1/{\sigma},1-{s}/2{\sigma}}(S_2^t(-a))}, \|\varrho_1\|_{W_r^{1,1}(S_2^t(-a))})\cdot\cr
&\quad\cdot[t^{\bar{a}}\|v_n\|_{V_{\sigma}^{2+{s}}(\Omega^t)} \|\nabla\varphi\|_{W_2^{2,1}(\Omega^t)}+\|f\|_{W_{\sigma}^{{s},{s}/2}(\Omega^t)}\cr
&\quad+|v(0)|_{2,\Omega}+\|\nabla\varphi\|_{W_2^{2,1}(\Omega^t)}].\cr}
\label{7.2}
\end{equation}
Similarly, Remark \ref{r7.1} and (\ref{6.3}) imply
\begin{equation}\eqal{
&\|p_{n+1}\|_{W_{\sigma}^{{s},{s}/2}(\Omega^t)}\le\varepsilon \|v_{n+1}\|_{W_{\sigma}^{2+{s},1+{s}/2}(\Omega^t)}\cr
&\quad+\phi(t^{\bar{a}}\|\varrho_n\|_{W_{r,\infty}^{1,1}(\Omega^t)},t^{\bar{a}} \|v_n\|_{V_{\sigma}^{2+{s}}(\Omega^t)})\cdot[|v_{n+1}|_{2,\infty,\Omega^t}\cr
&\quad+\|f\|_{W_{\sigma}^{{s},{s}/2}(\Omega^t)}+ \|d\|_{W_{\sigma}^{2+{s}-\frac{1}{{\sigma}},1+\frac{{s}}{2}-\frac{1}{2{\sigma}}}(S_2^t)}]\cr
&\quad+\phi(\|\varrho_1\|_{W_r^{1,1}(S_2^t(-a))}, \|d\|_{W_{\sigma}^{2+{s}-\frac{1}{{\sigma}},1+\frac{{s}}{2}-\frac{1}{2{\sigma}}}(S_2^t)}, \|f\|_{W_{\sigma}^{{s},{s}/2}(\Omega^t)}).\cr}
\label{7.3}
\end{equation}
Combining (\ref{7.2}), (\ref{7.3}) and (\ref{5.5}) yields (\ref{7.1}). This ends the proof.
\end{proof}

\begin{remark}\label{r7.3}
For solutions to (\ref{3.4}) we have the estimate
\begin{equation}
\|\nabla\varphi\|_{W_2^{1,1}(\Omega^t)}\le c\|d\|_{W_2^{3/2,3/4}(S_2^t)}.
\label{7.4}
\end{equation}
\end{remark}

\begin{corollary}\label{c7.4}
Let the assumptions of Lemma \ref{l7.2} be satisfied, then there exists a solution to problem (\ref{3.3}) such that $v_{n+1}\in V_{\sigma} ^{2+{s}}(\Omega^t)$, $\nabla p_{n+1}\in W_{\sigma}^{{s},{s}/2}(\Omega^t)$ and the estimate (\ref{7.1}) holds.\\
Since $W_{\sigma}^{2+{s},1+{s}/2}(\Omega^t)\subset C^\alpha(\Omega^t)$ for $\alpha$ such that $\frac{5}{\sigma}+\alpha<2+{s}$. Hence $v_{n+1}\in C^\alpha(\Omega^t)$. Since $v_{3;n+1}\in C^\alpha(\Omega^t)$ and $v_{3;n+1}|_{S_2}=d>0$ we have that $v_{3; n+1}$ in a neighborhood of $S_2$ is positive for a sufficiently small time.
\end{corollary}

\section{A lower bound for $v_{3;n+1}$}\label{s8}

 To prove a lower bound for $v_{3;n+1}$ we consider the problem
\begin{equation}\eqal{
&\varrho_n v_{3;n+1,t}+\varrho_n v_n\cdot\nabla v_{3;n+1}-\nu\Delta v_{3;n+1}+q_{n+1}=\varrho_n f_3,\cr
&v_{3;n+1}|_{x_3=-a}=-d_1,\ \ v_{3;n+1}|_{x_3=a}=d_2,\  d_i >0, \ i=1,2,\cr
&v_{3;n+1}|_{t=0}=v_{3;n}(0),\cr}
\label{8.1}
\end{equation}
where $q_{n+1}=p_{n+1,x_3}$.

We need also problem (\ref{3.2}) for $\varrho_n$,
$$ \eqal{
&\varrho_{n,t}+v_n\cdot\nabla\varrho_n=0\quad &{\rm in}\ \ \Omega^T,\cr
&\divv v_n=0\quad &{\rm in}\ \ \Omega^T,\cr
&\varrho_n|_{S_1^T(-a)}=\varrho_1,\cr
&v_n\cdot\bar n|_{S_1^T(-a)}=-d_1,\ \ v_n\cdot\bar n|_{S_2^T(a)}=d_2,\ \ d_i>0,\ \ i=1,2,\cr
&\varrho_n|_{t=0}=\varrho(0).\cr}
$$

\begin{lemma}\label{l8.1}
Assume that $\bar d_0\ge v_3(0)\ge d_0$, where $d_0$, $\bar d_0$ are positive constants and assume that $d_i\ge d_\infty$, $i=1,2$, where $d_\infty$ is also a positive constant and $f_3\in L_1(0,t;L_\infty(\Omega)).$
Let assumptions of Lemma~\ref{l7.2} hold, then $p_{n+1,x_3}\in L_1(0,t;L_\infty(\Omega))$ and there exists a positive number
$$
d_*(n+1)=\phi\bigg(\exp\bigg[-\frac{1}{\varrho_*}(|q_{n+1}|_{\infty,1,\Omega^t}+ |f_3|_{\infty,1,\Omega^t})\bigg],\frac{d_\infty d_0}{ 3\bar d_0+d_\infty}\bigg)
$$
such that
\begin{equation}
v_{3;n+1}\ge d_*(n+1),
\label{8.3}
\end{equation}
where $q_{n+1}=p_{n+1,x_3}$.
\end{lemma}

\begin{proof}
In this proof we drop indices $n$ for simplicity.
We multiply $(\ref{8.1})_1$ by $\frac{v_3}{|v_3|^{p+1}}$ and integrate over $\Omega$. Then we have
\begin{equation}\eqal{
&\intop_\Omega\varrho v_{3,t}\frac{v_3}{|v_3|^{p+1}}dx+\intop_\Omega\varrho v\cdot\nabla v_3\frac{v_3}{|v_3|^{p+1}}dx\cr
&\quad-\nu\intop_\Omega\Delta v_3\frac{v_3}{|v_3|^{p+1}}dx+\intop_\Omega q\frac{v_3}{|v_3|^{p+1}}dx=\intop_\Omega\varrho f_3\frac{v_3}{|v_3|^{p+1}}dx.\cr}
\label{8.4}
\end{equation}
Now, we examine the particular terms in (\ref{8.4}),
$$\eqal{
&v_{3,t}\frac{v_3}{|v_3|^{p+1}}=\frac{1}{2}\frac{\partial_tv_3^2}{|v_3|^{p+1}}= \partial_t|v_3|\frac{1}{|v_3|^p}=\frac{1}{-p+1}\partial_t|v_3|^{-p+1},\cr
&\nabla v_3\frac{v_3}{|v_3|^{p+1}}=\frac{1}{-p+1}\nabla|v_3|^{-p+1}.\cr}
$$
Therefore, the sum of the first two terms in (\ref{8.4}) takes the form
$$
I_1=\frac{1}{-p+1}\intop_\Omega(\varrho\partial_t|v_3|^{-p+1}+\varrho v\cdot\nabla|v_3|^{-p+1})dx.
$$
Using the equation of continuity $(\ref{1.1})_3$ in the form
$$
\varrho_t+\divv(\varrho v)=0
$$
in $I_1$ yields
$$
I_1=\frac{1}{-p+1}\frac{d}{dt}\intop_\Omega\varrho|v_3|^{-p+1}dx+\frac{1}{-p+1} \intop_\Omega\divv(\varrho v|v_3|^{-p+1})dx.
$$
Since $v\cdot\bar n|_{S_1}=0$ the second term in $I_1$ equals
$$\eqal{
&\frac{1}{-p+1}\intop_{S_2(-a)}\varrho v\cdot\bar n|v_3|^{-p+1}dS_2+\frac{1}{-p+1} \intop_{S_2(a)}\varrho v\cdot\bar n|v_3|^{-p+1}dS_2\cr
&=-\frac{1}{-p+1}\intop_{S_2(-a)}\varrho_1d_1|v_3|^{-p+1}dS_2+\frac{1}{-p+1} \intop_{S_2(a)}\varrho d_2|v_3|^{-p+1}dS_2,\cr}
$$
where $d_i\ge d_\infty$, $i=1,2$.

The third term on the l.h.s. of (\ref{8.4}) takes the form
$$\eqal{
&-\nu\intop_\Omega\divv(\nabla v_3v_3|v_3|^{-p-1})dx+\nu\intop_\Omega|\nabla v_3|^2|v_3|^{-p-1}dx\cr
&\quad+\nu\intop_\Omega v_3\nabla v_3\cdot\nabla|v_3|^{-p-1}dx\equiv I_2+I_3+I_4.\cr}
$$
Integral $I_2$ equals,
$$\eqal{
I_2&=-\nu\intop_\Omega\divv(\nabla|v_3|\,|v_3|^{-p})dx=-\frac{\nu}{-p+1} \intop_\Omega\divv(\nabla|v_3|^{-p+1})dx\cr
&=-\frac{\nu}{-p+1}\intop_S\bar n\cdot\nabla|v_3|^{-p+1}dS=-\nu\intop_S|v_3|^{-p} n\cdot\nabla|v_3|dS\cr
&=\nu\intop_{S_2(-a)}d_1^{-p}v_{3,x_3}dS_2-\nu\intop_{S_2(a)}d_2^{-p}v_{3,x_3}dS_2\cr
&\quad-\nu\intop_{S_1}|v_3|^{-p}\bar n\cdot\nabla|v_3|dS_1\equiv I_2^1+I_2^2+I_2^3,\cr}
$$
where in the $n+1$-step $v_{n+1}$ is positive and differentiable in a neighborhood of $S_2$. Using that $v$ is divergence free,
$$\eqal{
I_2^1&=-\nu\intop_{S_2(-a)}d_1^{-p}v_{\alpha,x_\alpha}dS_2=-\nu\intop_{S_2(-a)} (d_1^{-p}v_\alpha)_{,x_\alpha}dS_2\cr
&\quad-\nu p\intop_{S_2(-a)}d_1^{-p-1}d_{1,x_\alpha}v_\alpha dS_2=-\nu\intop_{\partial S_2(-a)}d_1^{-p}v_\alpha\cdot n_\alpha|_{S_1}dL_1\cr
&\quad-\nu p\intop_{S_2(-a)}d_1^{-p-1}d_{1,x_\alpha}v_\alpha dS_2,\cr}
$$
where the first integral vanishes because $L_1\subset S_1$ and $v\cdot\bar n|_{S_1}=0$.

Similarly,
$$
I_2^2=\nu p\intop_{S_2(a)}d_2^{-p-1}d_{2,x_\alpha}v_\alpha dS_2.
$$
To examine $I_2^3$ we recall that condition $(\ref{1.1})_5$ for $\bar\tau_\alpha=\bar e_3$ has the form on $S_1$
$$
\nu \bar n \cdot \nabla v_3 +\gamma v_3=0\quad{\rm so}\quad \nu\bar n\cdot\nabla|v_3|+\gamma|v_3|=0.
$$
Then
$$
I_2^3=\gamma\intop_{S_1}|v_3|^{-p+1}dS_1.
$$
Summarizing,
$$
I_2=-\nu p\intop_{S_2(-a)}d_1^{-p-1}d_{1,x_\alpha}v_\alpha dS_2+\nu p\intop_{S_2(a)}d_2^{-p-1}d_{2,x_\alpha}v_\alpha dS_2+\gamma\intop_{S_1}|v_3|^{-p+1}dS_1.
$$
Next, we calculate
$$\eqal{
I_3&=\nu\intop_\Omega|\nabla v_3|^2|v_3|^{-p-1}dx=\nu\intop_\Omega|v_3^{-\frac{p}{2}-\frac{1}{2}}\nabla v_3|^2dx\cr
&=\frac{4\nu}{(-p+1)^2}\intop_\Omega|\nabla v_3^{-p/2+1/2}|^2dx\cr}
$$
and
$$
I_4=-(p+1)\nu\intop_\Omega|\nabla v_3|^2v_3^{-p-1}dx=-\frac{4\nu(p+1)}{(-p+1)^2}\intop_\Omega|\nabla v_3^{-\frac{p}{ 2}+\frac{1}{2}}|^2dx.
$$
Hence,
$$
I_3+I_4=-\frac{4\nu p}{(-p+1)^2}\intop_\Omega|\nabla v_3^{-\frac{p}{2}+\frac{1}{2}}|^2dx.
$$
In view of Corollary \ref{c7.4} quantities $I_3$ and $I_4$ are well defined.

Using the above results in (\ref{8.4}) yields
\begin{equation}\eqal{
&\frac{1}{-p+1}\frac{d}{dt}\intop_\Omega\varrho |v_3|^{-p+1}dx-\frac{1}{-p+1}\intop_{S_2(-a)}\varrho_1d_1^{-p+2}dS_2\cr
&\quad+\frac{1}{-p+1}\intop_{S_2(a)}\varrho d_2^{-p+2}dS_2-\nu p\intop_{S_2(-a)}d_1^{-p-1}d_{1,x_\alpha}v_\alpha dS_2\cr
&\quad+\nu p\intop_{S_2(a)}d_2^{-p-1}d_{2,x_\alpha}v_\alpha dS_2+\gamma\intop_{S_1}|v_3|^{-p+1}dS_1\cr
&\quad-\frac{4\nu p}{(-p+1)^2}\intop_\Omega|\nabla v_3^{-p/2+1/2}|^2dx+\intop_\Omega q|v_3|^{-p}dx=\intop_\Omega f_3|v_3|^{-p}dx.\cr}
\label{8.5}
\end{equation}
In view of the assumptions of this lemma the last term on the l.h.s. of (\ref{8.5}) is bounded by
\begin{equation}
|q|_{\infty,\Omega}\frac{1}{\varrho_*d_*}\intop_\Omega\varrho|v_3|^{-p+1}dx
\label{8.6}
\end{equation}
and the r.h.s. term by
\begin{equation}
|f_3|_{\infty,\Omega}\frac{1}{\varrho_*d_*}\intop_\Omega\varrho|v_3|^{-p+1}dx,
\label{8.7}
\end{equation}
where $d_*=\min_\Omega v_3$. The existence of this quantity has not been proved yet. It will be found at the end of this proof.

Introducing the notation
\begin{equation}
X^p=\intop_\Omega\varrho|v_3|^{-p+1}dx,
\label{8.8}
\end{equation}
we multiply (\ref{8.5}) by $-p+1$ and exploit estimates (\ref{8.6}) and (\ref{8.7}). Then we obtain
\begin{equation}\eqal{
&\frac{d}{dt}X^p-\frac{4\nu p}{-p+1}\intop_\Omega|\nabla v_3^{-p/2+1/2}|^2dx\cr
&\le(|q|_{\infty,\Omega}+|f_3|_{\infty,\Omega})(p-1)\frac{1}{\varrho_*d_*}X^p\cr
&\quad+\intop_{S_2(-a)}\varrho_1 d_1^{-p+2}dS_2-\intop_{S_2(a)}\varrho d_2^{-p+2}dS_2\cr
&\quad-\nu p(p-1)\intop_{S_2(-a)}d_1^{-p-1}d_{1,x_\alpha}v_\alpha dS_2+\nu p(p-1)\intop_{S_2(a)}d_2^{-p-1}d_{2,x_\alpha}v_\alpha dS_2.\cr}
\label{8.9}
\end{equation}
Let
\begin{equation}
\alpha(t)=(|q(t)|_{\infty,\Omega}+|f_3(t)|_{\infty,\Omega})\frac{p-1}{d_*\varrho_*}.
\label{8.10}
\end{equation}
Then (\ref{8.9}) implies
\begin{equation}\eqal{
&\frac{d}{dt}\bigg(X^p\exp\bigg(-\intop_0^t\alpha(t')dt'\bigg)\bigg)\le\bigg[ \intop_{S_2(-a)}\varrho_1d_1^{-p+2}dS_2\cr
&\quad+\nu p(p-1)\bigg(\!\intop_{S_2(-a)}d_1^{-p-1}|d_{1,x'}|\,|v'|dS_2 \cr
&\quad+\intop_{S_2(a)} d_2^{-p-1}|d_{2,x'}|\,|v'|dS_2\bigg)\bigg]\cdot\exp\bigg(-\intop_0^t \alpha(t')dt'\bigg),\cr}
\label{8.11}
\end{equation}
where $d_{x'}=(d_{,x_1},d_{,x_2})$, $v'=(v_1,v_2)$.

Integrating (\ref{8.11}) with respect to time yields
\begin{equation}\eqal{
X^p&\le\exp\bigg(\intop_0^t\alpha(t')dt'\bigg)\intop_0^t\bigg[\intop_{S_2(-a)} \varrho_1d_1^{-p+2}dS_2\cr
&\quad+\nu p(p-1)\bigg(\intop_{S_2(-a)}d_1^{-p-1}|d_{x'}|\,|v'|dS_2 \cr &\quad +\intop_{S_2(a)} d_2^{-p-1}|d_{2,x'}|\,|v'|dS_2\bigg)\bigg]\cdot\exp\bigg(-\intop_0^{t'}\alpha(t'')dt''\bigg)dt'\cr  &\quad +\exp\bigg(\intop_0^t \alpha(t')dt'\bigg)X^p(0).\cr}
\label{8.12}
\end{equation}
Hence,
\begin{equation}\eqal{
X&\le\exp\bigg(\frac{1}{p}\intop_0^t\alpha(t')dt'\bigg)\frac{1}{ d_\infty}\bigg[\bigg(\intop_{S_2^t(-a)}\varrho_1d_1^2dS_2dt'\bigg)^{1/p}\cr
&\quad+(\nu p(p-1))^{1/p}\bigg(\intop_{S_2^t(-a)}d_1^{-1}|d_{1,x'}|\,|v'|dS_2dt'\bigg)^{1/p}\cr
&\quad+\bigg(\intop_{S_2^t(a)}d_2^{-1}|d_{2,x'}|\,|v'|dS_2dt'\bigg)^{1/p}\bigg]+ \exp\bigg(\frac{1}{p}\intop_0^t\alpha(t')dt'\bigg)X(0).\cr}
\label{8.13}
\end{equation}
Passing with $p\to\infty$ implies
$$
\sup_{x\in\Omega}\bigg|\frac{1}{v_3(x)}\bigg|\le\exp\bigg[\frac{1}{d_*\varrho_*}(|q|_{\infty,1,\Omega^t}+ |f_3|_{\infty,1,\Omega^t})\bigg]\bigg(\frac{3}{d_\infty}+\bigg|\frac{1}{v(0)}\bigg|_{\infty,\Omega}\bigg).
$$
Hence, we have
\begin{equation}\eqal{
&v_3\ge\inf_{x\in\Omega} v_3(x)\cr
&\ge\exp\bigg(-\bigg[\frac{1}{d_*\varrho_*}(|q|_{\infty,1,\Omega^t}+ |f_3|_{\infty,1,\Omega^t})\bigg]\bigg)\frac{d_\infty\inf_{x\in\Omega}|v_3(0)|}{3\inf_{x\in\Omega}|v_3(0)|+d_\infty}\cr
&\ge\exp\bigg(-\bigg[\frac{1}{d_*\varrho_*}(|q|_{\infty,1,\Omega^t}+ |f_3|_{\infty,1,\Omega^t})\bigg]\bigg)\frac{d_\infty d_0}{ 3\bar d_0+d_\infty}\cr
&\equiv d_*.\cr}
\label{8.14}
\end{equation}
Introduce the notation
$$\eqal{
&a(t)=\frac{1}{\varrho_*}(|q|_{\infty,1,\Omega^t}+|f_3|_{\infty,1,\Omega^t}),\cr
&b=\frac{d_\infty d_0}{ 3\bar d_0+d_\infty},\cr}
$$
where $\bar d_0>d_0$.

It is clear that $a(t)$ and $b$ are positive and $a(t)$ can be estimated by
$$
a(t)\le\frac{1}{\varrho_*}t^{1/{\sigma}'}(|q|_{\infty,\varkappa,\Omega^t}+ |f|_{\infty,\varkappa,\Omega^t}),\quad \frac{1}{\varkappa}+\frac{1}{\varkappa'}=1.
$$
Since $\varkappa>1$, $a(t)$ is small for small $t$.

Then (\ref{8.14}) implies the following equation for $d_*$,
\begin{equation}
\exp\bigg(-\frac{1}{d_*}a(t)\bigg)b=d_*.
\label{8.15}
\end{equation}
Therefore
\begin{equation}
\exp(-a(t))=\bigg(\frac{d_*}{b}\bigg)^{d_*}.
\label{8.16}
\end{equation}
The function
$$
h(x)=\bigg(\frac{x}{b}\bigg)^x
$$
equals 1 for $x=0$ and $x=b$.

In the interval $(0,b)$ we have $h(x)<1$. It attains minimum at $x=be^{-b}=b_*$ and $\frac{dh}{dx}<0$ for $x\in(0,b_*)$ and $\frac{dh}{ dx}>0$ for $x\in(b_*,b)$.

At the minimum
$$
h(be^{-b})=e^{-b^2e^{-b}}\equiv h_*.
$$
Since (\ref{8.16}) holds we have the restriction
$$
e^{-b^2e^{-b}}\le e^{-a(t)}.
$$
Hence
\begin{equation}
e^ba(t)\le b^2.
\label{8.17}
\end{equation}
Therefore, (\ref{8.17}) implies that $t$ must be sufficiently small.

Since $h_*<e^{-a}<1$ there exists a point $x_*$ such that
$$
\bigg(\frac{x_*}{b}\bigg)^{x_*}=e^{-a}.
$$
Hence $x_*=d_*$, which is a solution to (\ref{8.15}).

Moreover, $d_*=\phi(e^{-a(t)},b)$. Recalling the index $n,$ this ends the proof.
\end{proof}

\section{Existence of local solution to problem (\ref{1.1}) }\label{s9}
\begin{remark}\label{r9.1}
Let
\begin{equation}\eqal{
\bar d_2&=\|f\|_{W_{\sigma}^{{s},{s}/2}(\Omega^t)}+\sum_{i=1}^2 \|d_i\|_{W_{\sigma}^{2+{s}-1/{\sigma},1+{s}/2}(S_2^t(a_i))}\cr
&\quad+\|\varrho_1\|_{W_r^{1,1}(S_2^t(-a))}+\|v(0)\|_{W_{\sigma}^{2+{s}-2/{\sigma}}(\Omega)}\cr}
\label{9.1}
\end{equation}
and $\bar d_1$ be defined in (\ref{4.21}). Then equations (\ref{4.23}) and (\ref{7.1}) imply
\begin{equation}\eqal{
&\|v_{n+1}\|_{V_{\sigma}^{2+{s}}(\Omega^t)}+\|\nabla p_{n+1}\|_{W_{\sigma}^{{s},{s}/2}(\Omega^t)}\cr
&\le\phi(t^{\bar{a}} \|v_n\|_{V_{\sigma}^{2+{s}}(\Omega^t)},t^{\bar{a}}\|\nabla p_n\|_{W_{\sigma}^{{s},{s}/2}(\Omega^t)},\bar d_1,\bar d_2).\cr}
\label{9.2}
\end{equation}
\end{remark}

\begin{lemma} \label{l9.1}
Let $\bar d_1,$ $\bar d_2$ and $\phi$ be such as in Remark~\ref{r9.1}. Let $M$ be such a number that
$$
\phi(0,\bar d_1,\bar d_2)\le\frac{1}{2}M
$$
and
$$
\|\tilde v\|_{W_{\sigma}^{2+{s},1+{s}/2}(\Omega^t)}\le M,
$$
where $\tilde v$ is an extension of initial data such that $\tilde v|_{t=0}=v(0)$.\\
Then for $t$ sufficiently small
\begin{equation}
\|v_n\|_{V_{\sigma}^{2+{s}}(\Omega^t)}+\|\nabla p_n\|_{W_{\sigma}^{{s},{s}/2}(\Omega^t)}\le M\quad {\rm for\ any}\ \ n\in\N.
\label{9.3}
\end{equation}
\end{lemma}

\begin{proof}
We prove the lemma by the method of successive approximations. To continue the method we have to know that (\ref{8.3}) holds for any $n\in\N$, so
\begin{equation}
v_{3,n}\ge d_*(n).
\label{9.4}
\end{equation}
Let $n$ be fixed and (\ref{9.4}) holds. Introduce the quantity
\begin{equation}
X_n(t)=\|v_n\|_{V_{\sigma}^{2+{s}}(\Omega^t)}+\|\nabla p_n\|_{W_{\sigma}^{{s},{s}/2}(\Omega^t)}.
\label{9.5}
\end{equation}
Then (\ref{4.23}) and (\ref{7.1}) imply the inequality
\begin{equation}
X_{n+1}(t)\le\phi(X_n(t)t^{\bar{a}},\bar d_1,\bar d_2).
\label{9.6}
\end{equation}
Assuming that $X_n(t)\le M$, where $M$ is such that
$$
\phi(0,\bar d_1,\bar d_2)<M,
$$
we can obtain from (\ref{9.6}) for $t$ sufficiently small the estimate
\begin{equation}
X_{n+1}(t)\le M.
\label{9.7}
\end{equation}
Estimate (\ref{9.7}) implies that
$$
|q_{n+1}|_{\infty,1,\Omega^t}\le\|p_{n+1,x_3}\|_{W_{\sigma}^{{s},{s}/2}(\Omega^t))} \le M.
$$
It means that (\ref{9.4}) also holds for
\begin{equation}
v_{3;n+1}\ge d_*(n+1).
\label{9.8}
\end{equation}
Since $X_n(t)\le M$ implies that $X_{n+1}(t)\le M,$ then Lemma \ref{l8.1} gives that $d_*(n)=d_*$ for any $n$.

We also assume that
\begin{equation}
v_3(0)\ge d_0.
\label{9.9}
\end{equation}
Hence, the method of successive approximations can be continued and (\ref{9.3}) holds. This ends the proof.
\end{proof}

In order to show a convergence of the considered sequence $\{v_n\}_{n=1}^\infty$, we define differences
\begin{equation}
V_{n+1}=v_{n+1}-v_n,\quad P_{n+1}=p_{n+1}-p_n,\quad R_n=\varrho_n-\varrho_{n-1}.
\label{9.10}
\end{equation}
which are solutions to following problems
\begin{equation}\eqal{
&\varrho_nV_{n+1,t}+\varrho_nv_n\cdot\nabla V_{n+1}-\divv\T(V_{n+1},P_{n+1})\cr
&=-R_nv_{n,t}-(R_nv_n+\varrho_{n-1}V_n)\cdot\nabla v_n+R_nf,\cr
&\divv V_{n+1}=0,\cr
&\bar n\cdot V_{n+1}=0\quad &{\rm on}\ \ S^T,\cr
&\nu\bar n\cdot\D(V_{n+1})\cdot\bar\tau_\alpha+\gamma V_{n+1}\cdot\bar\tau_\alpha=0,\ \ \alpha=1,2\quad &{\rm on}\ \ S_1^T,\cr
&\bar n\cdot\D(V_{n+1})\cdot\bar\tau_\alpha=0,\ \ \alpha=1,2\quad &{\rm on}\ \ S_2^T,\cr
&V_{n+1}|_{t=0}=0,\cr}
\label{9.11}
\end{equation}
and
\begin{equation}\eqal{
&R_{n,t}+v_n\cdot\nabla R_n=-V_n\cdot\nabla\varrho_{n-1},\cr
&R_n|_{t=0}=0,\ \ R_n|_{S_2(-a)}=0.\cr}
\label{9.12}
\end{equation}

\begin{lemma}\label{l9.2}
Let assumptions of Lemma \ref{l9.1} hold.
Then for $t$ sufficiently small the sequence $\{v_n\}_{n=1}^\infty$ converges.
\end{lemma}

\begin{proof}
We multiply $(\ref{9.11})_1$ by $V_{n+1}$, integrate over $\Omega$ and use boundary conditions and the equation
$$
\varrho_{n,t}+\divv(\varrho_nv_n)=0.
$$
Then we obtain
\begin{equation}\eqal{
&\frac{1}{2}\frac{d}{dt}\intop_\Omega\varrho_nV_{n+1}^2dx+\nu\|V_{n+1}\|_{1,\Omega}^2=\intop_{S_2(-a)} \varrho_1d_1V_{n+1}^2dS_2\cr
&\quad-\intop_{S_2(a)}\varrho_nd_2V_{n+1}^2dS_2-\gamma\intop_{S_1}|V_{n+1}\cdot \bar\tau_\alpha|^2dS_1-\intop_\Omega R_nv_{n,t}V_{n+1}dx\cr
&\quad-\intop_\Omega(R_nv_n+\varrho_{n-1}V_n)\cdot\nabla v_nV_{n+1}dx+\intop_\Omega R_nf\cdot V_{n+1}dx.\cr}
\label{9.13}
\end{equation}
Since $\varrho_n>0$, $d_2>0$, $\gamma>0,$ the second and the third terms on the r.h.s. of (\ref{9.13}) can be dropped. The first term on the r.h.s. is bounded by
$$
|V_{n+1}|_{3,S_2}^2|\varrho_1d_1|_{3,S_2(-a)}\le\varepsilon^{1/6}|\nabla V_{n+1}|_{2,\Omega}^2+c\varepsilon^{-5/6}|\varrho_1d_1|_{3,S_2(-a)}^6|V_{n+1}|_{2,\Omega}^2.
$$
The fourth term on the r.h.s. of (\ref{9.13}) is bounded by
$$
\varepsilon|V_{n+1}|_{6,\Omega}^2+c(1/\varepsilon)|R_n|_{2,\Omega}^2 |v_{n,t}|_{3,\Omega}^2,
$$
the fifth by
$$\eqal{
&\varepsilon|V_{n+1}|_{6,\Omega}^2+c(1/\varepsilon)(|v_n|_{6,\Omega}^2|\nabla v_n|_{6,\Omega}^2|R_n|_{2,\Omega}^2+|\varrho_{n-1}|_{\infty,\Omega}^2|\nabla v_n|_{3,\Omega}^2|V_n|_{2,\Omega}^2).\cr}
$$
Finally, the last by
$$
\varepsilon|V_{n+1}|_{6,\Omega}^2+c(1/\varepsilon)|f|_{3,\Omega}^2|R_n|_{2,\Omega}^2.
$$
Using the above estimates in (\ref{9.13}) yields
\begin{equation}\eqal{
&\frac{d}{dt}|V_{n+1}|_{2,\Omega}^2+\nu\|V_{n+1}\|_{1,\Omega}^2\le c(\varrho^*)^6|d_1|_{3,\Omega}^6|V_{n+1}|_{2,\Omega}^2\cr
&\quad+c(|v_{n,t}|_{3,\Omega}^2+|v_n|_{6,\Omega}^2|\nabla v_n|_{6,\Omega}^2)|R_n|_{2,\Omega}^2\cr
&\quad+c|\varrho_{n-1}|_{\infty,\Omega}^2|\nabla v_n|_{3,\Omega}^2|V_n|_{2,\Omega}^2+c|f|_{3,\Omega}^2|R_n|_{2,\Omega}^2.\cr}
\label{9.14}
\end{equation}
Consider problem (\ref{9.12}). Multiplying $(\ref{9.12})_1$ by $R_n$ and integrating over $\Omega$ gives
\begin{equation}
\frac{1}{2}\frac{d}{dt}|R_n|_{2,\Omega}^2+\frac{1}{2}\intop_\Omega v_n\cdot\nabla R_n^2dx=-\intop V_n\nabla\varrho_{n-1}R_ndx+\frac{1}{2}\intop_{S_2(a)}d_2R_n^2dS_2.
\label{9.15}
\end{equation}
The last term on the l.h.s. of (\ref{9.15}) equals $\frac{1}{2}\intop_{S_2(a)}d_2R_n^2dS_2$ because $\divv v_n=0$, $v_n\cdot\bar n|_{S_1}=0$ and $R_n|_{S_2(-a)}=0$.

\noindent Then we obtain
$$
\frac{d}{dt}|R_n|_{2,\Omega}^2\le|\nabla\varrho_{n-1}|_{3,\Omega}|V_n|_{6,\Omega}|R_n|_{2,\Omega}.
$$
Integrating with respect to time implies
\begin{equation}
|R_n(t)|_{2,\Omega}\le|\nabla\varrho_{n-1}|_{3,\infty,\Omega^t}|V_n|_{6,1,\Omega^t}.
\label{9.16}
\end{equation}
Integrating (\ref{9.14}) with respect to time and using (\ref{9.16}) we obtain
\begin{equation}\eqal{
\|V_{n+1}\|_{V(\Omega^t)}^2&\le c\exp(\varrho_*^6|d_1|_{3,6,S_2^t(-a)}^6)\cdot [(|v_{n,t}|_{3,\infty,\Omega^t}^2\cr
&\quad+|v_n|_{6,\infty,\Omega^t}^2|\nabla v_n|_{6,\infty,\Omega^t}^2)t|\nabla\varrho_{n-1}|_{3,\infty,\Omega^t}^2 |V_n|_{6,1,\Omega^t}^2\cr
&\quad+|\varrho_{n-1}|_{\infty,\Omega^t}^2+|\nabla v_n|_{3,\infty,\Omega^t}^2t|V_n|_{2,\infty,\Omega^t}^2\cr
&\quad+c|f|_{3,\infty,\Omega^t}^2|\nabla\varrho_{n-1}|_{3,\infty,\Omega^t}^2t |V_n|_{6,2,\Omega^t}^2].\cr}
\label{9.17}
\end{equation}
Through imbeddings and estimate (\ref{9.3}) we get
$$\eqal{
\|V_{n+1}\|_{V(\Omega^t)}^2&\le t \exp(\varrho_*^6|d_1|_{3,6,S_2^t(-a)}^6)[M^6+M^4+ M^2|f|_{3,\infty,\Omega^t}^2]\|V_n\|_{V(\Omega^t)}^2.\cr}
$$
Hence, for $t$ sufficiently small the sequence $\{v_n\}_{n=1}^\infty$ converges. This ends the proof.
\end{proof}
\goodbreak

\begin{remark}\label{r9.3}
Lemmas \ref{l9.1} and \ref{l9.2} imply Theorem \ref{t1.1}.
\end{remark}

\end{document}